\documentclass[11pt,leqno]{amsart}
\usepackage{amssymb,amsmath,amsthm,dsfont,fourier}

\usepackage{mathtools}
\usepackage{array}
\usepackage[text={6.5in,8.5in},centering]{geometry}
\usepackage[dvipsnames]{xcolor}
\usepackage{enumitem}
\usepackage{tikz-cd}
\usetikzlibrary{decorations.pathmorphing}
\usepackage[colorlinks=true,citecolor=NavyBlue,linkcolor=NavyBlue,urlcolor=NavyBlue]{hyperref}

\newtheorem{theorem}{Theorem}[section]
\newtheorem{proposition}[theorem]{Proposition}
\newtheorem{lemma}[theorem]{Lemma}
\newtheorem{corollary}[theorem]{Corollary}
\theoremstyle{definition}
\newtheorem{definition}[theorem]{Definition}

\newtheorem{example}[theorem]{Example}
\newtheorem*{theoremA}{Theorem A}
\newtheorem*{theoremB}{Theorem B}
\newtheorem*{theoremC}{Theorem C}
\newtheorem*{theoremD}{Theorem D}

\newcommand{\M}{\mathcal{M}}
\newcommand{\N}{\mathcal{N}}

\newcommand{\E}{\mathcal{E}}
\newcommand{\R}{\mathcal{R}}
\renewcommand{\P}{\mathcal{P}}

\newcommand{\Z}{\mathds{Z}}
\newcommand{\F}{\mathds{F}}
\newcommand{\trnl}{\trianglelefteq}
\newcommand{\Soc}{\operatorname{Soc}}
\newcommand{\MinNorm}{\operatorname{MinNorm}}
\newcommand{\End}{\operatorname{End}}
\newcommand{\Hom}{\operatorname{Hom}}
\newcommand{\Aut}{\operatorname{Aut}}
\newcommand{\Sub}{\operatorname{Sub}}
\newcommand{\md}{\operatorname{md}}
\newcommand{\clM}{\operatorname{cl}_{\M}}
\newcommand{\supp}{\operatorname{supp}}
\newcommand{\elln}{\ell_{\mathrm n}}
\newcommand{\qbinom}[2]{\genfrac[]{0pt}{}{#1}{#2}}
\newcommand{\id}{\text{id}}
\newcommand{\im}{\operatorname{im}}
\newcommand{\ord}{\operatorname{ord}}

\makeatletter

\def\twoheadrightarrowfill@{%
    \arrowfill@\relbar\relbar\twoheadrightarrow
}

\newcommand{\xtwoheadrightarrow}[2][]{%
    \ext@arrow 0359\twoheadrightarrowfill@{#1}{\mkern8mu#2\mkern8mu}%
}

\makeatother

\title[$\M$-subgroups and normal-subgroup lattices]{$\M$-subgroups and normal-subgroup lattices of finite groups}

\author{Aneesh Nag}
\address{[A. Nag] \; Department of Pure Mathematics, University of Calcutta, 35 Ballygunge Circular Road, Kolkata, West Bengal, 700019, India}
\email{aneesh89884@gmail.com}
\urladdr{ORCID: \url{https://orcid.org/0009-0007-3915-5316}}

\author{Pronay Biswas}
\address{[P. Biswas] \; Department of Mathematics, Jadavpur University, 188 Raja Subodh Chandra Mallick Road, Jadavpur, Kolkata, West Bengal 700032, India}
\email{pronayb.math.rs@jadavpuruniversity.in, pronaybiswas1729@gmail.com}
\urladdr{ORCID: \url{https://orcid.org/0009-0008-2710-1748}}

\author{Amartya Goswami}
\address{[A. Goswami] \;Department of Mathematics and Applied Mathematics, University of
Johannesburg, P.O.\ Box 524, Auckland Park 2006, South Africa}
\address{National Institute for Theoretical and Computational Sciences
(NITheCS), South Africa}
\email{agoswami@uj.ac.za}
\urladdr{ORCID: \url{https://orcid.org/0000-0003-4829-0847}}


\begin{document}
\begin{abstract}
Let $G$ be a group. A normal subgroup $E$ is an $\M$-subgroup if
$E\cap N_1\cap N_2\neq1$ whenever $N_1\cap N_2$, $E\cap N_1$, and
$E\cap N_2$ are non-trivial for normal subgroups $N_1,N_2$. We
characterise these subgroups as the essential or relatively uniform
normal subgroups. For finite $G$, an $\M$-subgroup is either trivial,
contains exactly one minimal normal subgroup, or contains the socle.
Consequently, every normal subgroup is an $\M$-subgroup precisely when
the socle has $G$-chief length at most two. We obtain an explicit
$\M$-closure operation and describe the minimal normal subgroups of
finite direct products. For crown-based powers, we use the known
comparability of normal subgroups with the socle to describe the normal
lattice in terms of Boolean or subspace lattices and the quotient above
the socle. We then determine when the $\M$-lattice retains the normal
lattice, classify normal quotients and their automorphism orbits, and
characterise quotient maps compatible with $\M$-closure. Finally, we
compute the maximal length of chains of $\M$-subgroups in terms of the
chief length and socle length of $G$.
\end{abstract}
\maketitle
\tableofcontents
\section{Introduction}
\label{sec:intro}

For a finite group $G$, the normal subgroups form a finite modular lattice
$\N(G)$ whose atoms are the minimal normal subgroups of $G$. Their product
is the socle $\Soc(G)$, which is also the least essential normal subgroup:
it meets every non-trivial normal subgroup, and every essential normal
subgroup contains it. Thus the socle provides a natural point of departure
for studying intersection properties in $\N(G)$. In this paper we identify
a class of normal subgroups defined by such a property and show that, for
finite groups, its structure is governed by the socle and by the geometry
of minimal normal subgroups.

A normal subgroup $E\trnl G$ is called an \emph{$\M$-subgroup} if, for all
normal subgroups $N_1,N_2\trnl G$,
\[
 N_1\cap N_2\neq1,\qquad E\cap N_1\neq1,\qquad E\cap N_2\neq1
 \quad\Longrightarrow\quad E\cap N_1\cap N_2\neq1.
\]
As we prove in Proposition~\ref{prop:pairs}, this condition is equivalent
to the analogous implication for any finite family of normal subgroups.
We write $\N^{\M}(G)$ for the collection of $\M$-subgroups and call $G$
\emph{$\M$-full} when $\N^{\M}(G)=\N(G)$. Our principal questions are to
characterise these subgroups, determine when $G$ is $\M$-full, and describe
the resulting lattice for direct products and crown-based powers.

The intersection condition is motivated by the theory of $\M$-ideals in
Banach spaces, initiated by Alfsen and Effros~\cite{AlfsenEffros} (see also~\cite{HWW}), and
its later analogues for ring ideals and quantales. In particular, ~\cite{BlecherGoswami} introduced
ring $\M$-ideals, while ~\cite{CavigliaGoswamiMesiti}
studied $\mu$-elements of quantales. The dichotomy between essentiality and a
relative uniformity condition in these settings motivates the group case.
Our arguments use the modular normal-subgroup lattice directly, rather
than applying the quantale results.

Let us now discuss our main results.
We first show, for arbitrary groups, that a normal subgroup is an
$\M$-subgroup if and only if it is essential or relatively uniform
(Theorem~\ref{thm:dichotomy}); the latter means that any two non-trivial
normal subgroups contained in it meet non-trivially. For finite groups,
this becomes a classification in terms of minimal normal subgroups.
Write $\MinNorm(G)$ for their set and put
$\supp(N)=\{A\in\MinNorm(G):A\leq N\}$ for $N\trnl G$.

\begin{theoremA}
Let $G$ be finite and let $E\trnl G$.  Then $E$ is an $\M$-subgroup if and
only if
\[
 E=1,\qquad |\supp(E)|=1,\qquad\text{or}\qquad \Soc(G)\leq E.
\]
The $\M$-subgroups are closed under intersections.  Moreover, the least
$\M$-subgroup containing a normal subgroup $N$ is
\[
\clM(N)=
\begin{cases}
N,& |\supp(N)|\leq 1,\\
N\Soc(G),& |\supp(N)|\geq 2.
\end{cases}
\]
Consequently they form a lattice whose join is
$E\mathbin{\vee_{\M}}F=\clM(EF)$.
\end{theoremA}

The classification extends to normally atomic groups
(Theorem~\ref{thm:atomic-classification}). In the finite case it also
provides an explicit closure operator on $\N(G)$, and hence a lattice
structure on $\N^{\M}(G)$ (Proposition~\ref{prop:Mclosure} and
Corollary~\ref{cor:Mlattice}). In particular, the normal subgroups excluded
from $\N^{\M}(G)$ are precisely those containing more than one, but not
all, of the minimal normal subgroups of $G$.

Let $\sigma(G)$ denote the length of the interval
$[1,\Soc(G)]_{\N(G)}$; equivalently, it is the number of factors in a
decomposition of the socle into a direct product of minimal normal
subgroups. The criterion for $\M$-fullness is then numerical.

\begin{theoremB}
Let $G$ be a non-trivial finite group and put $S=\Soc(G)$.  The following are
equivalent:
\begin{enumerate}
\item every normal subgroup of $G$ is an $\M$-subgroup;
\item $\sigma(G)\leq2$;
\item one of the following holds:
      \begin{enumerate}[label=\upshape(\alph*)]
      \item $G$ has one minimal normal subgroup;
      \item $G$ has exactly two minimal normal subgroups;
      \item for some prime $p$ and some irreducible $\F_pG$-module $A$,
            one has $S\cong A\oplus A$ as $G$-modules.
      \end{enumerate}
\end{enumerate}
In alternative \textup{(c)}, if $D=\End_{\F_pG}(A)$, then
$\MinNorm(G)$ is naturally the projective line $\mathds P^1(D)$ and has
$|D|+1$ elements.
\end{theoremB}

Thus the relevant invariant is the \emph{length} of the socle, not the
number of its minimal normal subgroups. For example, $C_p\times C_p$
has $p+1$ minimal normal subgroups, but its socle has length two, so it
is $\M$-full. More generally, the abelian isotypic components of the
socle give rise to finite projective geometries
(Proposition~\ref{prop:socle-geometry}). The classification also has a
lattice-theoretic formulation: for arbitrary $G$, $\M$-fullness is
equivalent to the normal-subgroup lattice having uniform (Goldie)
dimension at most two (Theorem~\ref{thm:Mfull-udim}).

Direct products and crown-based powers make these results concrete. The
minimal normal subgroups of a direct product need not come from the
factors: central elementary abelian subgroups also produce diagonal
minimal normal subgroups. A crown-based power, introduced by Dalla Volta
and Lucchini~\cite{DallaVoltaLucchini}, is an important construction in
the structure and generation theory of finite groups; in particular,
certain quotients $G/R_G(B)$ attached to non-Frattini chief factors are
crown-based powers~\cite{DetomiLucchiniCrowns,LucchiniCrowns}.
Burness, Garonzi and Lucchini~\cite[Section~5]{BurnessGaronziLucchini}
use crown theory to bound the intersection number of a finite group,
which is defined using intersections of maximal subgroups. Our
intersection condition, by contrast, concerns normal subgroups. We
describe the atoms of the corresponding normal-subgroup lattices and
determine $\M$-fullness.

\begin{theoremC}
Let $G=G_1\times\cdots\times G_t$ be finite.  For each prime $p$, set
\[
 V_{i,p}=\{z\in Z(G_i):z^p=1\},
 \qquad V_p=\prod_{i=1}^tV_{i,p}.
\]
The minimal normal subgroups of $G$ are precisely the coordinate copies of
the noncentral minimal normal subgroups of the factors, together with the
one-dimensional $\F_p$-subspaces of the spaces $V_p$.  If $t\geq2$ and all
factors are non-trivial, then $G$ is $\M$-full exactly when $t=2$ and both
factors are monolithic.

If $L$ is a finite monolithic primitive group with socle $A$ and $L_k$ is its
crown-based power, then $\Soc(L_k)=A^k$.  If $A$ is non-abelian, the atoms are
the $k$ coordinate copies of $A$.  If $A$ is abelian, then $A$ is an
irreducible $\F_p[L/A]$-module for a prime $p$; with
$D=\End_{\F_p[L/A]}(A)$, the atoms are parametrized by
$\mathds P^{k-1}(D)$.  In either case every normal subgroup of $L_k$ is an
$\M$-subgroup exactly when $k\leq2$.
\end{theoremC}

The comparability of every normal subgroup of a crown-based power
$L_k$ with its socle $V=A^k$ is known by a lemma of Lucchini,
Moscatiello and Spiga~\cite[Lemma~2.2]{LucchiniMoscatielloSpiga}.
Using this result, we identify the lower interval of $\N(L_k)$ with
a Boolean lattice or a subspace lattice and the upper interval with
$\N(L/A)$ (Theorem~\ref{thm:crown-normal-lattice}). Our subsequent
results concern the information retained by $\N^\M(L_k)$: we give
a rigidity criterion comparing normal-subgroup lattices and
$\M$-lattices (Theorem~\ref{thm:crown-lattice-rigidity}), classify the
quotients by arbitrary normal subgroups contained in $V$ together
with their automorphism orbits (Theorem~\ref{thm:crown-quotients}),
and determine precisely which quotient maps preserve $\M$-closure
(Proposition~\ref{prop:crown-closure-quotients}). The description also
gives enumeration formulas (Corollary~\ref{cor:crown-M-lattice}) and
applies to crowns associated with arbitrary finite groups
(Corollary~\ref{cor:arbitrary-crown}).

Finally, the omission of normal subgroups from $\N^{\M}(G)$ can be
measured by a chain-length invariant. Define the \emph{$\M$-depth}
$\md(G)$ to be the maximum length of a chain from $1$ to $G$ in
$\N^{\M}(G)$, and let $\elln(G)$ be the chief length of $G$.

\begin{theoremD}
Let $G$ be a non-trivial finite group.  If $G$ has one minimal normal subgroup,
then
\[
\md(G)=\elln(G).
\]
If $G$ has at least two minimal normal subgroups, then
\[
\md(G)=\elln(G)-\sigma(G)+2.
\]
In particular, if $G_1,\ldots,G_t$ are non-trivial finite groups and $t\geq2$,
then
\[
\md(G_1\times\cdots\times G_t)
 =\sum_{i=1}^t\bigl(\elln(G_i)-\sigma(G_i)\bigr)+2.
\]
\end{theoremD}

In fact, every maximal chain of $\M$-subgroups has the same length.
When $G$ has at least two minimal normal subgroups, such a chain follows
a chief series below a normal subgroup containing a unique minimal
normal subgroup, crosses the socle in one step, and continues along a
chief series above the socle
(Corollary~\ref{cor:maximum-chain-structure}). Consequently,
$\md(G)=\elln(G)$ holds precisely for $\M$-full finite groups
(Corollary~\ref{cor:depth-Mfull}). We also compute the depth for several
standard families, including finite abelian and nilpotent groups,
dihedral groups, and crown-based powers.
\smallskip

The organization of the paper is as follows. In
Section~\ref{sec:definitions}, we fix notation. In Section~\ref{sec:intersections}, we
prove the basic dichotomy and the classification of $\M$-subgroups. The purpose of 
Section~\ref{sec:noloss} is to treat $\M$-fullness and its relation to the
socle, whereas in Sections~\ref{sec:products} and~\ref{sec:crowns} we address direct
products and crown-based powers. In Section~\ref{sec:abelian}, we specialise
the results to finite abelian groups, and finally, in Section~\ref{sec:depth}, we study
$\M$-depth.

\section{Definitions and notation}
\label{sec:definitions}

Groups are arbitrary unless finiteness is stated explicitly. We begin with the lattice
of normal subgroups and the standard lattice-theoretic vocabulary that we
use for it.

\begin{definition}\label{def:lattice-notation}
For a group $G$, write
\[
\N(G)=\{N\leq G:N\trnl G\}
\]
for the poset of normal subgroups ordered by inclusion, and write $\Sub(G)$
for the poset of all subgroups.  The least element of $\N(G)$ is the trivial
subgroup $1$ and the greatest element is $G$.  For $H,K\trnl G$, their meet
is $H\cap K$ and their join is the subgroup product
$HK=\{hk:h\in H,\ k\in K\}$.  For $H\leq K$ in a poset $P$, the interval
$[H,K]_P$ is $\{X\in P:H\leq X\leq K\}$, and the principal filter generated
by $H$ is $[H,G]_P$ when $G$ is the greatest element.  An \emph{atom} of a
poset with least element is a minimal element strictly larger than the least
one.  A lattice is \emph{modular} if
\[
x\vee(y\wedge z)=(x\vee y)\wedge z\qquad(x\leq z).
\]
A \emph{meet-subsemilattice} is a subset closed under finite nonempty meets.
The length of a finite chain $X_0<X_1<\cdots<X_r$ is $r$; the length of a
finite interval is the maximum length of a chain joining its endpoints.
\end{definition}

The group-theoretic notation is standard~\cite{Isaacs,Robinson}; we record
the conventions that matter for computations, notably the order of factors in
commutators.

\begin{definition}\label{def:group-notation}
For subgroups $A,B\leq G$, our commutator convention is
$[a,b]=aba^{-1}b^{-1}$, and $[A,B]$ denotes the subgroup generated by
these commutators. We write $Z(G)$ for the centre of $G$ and $C_n$
for the cyclic group of order $n$. As usual, an elementary abelian
$p$-group is regarded as a vector space over $\F_p$.
All other elementary group-theoretic terminology is standard.
\end{definition}

The atoms of $\N(G)$ are the minimal normal subgroups, and they are the
central objects of this paper. The following notions are classical for
finite groups~\cite{Isaacs,Robinson}; we state them for arbitrary groups,
since some of our results hold in that generality.

\begin{definition}\label{def:socle-support}
A \emph{minimal normal subgroup} of $G$ is a non-trivial normal subgroup
minimal under inclusion.  Thus the minimal normal subgroups are exactly the
atoms of $\N(G)$.  Write $\MinNorm(G)$ for their set and define the
\emph{socle}
\[
\Soc(G)=\bigl\langle A:A\in\MinNorm(G)\bigr\rangle.
\]
For $N\trnl G$, define its atom support by
\[
\supp(N)=\{A\in\MinNorm(G):A\leq N\}.
\]
The group $G$ is \emph{normally atomic} if every non-trivial normal subgroup
contains a minimal normal subgroup.  It is \emph{monolithic} if it has a
unique minimal normal subgroup.
\end{definition}

The group $\mathds{Z}$ is not normally atomic, whereas all finite groups are. The
Pr\"ufer $p$-group $\Z(p^\infty)\leq\mathds{Q}/\Z$ is normally atomic and
monolithic, with $\langle \frac{1}{p}+\mathds{Z}\rangle$ as its unique minimal
normal subgroup.

A finite group $G$ is monolithic if and only if the intersection of all
non-trivial normal subgroups is non-trivial; equivalently, it is
\emph{subdirectly irreducible} in the sense of universal algebra
\cite[Chapter~II, Section~8]{BurrisSankappanavar}. In general, a non-trivial group is
subdirectly irreducible precisely when it has a least non-trivial normal
subgroup. Thus the equivalence with monolithicity, as defined above, fails for
infinite groups. For example, in $C_p\times \mathds{Z}$, where $p$ is prime,
$C_p\times 0$ is the unique minimal normal subgroup, but the intersection of
all non-trivial normal subgroups is trivial; of course $0\times\Z$ contains
no minimal normal subgroup, so this group is not normally atomic.

The two halves of the dichotomy proved in Section~\ref{sec:intersections} are
the following notions. Essential normal subgroups of arbitrary groups are studied
in~\cite{KonerMitra}; relative uniformity is the group analogue of the
relatively irreducible ideals of~\cite{BlecherGoswami} and of the irreducible
elements of~\cite{CavigliaGoswamiMesiti}.

\begin{definition}\label{def:essential-uniform}
A normal subgroup $E\trnl G$ is \emph{essential} in $G$ if
$E\cap N\neq1$ for every non-trivial $N\trnl G$.  It is \emph{relatively
uniform} in $G$ if any two non-trivial normal subgroups of $G$ contained in
$E$ have non-trivial intersection.
\end{definition}

The central definition of the paper transfers the ring $\M$-ideals
of~\cite{BlecherGoswami} and the $\mu$-elements of~\cite{CavigliaGoswamiMesiti}
to the lattice $\N(G)$.

\begin{definition}\label{def:M}
A normal subgroup $E\trnl G$ is an \emph{$\M$-subgroup} if, for every finite
family $N_1,\ldots,N_r$ of non-trivial normal subgroups with $r\geq2$
satisfying
\[
\bigcap_{i=1}^rN_i\neq1
\quad\text{and}\quad
E\cap N_i\neq1\quad(1\leq i\leq r),
\]
one has
\[
E\cap\bigcap_{i=1}^rN_i\neq1.
\]
The trivial subgroup satisfies this condition vacuously.  Write
\[
\N^{\M}(G)=\{E\trnl G:E\text{ is an $\M$-subgroup}\}.
\]
A group $G$ is \emph{$\M$-full} if $\N^{\M}(G)=\N(G)$.
\end{definition}

Trivially, for every group $G$ the subgroups $1$ and $G$ are $\M$-subgroups,
and so is every minimal normal subgroup (for instance every normal subgroup of
prime order), because two normal subgroups meeting it non-trivially both
contain it. The first part of the following example computes $\N^{\M}(G)$
directly from the definition in a group with three minimal normal subgroups;
the second part exhibits an infinite family of $\M$-full groups.

\begin{example}\label{ex:D60-Sn}
(1)
Let $G=D_{60}$ be the dihedral group of the $30$-gon, generated by $r$ and $s$, satisfying $r^{30}=1=s^2$ and $srs=r^{-1}$. %Observe that $[r^i,s]=r^isr^{-i}s=r^{2i}$, for $1\leq i\leq 30$. 
Every non-trivial normal subgroup of $\langle r\rangle$ contains one of $\langle r^{15}\rangle$, $\langle r^{10}\rangle$, $\langle r^6\rangle$. Moreover, if a normal subgroup $N$ is not contained in $\langle r\rangle$, then $sr^i\in N$ for some $0\leq i<30$.  Since $N$ is normal, $r(sr^i)r^{-1}=sr^{i-2}\in N$ and we get
    \[
    (sr^i)(sr^{i-2})^{-1}=sr^2s=r^{-2}\in N.
    \]

Therefore $\langle r^2\rangle\subseteq N$ and hence $\langle r^{10}\rangle\subseteq N$. Such a normal subgroup $N$ either contains $\langle r^2,s\rangle$ or $\langle r^2,rs\rangle$, according as $i$ is even or odd. The group $D_{60}$ has $11$ normal subgroups
\[
\N(D_{60})=\{1, \langle r^{15}\rangle,\langle r^{10}\rangle,\langle r^6\rangle,\langle r^5\rangle,\langle r^3\rangle,\langle r^2\rangle, \langle r^2,s\rangle,\langle r^2,rs\rangle,\langle r\rangle, D_{60}\}.
\]
The collection of $\M$-subgroups is precisely
\[
\N^\M(D_{60})=\{1,\langle r^{15}\rangle,\langle r^{10}\rangle,\langle r^6\rangle, \langle r\rangle, D_{60}\}.
\]
The obvious $\M$-subgroups are $\{1, \langle r^{15}\rangle, D_{60}\}$. If $E\in \{\langle r^{10}\rangle,\langle r^6\rangle\}$, then for any normal subgroup $N_i$ with $E\cap N_i\neq1$, $E\cap N_i$ is a non-trivial subgroup of $E$, which has prime order, and so $E\leq N_i$. Hence $E\cap \bigcap_i N_i=E\neq 1$. Every non-trivial normal subgroup contains one of the three minimal normal subgroups $\langle r^{15}\rangle$, $\langle r^{10}\rangle$, $\langle r^6\rangle$, all of which lie in $\langle r\rangle$; hence $\langle r\rangle$ meets every non-trivial normal subgroup, that is, $\langle r\rangle$ is essential, and an essential normal subgroup is clearly an $\M$-subgroup. Finally, each of $\langle r^5\rangle$, $\langle r^3\rangle$, $\langle r^2\rangle$, $\langle r^2,s\rangle$, $\langle r^2,rs\rangle$ contains exactly two of the three minimal normal subgroups (the last three contain $\langle r^{10}\rangle$ and $\langle r^{6}\rangle$ but not $\langle r^{15}\rangle$); if $E$ is one of them and $B$ is the minimal normal subgroup not contained in $E$, while $A_1,A_2\leq E$ are the other two, then $N_1=A_1B$ and $N_2=A_2B$ satisfy $N_1\cap N_2=B\neq1$ and $E\cap N_i\geq A_i$, but $E\cap N_1\cap N_2=E\cap B=1$. So none of these five is an $\M$-subgroup. For odd $n\geq 3$, a general formula for $\N^\M(D_{2n})$ can be found in the proof of Proposition~\ref{prop:dihedral-depth}.

(2) For $n\geq 2$, the group $S_n$ is $\M$-full. Indeed, the lattice $\N(S_n)$ is a chain:
\[
\begin{gathered}
\N(S_2)=\{1<S_2\},\qquad
\N(S_3)=\{1<A_3<S_3\},\\
\N(S_4)=\{1<V_4<A_4<S_4\},\qquad
\N(S_n)=\{1<A_n<S_n\}\ \ (n\ge5),
\end{gathered}
\]
where $V_4=\{1,(12)(34),(13)(24),(14)(23)\}$. For each case, $S_n$ is monolithic and $\M$-fullness follows from Theorem \ref{thm:Mfull-classification}, case (4a). 
\end{example}

Since $G$ itself is always an $\M$-subgroup, every normal subgroup lies in at least one $\M$-subgroup, so the
intersection of all $\M$-subgroups containing it is defined. This gives a closure
operator, and it also leads us to name the pieces into which $\N^{\M}(G)$
decomposes.

\begin{definition}\label{def:Mclosure}
For a group $G$ and $N\trnl G$, define
\[
\clM(N)=\bigcap\{E\trnl G:E\text{ is an $\M$-subgroup and }N\leq E\}.
\]
When $\N^{\M}(G)$ is regarded as an ordered set, it is called the
\emph{$\M$-poset}; when it is a lattice, $\vee_{\M}$ denotes its join.  Put
$S=\Soc(G)$ and, for $A\in\MinNorm(G)$, define the \textit{$A$-branch}
\[
\R_A=\{E\trnl G:\supp(E)=\{A\}\},
\qquad
\E(G)=[S,G]_{\N(G)}.
\]
\end{definition}

The numerical invariants we use are lengths in $\N(G)$. Chief series and the
Jordan--H\"older theorem for them are standard~\cite{Isaacs,Robinson}; the length $\sigma(G)$ of the socle as a $G$-group
appears, for instance, in the study of generation of finite
groups~\cite{DallaVoltaLucchiniSocle}.

\begin{definition}(\cite[Chapter 8, Section B]{IsaacsAlgebra})\label{def:lengths}
For a finite group $G$, a \emph{chief series} is a chain
\[
1=G_0<G_1<\cdots<G_n=G
\]
of normal subgroups of $G$ such that no normal subgroup of $G$ lies strictly
between $G_{i-1}$ and $G_i$.  Its number $n$ of factors is the \emph{chief
length} (or normal length) and is denoted $\elln(G)$.  For $N\trnl G$, let
$\lambda_G(N)$ be the length of the interval $[1,N]_{\N(G)}$.  Define the
\emph{socle length}
\[
\sigma(G)=\lambda_G(\Soc(G)).
\]
An \emph{$\M$-chain} is a chain all of whose terms belong to $\N^{\M}(G)$,
and the \emph{$\M$-depth} of $G$ is
\[
\md(G)=\max\{r:1=E_0<E_1<\cdots<E_r=G,\ E_i\in\N^{\M}(G)\}.
\]
\end{definition}

Abelian minimal normal subgroups are irreducible modules, and the
combinatorics of the socle is that of projective spaces over endomorphism
rings. We fix the corresponding module-theoretic language.

\begin{definition}\label{def:module-notation}
For a prime $p$, $\F_p$ denotes the field with $p$ elements and $\F_pG$
denotes the group algebra.  A left module is \emph{simple} (or
\emph{irreducible}) if it is non-zero and has no proper non-zero submodule.  For
a ring $R$ and a left $R$-module $A$, $\End_R(A)$ denotes its endomorphism
ring.  A module is \emph{semisimple} if it is a direct sum of simple
submodules.  A semisimple module is \emph{homogeneous} if all of its simple
summands are isomorphic; an \emph{isotypic decomposition} is a direct-sum
decomposition into homogeneous components belonging to pairwise
nonisomorphic simple-module types.  For a division ring $D$,
$\mathds P^{n-1}(D)$ denotes the set of one-dimensional right $D$-subspaces
of $D^n$; in particular $\mathds P^1(D)$ is the projective line.
\end{definition}

For direct products the central $p$-torsion of the factors plays a special
role, and we introduce notation to count the minimal normal subgroups lying in
a given normal subgroup.

\begin{definition}\label{def:product-notation}
For $G=G_1\times\cdots\times G_t$, let $G_i^*$ be the $i$-th coordinate
subgroup, and for $A\leq G_i$ let $\widehat A\leq G_i^*$ be its coordinate
copy.  For each prime $p$, put
\[
V_{i,p}=\{z\in Z(G_i):z^p=1\},
\qquad
V_p=\prod_{i=1}^tV_{i,p}.
\]
Each $V_p$ is viewed as an $\F_p$-vector space, and a \emph{line} in $V_p$
means a one-dimensional $\F_p$-subspace.  Such a line is \emph{coordinate}
if it is contained in the coordinate copy of some $V_{i,p}$, and it is
\emph{non-coordinate} otherwise.  If
$H\trnl G_1\times\cdots\times G_t$, define
\[
a(H)=\#\{(i,A):A\in\MinNorm(G_i),\ A\nleq Z(G_i),\ \widehat A\leq H\}
\]
and
\[
d_p(H)=\dim_{\F_p}(H\cap V_p).
\]
For a direct power $X^k$, $X_i$ denotes the $i$-th coordinate copy of $X$.
\end{definition}

Crown-based powers were introduced by Dalla Volta and
Lucchini~\cite{DallaVoltaLucchini}, building on the crowns of
Gasch\"utz~\cite{Gaschutz}; for primitive groups see~\cite{DixonMortimer}.

\begin{definition}\label{def:crowns}
For a subgroup $M\leq L$, its core is
\[
\operatorname{core}_L(M)=\bigcap_{x\in L}xMx^{-1}.
\]
The subgroup $M$ is \emph{core-free} if $\operatorname{core}_L(M)=1$.  A
finite group $L$ is \emph{primitive} if it has a core-free maximal subgroup.
If $L$ is finite, monolithic, and primitive with unique minimal normal
subgroup $A=\Soc(L)$, then for $k\geq1$ its \emph{crown-based power of size
$k$} is
\[
L_k=\{(x_1,\ldots,x_k)\in L^k:x_1A=\cdots=x_kA\}.
\]
Throughout statements concerning $L_k$, write $V=A^k$ and let
$A_1,\ldots,A_k$ be the coordinate copies of $A$ in $V$.
\end{definition}

Finally, for the counting results on abelian, nilpotent, cyclic and dihedral
groups we use the following notation.

\begin{definition}\label{def:abelian-notation}
For a finite abelian group $G$, additive notation may be used.  Write
\[
G=\bigoplus_{p\mid |G|}G_p,
\qquad
G_p\cong C_{p^{\lambda_{p,1}}}\oplus\cdots\oplus
C_{p^{\lambda_{p,r_p}}},
\]
where $G_p$ is the Sylow $p$-subgroup and
$\lambda_{p,1}\geq\cdots\geq\lambda_{p,r_p}\geq1$.  Put
\[
G_p[p]=\{x\in G_p:px=0\},
\qquad
u_p(k)=\sum_{i=1}^{r_p}\min\{k,\lambda_{p,i}\}\quad(k\geq0),
\]
with $u_p(0)=0$.  Let $c_p(G)$ be the number of non-trivial cyclic
$p$-subgroups of $G$, and set
\[
r(G)=\sum_{p\mid |G|}r_p.
\]
\end{definition}

\begin{definition}\label{def:arithmetic-notation}
For a positive integer $n$, $\Omega(n)$ is the number of prime factors of
$n$ counted with multiplicity and $\omega(n)$ is the number of distinct prime
factors.  If $G$ is a finite $p$-group, set
\[
\Omega_1(Z(G))=\{z\in Z(G):z^p=1\}.
\]
If $G$ is a finite nilpotent group and $P_p$ is its Sylow $p$-subgroup, set
\[
\rho(G)=\sum_{p\mid |G|}\dim_{\F_p}\{z\in Z(P_p):z^p=1\}.
\]
For $n\geq3$, $D_{2n}$ denotes the dihedral group of order $2n$ with
presentation
\[
D_{2n}=\langle r,s:r^n=s^2=1,\ srs=r^{-1}\rangle.
\]
\end{definition}

\section{The dichotomy and the finite classification}
\label{sec:intersections}

In this section we prove that $\M$-subgroups are exactly the essential and
the relatively uniform normal subgroups, and deduce Theorem~A. Everything
rests on the modularity of $\N(G)$ and on its consequence for lengths.

The following lemma is well known (\textit{cf.} \cite[Chapter~4, Example~4.6(5)]{Davey}) and will be used throughout. We give a short proof to make the paper self-contained. 
\begin{lemma}\label{lem:modular}
If $A,B,C\trnl G$ and $A\leq C$, then
\[
 A(B\cap C)=AB\cap C.
\]
In particular, $\N(G)$ is a modular lattice.
\end{lemma}

\begin{proof}
The inclusion $A(B\cap C)\leq AB\cap C$ is immediate. If
$ab\in AB\cap C$, with $a\in A$ and $b\in B$, then $a\in C$ implies
$b=a^{-1}(ab)\in B\cap C$. Thus $ab\in A(B\cap C)$, proving the
identity and hence modularity.
\end{proof}

For finite $G$, modularity has a numerical consequence that we use
repeatedly: lengths in $\N(G)$ are additive along chief series and satisfy the rank identity of
modular lattices~\cite{Birkhoff}. In particular, every maximal chain in an
interval of $\N(G)$ has the same length (the Jordan--H\"older theorem for
chief series).

\begin{lemma}\label{lem:relative-length}
Let $G$ be a finite group and let $H,K\trnl G$.  Then
\begin{align}
\elln(G)&=\lambda_G(H)+\elln(G/H),\label{eq:length-quotient}\\
\lambda_G(HK)+\lambda_G(H\cap K)
 &=\lambda_G(H)+\lambda_G(K).\label{eq:modular-rank}
\end{align}
\end{lemma}

\begin{proof}
A maximal chain in $[1,H]_{\N(G)}$ concatenated with the inverse image
of a chief series of $G/H$ is a chief series of $G$. The
Jordan--H\"older theorem for chief series therefore gives
\eqref{eq:length-quotient}.
For \eqref{eq:modular-rank}, the correspondence (diamond) isomorphism
$[H\cap K,K]_{\N(G)}\cong[H,HK]_{\N(G)}$ gives
\[
\lambda_G(K)-\lambda_G(H\cap K)
 =\lambda_G(HK)-\lambda_G(H).
\]
Rearranging proves the claim.
\end{proof}

We first observe that in Definition~\ref{def:M} it suffices to test pairs of
normal subgroups; the ring-theoretic analogue is proved
in~\cite{BlecherGoswami}, where it is also observed that this contrasts with
Banach spaces: the Alfsen--Effros characterisation of $\M$-ideals by
intersections of balls can be restricted to three balls, but not to
two~\cite{AlfsenEffros,HWW}.

\begin{proposition}\label{prop:pairs}
A normal subgroup $E\trnl G$ is an $\M$-subgroup if and only if it satisfies
the defining implication for pairs $N_1,N_2$.
\end{proposition}

\begin{proof}
Necessity follows by taking $r=2$. Conversely, assume the implication
holds for pairs, and proceed by induction on $r\geq2$. Given
$N_1,\ldots,N_r$ as in Definition~\ref{def:M}, put
$N=\bigcap_{i=1}^{r-1}N_i$. Then $N\cap N_r\neq1$; by induction
$E\cap N\neq1$, and by assumption $E\cap N_r\neq1$. Applying the
pairwise implication to $N$ and $N_r$ yields
$E\cap\bigcap_{i=1}^rN_i\neq1$.
\end{proof}

With pairs at hand, the $\M$-condition splits into two classical conditions.
For rings this is the main theorem of~\cite{BlecherGoswami}, and for modular
quantales it is the characterisation of $\mu$-elements
in~\cite{CavigliaGoswamiMesiti}; as $\N(G)$ is in general not an integral commutative quantale (see the
discussion before Theorem~\ref{thm:Mfull-udim}), we give the short direct
proof.

\begin{theorem}\label{thm:dichotomy}
A normal subgroup $E\trnl G$ is an $\M$-subgroup if and only if it is
essential in $G$ or relatively uniform in $G$.
\end{theorem}

\begin{proof}
    Suppose $E$ is an $\M$-subgroup that is neither essential nor relatively uniform. Then there are non-trivial $N_1$, $N_2$, $N_3\trianglelefteq G$ such that $E\cap N_1=1$, $N_2$, $N_3\leq E$ and $N_2\cap N_3=1$. Now, by the modular law, for $i=2,3$, we have $$E\cap N_iN_1=N_i(E\cap N_1)=N_i.$$ Therefore, for the pairs of normal subgroups $N_2N_1$ and $N_3N_1$, $E$ fails to be an $\M$-subgroup: indeed $N_2N_1\cap N_3N_1\geq N_1\neq1$ and $E\cap N_iN_1=N_i\neq1$, whereas $E\cap N_2N_1\cap N_3N_1=N_2\cap N_3=1$.

   Conversely, if $E$ is essential, then for any family $\{N_i\}_{i\in \Lambda}$ of normal subgroups with $\bigcap_{i\in \Lambda}N_i\neq 1$, we have $E\cap N_i\neq 1$ and $E\cap \bigcap_{i\in \Lambda}N_i\neq 1$ by the definition of essentiality. Moreover, if $E$ is relatively uniform and $N_1,\ N_2$ are as in Proposition~\ref{prop:pairs}, then $E\cap N_1$ and $E\cap N_2$ are non-trivial normal subgroups contained in $E$, so $E\cap N_1\cap N_2\neq 1$. Therefore, in both cases $E$ is an $\M$-subgroup.
\end{proof}

Theorem \ref{thm:dichotomy} confirms that each $\M$-subgroup is either essential or relatively uniform. The next example shows that neither class contains the other, that they can overlap, and that a normal subgroup may lie in neither.
\begin{example} In the group $C_{p^n}$ the subgroups form a chain, so any two non-trivial subgroups
      meet non-trivially. Hence every non-trivial subgroup of $C_{p^n}$ is both
      essential and relatively uniform.
Now let $G=C_{30}$, and for any divisor $d$ of $30$ let $C_d$ be its cyclic subgroup of order $d$. Clearly,
$C_{d_1}\cap C_{d_2}=C_{\gcd(d_1,d_2)}$.
\begin{enumerate}
\item The subgroup $C_2$ is relatively uniform but not essential, since $C_2\cap C_3=1$.
\item The whole group $G$ is essential, but not relatively uniform, since $C_2,C_3\le G$ and $C_2\cap C_3=1$.
\item The subgroup $C_6$ is neither essential nor relatively uniform, since $C_6\cap C_5=1$, and $C_2,C_3\le C_6$ with $C_2\cap C_3=1$. So $C_6$ is not an $\M$-subgroup.

\end{enumerate}
\end{example}

In normally atomic groups relative uniformity and essentiality can be read off
from the minimal normal subgroups. We need the following elementary fact.

\begin{lemma}\label{lem:distinct-atoms}
Distinct minimal normal subgroups $A,B$ of a group satisfy
$A\cap B=1$ and $[A,B]=1$.
\end{lemma}

\begin{proof}
The normal subgroup $A\cap B$ is contained in both $A$ and $B$, and minimality forces $A\cap B=1$. Furthermore, the subgroup $[A,B]$ is generated by $[a,b]$ for $a\in A,\ b\in B$, and each $[a,b]\in A\cap B$. Therefore $[A,B]\leq A\cap B=1.$
\end{proof}

\begin{theorem}\label{thm:atomic-classification}
Let $G$ be normally atomic and let $E\trnl G$.  Then $E$ is an
$\M$-subgroup if and only if one of the following holds:
\begin{enumerate}
\item $E=1$;
\item $E$ contains exactly one minimal normal subgroup of $G$;
\item $E$ contains every minimal normal subgroup of $G$.
\end{enumerate}
\end{theorem}

\begin{proof}
    %Let $E$ be an $\M$-subgroup of $G$. If $E=1$, we are done. Suppose that $E\ne1$. As $G$ is normally atomic, there exists $N\in\MinNorm(G)$ such that $N\leq E$. If $N$ is the only such minimal normal subgroup, then we are done. So, now suppose that there exists another $N'\in\MinNorm(G)$ such that $N'\leq E$. Now, consider $M\in\MinNorm(G)$. Then, by Lemma \ref{lem:distinct-atoms}, we have $M\cap N=1$, $M\cap N'=1$ and $N\cap N'=1$. Also, by \autoref{lem:distinct-atoms}, they commute, i.e., $[N,M]=1$ and $[N',M]=1$. Now let $H=NM$ and $K=N'M$. Then, $H,K\trnl G$. Then, 
   % \[
   %    H=N\times M \quad \text{and}\quad K=N'\times M
   % \]
   % Furthermore, it is clear that $H\cap K=M$. Now, as $E$ is an $\M$-subgroup, by \autoref{prop:pairs} we have
   % \[
   %    E\cap(H\cap K)=E\cap M\ne1
   % \]
   % But $E\cap M\trnl G$ and $E\cap M\leq M$, so $E\cap M=1$ or $E\cap M=M$. As $E\cap M\ne1$ so $E\cap M=M$. Hence, $M\leq E$. As $M$ was arbitrary, it follows that, $E$ contains every minimal normal subgroup of $G$.\par

     Let $E\neq1$ be an $\M$-subgroup. As $G$ is normally atomic, $E$ contains at least one minimal normal subgroup, so if $E$ does not contain two distinct ones, then (2) holds. Suppose the $\M$-subgroup $E\neq 1$ contains two distinct minimal normal subgroups $N$ and $N'$. By Lemma \ref{lem:distinct-atoms}, $N\cap N'=1$. Then $E$ is not relatively uniform, and so Theorem \ref{thm:dichotomy} forces $E$ to be essential. Then $E\cap M\neq 1$ for every $M\in \MinNorm(G)$ and so $E\cap M=M$, that is, $M$ is a subgroup of $E$. Hence $E$ contains every minimal normal subgroup of $G$.
    
    Conversely, if $E=1$, then vacuously $E$ is an $\M$-subgroup. Now, let $E$ contain exactly one minimal normal subgroup, say $N$. Consider two non-trivial normal subgroups $N_{1},N_{2}\trnl G$ such that $N_{1},N_{2}\leq E$. As $G$ is normally atomic, there exist $M_{1},M_{2}\in\MinNorm(G)$ such that $M_{1}\leq N_{1}$ and $M_{2}\leq N_{2}$. Thus, $M_{1},M_{2}\leq E$. But $E$ contains a unique minimal normal subgroup $N$. Therefore, $M_{1}=N=M_{2}$ and hence $N\leq N_{1}$ and $N\leq N_{2}$. Thus, $N\leq N_{1}\cap N_{2}$ and hence $N_{1}\cap N_{2}\ne1$. Thus $E$ is relatively uniform in $G$ and so by Theorem \ref{thm:dichotomy}, $E$ is an $\M$-subgroup. Now, suppose $E$ contains every minimal normal subgroup of $G$. Now, consider a normal subgroup $H\trnl G$, where $H\ne1$. Now, as $G$ is normally atomic, there exists $K\in\MinNorm(G)$ such that $K\leq H$. But by hypothesis, $K\leq E$. So, $K\leq E\cap H$ and hence $E\cap H\ne1$. Thus $E$ is essential in $G$ and so by Theorem \ref{thm:dichotomy} $E$ is an $\M$-subgroup.
\end{proof}

\begin{corollary}\label{cor:finite-classification}
Let $G$ be finite and let $E\trnl G$.  Then $E$ is an $\M$-subgroup if and
only if
\[
E=1,\qquad |\supp(E)|=1,
\qquad\text{or}\qquad \Soc(G)\leq E.
\]
\end{corollary}

\begin{proof}
As $G$ is finite, it is normally atomic. Now the result follows immediately from Theorem~\ref{thm:atomic-classification}, as $|\supp(E)|=1$ means that $E$ contains a unique minimal normal subgroup and $\Soc(G)\leq E$ means that $E$ contains every minimal normal subgroup of $G$, as $\Soc(G)$ is the smallest subgroup containing all minimal normal subgroups of $G$.
\end{proof}

For a finite group $G$, Corollary~\ref{cor:finite-classification} gives us a tool to check whether a subgroup is an $\M$-subgroup. In the next example, we use it to show that the centre and the derived subgroup of a group need not be $\M$-subgroups. 
\begin{example}
   (1) Let $G=K_4\times S_3$, where $K_4=\{1,a,b,ab\}\cong C_2\times C_2$. Then $Z(G)=K_4\times 1\neq 1$, $|\supp(Z(G))|=3$ (the three subgroups of order $2$ of $Z(G)$ are central of prime order, hence minimal normal), and $Z(G)$ does not contain the minimal normal subgroup $1\times A_3$. Therefore, by Corollary~\ref{cor:finite-classification}, $Z(G)$ is not an $\M$-subgroup.
   
        (2) Take $G=D_{60}$ (\textit{cf.} Example~\ref{ex:D60-Sn}), then for any $1\leq i\leq 30$, 
        \[
        [r^i,s]=r^isr^{-i}s=r^{2i},
        \] and hence $[D_{60},D_{60}]\geq\langle r^2\rangle$; as $D_{60}/\langle r^2\rangle$ has order $4$ and is therefore abelian, $[D_{60},D_{60}]=\langle r^2\rangle$. But $\langle r^2\rangle$ is not an $\M$-subgroup, since it contains the two minimal normal subgroups $\langle r^{10}\rangle$ and $\langle r^{6}\rangle$, while the minimal normal subgroup $\langle r^{15}\rangle$ is not inside $\langle r^2\rangle$. Hence by Corollary~\ref{cor:finite-classification} the derived subgroup of $D_{60}$ is not an $\M$-subgroup. 
\end{example}
We now turn to the closure system formed by $\M$-subgroups. Finite
intersections are always harmless.

\begin{proposition}\label{prop:meet}
The intersection of finitely many $\M$-subgroups of an arbitrary group is an
$\M$-subgroup.
\end{proposition}

\begin{proof}
 Let $E_1$, $E_2$ be $\M$-subgroups and for $r\ge 2$, let $N_1,\ldots, N_r$ be non-trivial normal subgroups such that $\bigcap_i N_i\ne 1$ and $E_1\cap E_2\cap N_i\ne 1$ for all $i$. Since $E_1\cap N_i\neq 1$ and $E_1$ is an $\M$-subgroup, then $E_1\cap\bigcap_i N_i\neq 1$. Now each $E_1\cap N_i$ meets $E_2$ non-trivially, and the non-trivial normal subgroups $E_1\cap N_1,\ldots,E_1\cap N_r$ have non-trivial intersection $E_1\cap\bigcap_iN_i$; then since $E_2$ is an $\M$-subgroup, $E_1\cap E_2\cap \bigcap_i N_i\neq 1$. Therefore $E_1\cap E_2$ is an $\M$-subgroup, and the case of finitely many $\M$-subgroups follows by induction. 
\end{proof}

\begin{proposition}\label{prop:Mclosure}
The map $\clM:\N(G)\to\N(G)$ is extensive, order preserving, and idempotent,
and every $\M$-subgroup is a fixed point of $\clM$. If $G$ is finite, then
the fixed points are exactly the $\M$-subgroups, and explicitly
\[
\clM(N)=
\begin{cases}
N,&N=1\text{ or }|\supp(N)|=1,\\
N\Soc(G),&|\supp(N)|\geq2.
\end{cases}
\]
\end{proposition}

\begin{proof}
    Let $N\trnl G$. Now, if $E\in\N^\M(G)$ such that $N\leq E$, then 
    \[
       N\leq\bigcap\{E\in\N^\M(G):N\leq E\}=\clM(N)
    \]
    Hence, $\clM$ is extensive. Now suppose $N,M\in\N(G)$ such that $N\leq M$. If $E\in\N^\M(G)$ such that $M\leq E$, then $N\leq E$. Thus, $\{E\in\N^\M(G):M\leq E\}\subseteq\{E\in\N^\M(G):N\leq E\}$. As intersection reverses the order, we get
    \[
       \bigcap\{E\in\N^\M(G):N\leq E\}\leq\bigcap\{E\in\N^\M(G):M\leq E\}
    \]
    Thus, $\clM(N)\leq\clM(M)$. Hence, $\clM$ is order preserving.\par
    Now let $N\in\N(G)$. Define, $\mathcal{A}(N)=\{E\in\N^\M(G):N\leq E\}$ and $\mathcal{C}=\clM(N)$. If $H\in\mathcal{A}(N)$ then $\mathcal{C}\leq H$. Thus, $H\in\mathcal{A}(\mathcal{C})$. Conversely if $H\in\mathcal{A}(\mathcal{C})$, then $N\leq\mathcal{C}\leq H$, where the first inequality comes from extensiveness of $\clM$. Thus, $H\in\mathcal{A}(N)$. Thus, $\mathcal{A}(N)=\mathcal{A}(\mathcal{C})$. Thus,
    \[
    \clM(N)=\bigcap\mathcal{A}(N)=\bigcap\mathcal{A}(\mathcal{C})=\clM(\clM(N))=\clM^{2}(N)
    \]
    Thus, $\clM$ is idempotent.\par
    Now suppose, $N$ is an $\M$-subgroup. Then $N\in\{E\in\N^\M(G):N\leq E\}$. Thus, $\clM(N)\leq N$. But by extensiveness we get $N\leq\clM(N)$. Thus, $N=\clM(N)$.
    
   Now let $G$ be finite. Then $\clM(N)$ is the intersection of finitely many $\M$-subgroups, hence an $\M$-subgroup by Proposition~\ref{prop:meet}. In particular, if $\clM(N)=N$, then $N$ is an $\M$-subgroup, so the fixed points are exactly the $\M$-subgroups. If $N=1$ or $|\supp(N)|=1$, then $N$ is an $\M$-subgroup by Corollary~\ref{cor:finite-classification}, and so $\clM(N)=N$. Lastly, if $|\supp(N)|\ge 2$, every $\M$-subgroup $E\ge N$ has $|\supp(E)|\ge 2$, so by Corollary~\ref{cor:finite-classification} such an $E$ must contain $\Soc(G)$ and hence $N\Soc(G)\le E$. Clearly, $N\Soc(G)$ is itself an $\M$-subgroup and the smallest $\M$-subgroup containing both $N$ and $\Soc(G)$.
\end{proof}

Finiteness of the group is assumed in Corollary~\ref{cor:finite-classification} and in the second part of Proposition~\ref{prop:Mclosure}, and Proposition~\ref{prop:meet} is about finitely many $\M$-subgroups. The next example shows that these assumptions cannot be dropped, even for finitely generated abelian groups.

\begin{example}\label{ex:Z3}
	Let $G=\Z^3$. Consider the subgroups
	\[
	P=\Z\times\Z\times0\qquad\text{and}\qquad H_m=\Z\times\Z\times m\Z\quad(m\ge1).
	\]
	
	Each $H_m$ is an essential subgroup. Indeed, let $N$ be a non-zero subgroup and $0\ne x\in N$. Then $mx\in N\cap H_m$. So $H_m$ is essential, and hence an $\M$-subgroup. Moreover, $P$ is not an $\M$-subgroup. To this end, let $N_1=\Z\times0\times\Z$ and $N_2=0\times\Z\times\Z$. Then
	\[
	N_1\cap N_2=0\times0\times\Z\ne0,\qquad P\cap N_1=\Z\times0\times0\ne0,\qquad P\cap N_2=0\times\Z\times0\ne0,
	\]
	but $P\cap N_1\cap N_2=0$. We observe the following:
	
	\begin{enumerate}
		\item Clearly, $\bigcap_{m\ge1}H_m=P$. Therefore, an intersection of infinitely many $\M$-subgroups need not be an $\M$-subgroup. Therefore, Proposition~\ref{prop:meet} holds only for finitely many.
		\item The only subgroups containing $P$ are $P$ and the $H_m$ with $m\geq1$. Indeed, if $H$ is any subgroup containing $P$, then $(a,b,c)\in H$ if and only if $(0,0,c)\in H$. Therefore $\{c\mid (0,0,c)\in H\}$ is a subgroup of $\Z$ and so it is of the form $m\Z$ for some $m\ge 0$. Hence $H=H_m$ for $m>0$ and $H=P$ for $m=0$. Now, each $H_m$ with $m\geq1$ is an $\M$-subgroup containing $P$, and $P$ is not, so by (1)
		\[
		P\le\clM(P)\le\bigcap_{m\ge1}H_m=P,
		\]
		that is, $\clM(P)=P$.
		Thus $P$ is a fixed point of the closure operator $\clM$, but it is not an $\M$-subgroup. Therefore, the assertion of Proposition~\ref{prop:Mclosure} that the fixed points are exactly the $\M$-subgroups needs $G$ to be finite.

        \item Let $N$ be a non-trivial subgroup of $G$. If $x\in N\setminus\{0\}$, then $\langle2x\rangle$ is a non-zero subgroup of $N$ not containing $x$. Hence $G$ has no minimal subgroups, so $\Soc(G)=0\le P$, so the condition $\Soc(G)\le E$ of Corollary~\ref{cor:finite-classification} holds for $E=P$, although $P$ is not an $\M$-subgroup. 
		
	\end{enumerate}
\end{example}

For finite groups we can now assemble Theorem~A.

\begin{corollary}\label{cor:Mlattice}
For every finite group $G$, the ordered set $\N^{\M}(G)$ is a lattice.  Its
meet is ordinary intersection, and its join is
\[
 E\mathbin{\vee_{\M}}F=\clM(EF).
\]
Thus the $\M$-join agrees with the product whenever $EF$ is an
$\M$-subgroup, and otherwise it is $EF\Soc(G)$.
\end{corollary}

\begin{proof} 
This follows as a consequence of Proposition~\ref{prop:meet} and Proposition~\ref{prop:Mclosure}: the meet of $E,F\in\N^{\M}(G)$ in $\N^{\M}(G)$ is $E\cap F$, and $\clM(EF)$ is the least $\M$-subgroup containing $E$ and $F$. If $EF$ is not an $\M$-subgroup, then $EF\neq1$ and $|\supp(EF)|\neq1$, so $|\supp(EF)|\geq2$ and $\clM(EF)=EF\Soc(G)$.
\end{proof}
\section{The \texorpdfstring{$\M$}{M}-poset and the no-loss problem}
\label{sec:noloss}

Passing from $\N(G)$ to $\N^{\M}(G)$ discards the normal subgroups that
contain at least two, but not all, minimal normal subgroups. In this section we show
that $\N^{\M}(G)$ remembers the minimal normal subgroups, the socle and the
filter above it, and we determine exactly when nothing is lost
(Theorem~B). In the rest of this section $G$ is finite.

By Theorem~A, the $\M$-poset of a finite group splits into $1$, the branches
lying over the individual minimal normal subgroups, and the filter above the
socle.

\begin{proposition}\label{prop:poset-decomposition}
Suppose that $G$ is finite and has at least two minimal normal subgroups.
Then
\[
\N^{\M}(G)=\{1\}\sqcup
\coprod_{A\in\MinNorm(G)}\R_A\sqcup\E(G).
\]
Each $\{1\}\cup\R_A$ is closed under intersections, and $\E(G)$ is a
sublattice of $\N(G)$.  If $E\in\R_A$ and $F\in\R_B$ with $A\neq B$, then
$EF$ is an $\M$-subgroup if and only if $\Soc(G)\leq EF$. 
\end{proposition}

\begin{proof}
   By Corollary~\ref{cor:finite-classification} the $\M$-subgroups are $1$, the members of the $\R_A$, and the members of $\E(G)$, and they are disjoint. Clearly, $\R_A\cap\R_B=\emptyset$ for $A\ne B$, and each member of $\E(G)$ contains $S$, hence at least two minimal normal subgroups. Since $\supp(E\cap F)=\supp(E)\cap\supp(F)$, the intersection of two members of $\R_A$ lies in $\R_A$. Moreover, $\E(G)$ is an interval of $\N(G)$, hence a sublattice. Finally, if $E\in\R_A$ and $F\in\R_B$ with $A\ne B$, then $A,B\le EF$, so $|\supp(EF)|\ge2$, and by Corollary~\ref{cor:finite-classification}, $EF$ is an $\M$-subgroup if and only if $\Soc(G)\le EF$.
\end{proof}

The decomposition is visible inside the abstract poset $\N^{\M}(G)$: the
minimal normal subgroups, the socle and the branches can all be recovered from
it. The last part of the next theorem answers the no-loss problem in terms of
products of minimal normal subgroups.

\begin{theorem}\label{thm:reconstruction}
Let $G$ be finite, put $P=\N^{\M}(G)$, ordered by inclusion, and put $S=\Soc(G)$.
\begin{enumerate}
\item The atoms of $P$ are exactly the minimal normal subgroups of $G$.
\item If $G$ has at least two atoms, then $S$ is the least element of $P$
      lying above every atom of $P$.
\item Under the same hypothesis, the interval $[S,G]_{\N(G)}$ is the
      principal filter $[S,G]_P$, and the branch $\R_A$ is recovered from
      $P$ as
      \[
      \{E\in P:A\leq E\text{ and }S\nleq E\}.
      \]
\item The following conditions are equivalent:
      \begin{enumerate}[label=\upshape(\alph*)]
      \item $P=\N(G)$;
      \item $P$ is closed under products of its elements;
      \item $AB=S$ for every two distinct $A,B\in\MinNorm(G)$.
      \end{enumerate}
\end{enumerate}
\end{theorem}

\begin{proof}
    (1) Each $A\in\MinNorm(G)$ lies in $P$ by Corollary~\ref{cor:finite-classification}, and no normal subgroup lies strictly between $1$ and $A$. Therefore, $A$ is an atom of $P$. Since $G$ is normally atomic, an atom $B$ of $P$ contains some $A\in\MinNorm(G)\subseteq P$, forcing $B=A$.
 
(2) $S\in P$ by Corollary~\ref{cor:finite-classification}, and by (1) an element of $P$ lies above every atom if and only if it contains $S$.
 
(3) By Corollary~\ref{cor:finite-classification} every normal subgroup containing $S$ lies in $P$. Therefore, $[S,G]_{\N(G)}=[S,G]_P$. If $E\in\R_A$, then $A\le E$, and
$S\nleq E$ because $E$ misses the other atoms. Conversely, if $E\in P$, $A\le E$ and $S\nleq E$, then $|\supp(E)|=1$ by Corollary~\ref{cor:finite-classification}, and so $\supp(E)=\{A\}$.
 
(4) (a)$\Rightarrow$(b) is clear.

(b)$\Rightarrow$(c): for distinct $A,B\in\MinNorm(G)\subseteq P$ we get $AB\in P$ and $|\supp(AB)|\ge2$, and hence $S\le AB\le S$ by Corollary~\ref{cor:finite-classification}.

(c)$\Rightarrow$(a): a normal subgroup
containing two distinct minimal normal subgroups $A,B$ contains $AB=S$, so
every normal subgroup lies in $P$ by Corollary~\ref{cor:finite-classification}.
\end{proof}

Condition (c) of Theorem~\ref{thm:reconstruction}(4) involves products of
minimal normal subgroups, and to make it effective we need to know how many
minimal normal subgroups the socle has. Abelian minimal normal subgroups are
simple modules, and the simple submodules of a homogeneous semisimple module
form a projective space. We record this classical fact with a proof, since
the explicit parametrisation is used again in Section~\ref{sec:crowns}.

\begin{lemma}\label{lem:projective-submodules}
Let $R$ be a finite ring, let $A$ be a simple left $R$-module, and put
$D=\End_R(A)$.  Then $D$ is a finite division ring.  For $n\geq1$, the simple submodules of $A^n$ are naturally
parametrized by the one-dimensional right $D$-subspaces of $D^n$.
Consequently, their number is
\[
\frac{|D|^n-1}{|D|-1}.
\]
For $n=2$, any two distinct simple submodules have direct sum $A^2$.
\end{lemma}

\begin{proof}
    
    By Schur's Lemma, $D$ is a division ring and every non-zero $R$-module homomorphism between simple $R$-modules is an isomorphism. If $0\ne a\in A$, then $Ra$ is a non-zero submodule of $A$, so $A=Ra$. As $R$ is finite, $A$ is finite, and hence $D$ is finite. Here the product in $D$ is the composition, and $D^{n}$ is a right $D$-vector space via $(f_{1},\cdots,f_{n})g=(f_{1}g,\cdots,f_{n}g)$.\par
    Let $\pi_{i}:A^{n}\rightarrow A$ be the $i$-th projection. For $f=(f_{1},\cdots,f_{n})\in D^{n}$ define $\Phi_{f}\in\Hom_{R}(A,A^{n})$ by $\Phi_{f}(a)=(f_{1}(a),\cdots,f_{n}(a))$. As $\pi_{i}\circ\Phi_{f}=f_{i}$, the map $f\mapsto\Phi_{f}$ is a bijection from $D^{n}$ onto $\Hom_{R}(A,A^{n})$, with inverse $\varphi\mapsto(\pi_{1}\circ\varphi,\cdots,\pi_{n}\circ\varphi)$. Moreover, for $g\in D$ we have
    \[
       \Phi_{fg}=\Phi_{f}\circ g.
    \]
    Now let $f\ne0$. Then $\Phi_{f}\ne0$, so $\ker\Phi_{f}=\{0\}$ as $A$ is simple. Hence $\text{im}\,\Phi_{f}\cong A$ is a simple submodule of $A^{n}$. Conversely, let $S$ be a simple submodule of $A^{n}$. As $S\ne\{0\}$, $\pi_{i}(S)\ne\{0\}$ for some $i$, and by Schur's Lemma $\pi_{i}|_{S}:S\rightarrow A$ is an isomorphism. In particular, $S\cong A$. If $\iota:S\hookrightarrow A^{n}$ is the inclusion, then $\iota\circ(\pi_{i}|_{S})^{-1}=\Phi_{f}$ for some $0\ne f\in D^{n}$, and $S=\text{im}\,\Phi_{f}$.\par
    Now let $0\ne f,f'\in D^{n}$. We claim that $\text{im}\,\Phi_{f}=\text{im}\,\Phi_{f'}$ if and only if $f'=fg$ for some $g\in D^{\times}$, that is, if and only if $fD=f'D$. If $f'=fg$ with $g\in D^{\times}$, then $g$ is onto and so $\text{im}\,\Phi_{f'}=\text{im}\,(\Phi_{f}\circ g)=\text{im}\,\Phi_{f}$. Conversely, if $\text{im}\,\Phi_{f}=\text{im}\,\Phi_{f'}=S$, then $\Phi_{f},\Phi_{f'}:A\rightarrow S$ are isomorphisms. So $g=\Phi_{f}^{-1}\circ\Phi_{f'}\in D^{\times}$ and $\Phi_{fg}=\Phi_{f}\circ g=\Phi_{f'}$, which gives $f'=fg$. Hence
    \[
       fD\longmapsto\text{im}\,\Phi_{f}
    \]
    is a well-defined bijection from $\mathds{P}^{n-1}(D)$ onto the set of simple submodules of $A^{n}$.\par
    Every non-zero $f\in D^{n}$ lies in exactly one line, namely $fD$. Also $g\mapsto fg$ is injective, since $fg=fg'$ with $g\ne g'$ would give $f=f(g-g')(g-g')^{-1}=0$. So each line has exactly $|D|-1$ non-zero elements, and thus the number of simple submodules of $A^{n}$ is
    \[
       |\mathds{P}^{n-1}(D)|=\frac{|D|^{n}-1}{|D|-1}.
    \]
    Finally, let $A_{1}\ne A_{2}$ be simple submodules of $A^{2}$. If $A_{1}\cap A_{2}\ne\{0\}$, then $A_{1}\cap A_{2}=A_{1}$ as $A_{1}$ is simple, so $A_{1}\leq A_{2}$ and hence $A_{1}=A_{2}$ as $A_{2}$ is simple, which is a contradiction. So $A_{1}\cap A_{2}=\{0\}$, and as $A_{1}\cong A\cong A_{2}$, we get
    \[
       |A_{1}+A_{2}|=|A_{1}|\,|A_{2}|=|A|^{2}=|A^{2}|.
    \]
    Hence, $A_{1}\oplus A_{2}=A^{2}$.
\end{proof}

Applying Lemma~\ref{lem:projective-submodules} to the isotypic components of
the abelian part of the socle, we obtain the complete geometry of the set of
minimal normal subgroups of a finite group and an expression for the socle
length.

\begin{proposition}\label{prop:socle-geometry}
Let $G$ be a non-trivial finite group.  Let $B_1,\ldots,B_b$ be its non-abelian minimal normal subgroups.  For each prime $p$, let $S_p$ be the subgroup generated by the abelian minimal normal
subgroups that are $p$-groups.  Then
\[
\Soc(G)=B_1\times\cdots\times B_b\times\prod_p S_p,
\]
and each $S_p$ is a semisimple $\F_pG$-module.  Let $A_{p,1},\ldots,A_{p,m_p}$ be representatives of the isomorphism classes of simple submodules of $S_p$, and let $S_{p,j}$ be the sum of all simple submodules of $S_p$ isomorphic to $A_{p,j}$. Then $S_p$ has the isotypic decomposition
\[
S_p=\bigoplus_{j=1}^{m_p}S_{p,j},\qquad S_{p,j}\cong A_{p,j}^{\,n_{p,j}},
\]
where $n_{p,j}\geq1$.  Set
\[
D_{p,j}=\End_{\F_pG}(A_{p,j}).
\]
Then:
\begin{enumerate}
\item the minimal normal subgroups of $G$ are the groups $B_i$ together with
      the simple submodules of the homogeneous components $S_{p,j}$;
\item their number is
      \[
      b+\sum_p\sum_{j=1}^{m_p}
      \frac{|D_{p,j}|^{n_{p,j}}-1}{|D_{p,j}|-1};
      \]
\item the socle length is
      \[
      \sigma(G)=b+\sum_p\sum_{j=1}^{m_p}n_{p,j}.
      \]
\end{enumerate}
\end{proposition}

\begin{proof}
     Let $N\in\MinNorm(G)$ be abelian and let $p$ be a prime dividing $|N|$. The Sylow $p$-subgroup of $N$ and $\{x\in N:x^{p}=1\}$ are non-trivial characteristic subgroups of $N$, and hence normal in $G$. So both are equal to $N$, that is, $N$ is an elementary abelian $p$-group. Thus every minimal normal subgroup of $G$ is some $B_{i}$ or lies in some $S_{p}$, and hence $\Soc(G)=UB_{1}\cdots B_{b}$, where $U=\prod_{p}S_{p}$.\par
    By Lemma~\ref{lem:distinct-atoms}, distinct minimal normal subgroups commute elementwise. So $U$ is abelian and each $S_{p}$ is an elementary abelian $p$-group. As the $S_{p}$ have pairwise coprime orders, $U=\prod_{p}S_{p}$ is a direct product. Now $Z(B_{i})$ is characteristic in $B_{i}$, hence normal in $G$, and $Z(B_{i})\ne B_{i}$ as $B_{i}$ is non-abelian. So $Z(B_{i})=1$. Again by Lemma~\ref{lem:distinct-atoms}, $UB_{1}\cdots B_{i-1}$ centralizes $B_{i}$, and so
    \[
       B_{i}\cap UB_{1}\cdots B_{i-1}\leq Z(B_{i})=1.
    \]
    As all these subgroups are normal in $G$, it follows that $\Soc(G)=B_{1}\times\cdots\times B_{b}\times\prod_{p}S_{p}$.\par
    Now $S_{p}\trnl G$ and $S_{p}$ is elementary abelian, so $G$ acts $\F_{p}$-linearly on $S_{p}$ by conjugation, and $S_{p}$ is an $\F_{p}G$-module. As every subgroup of $S_{p}$ is an $\F_{p}$-subspace, the normal subgroups of $G$ contained in $S_{p}$ are exactly the $\F_{p}G$-submodules of $S_{p}$. In particular, the abelian minimal normal $p$-subgroups of $G$ are exactly the simple submodules of $S_{p}$. As $S_{p}$ is generated by them, $S_{p}$ is a sum of simple submodules, and hence it is semisimple.\par
    So let $S_{p}=T_{1}\oplus\cdots\oplus T_{s}$ with simple submodules $T_{k}$. If $T$ is a simple submodule of $S_{p}$ with $T\cong A_{p,j}$, then by Schur's Lemma the projection of $T$ onto $T_{k}$ is zero whenever $T_{k}\not\cong A_{p,j}$. So $T$ lies in the sum of those $T_{k}$ with $T_{k}\cong A_{p,j}$. Hence
    \[
       S_{p,j}=\bigoplus_{T_{k}\cong A_{p,j}}T_{k}\cong A_{p,j}^{\,n_{p,j}},\quad\text{where }n_{p,j}=\#\{k:T_{k}\cong A_{p,j}\}
    \]
    and $n_{p,j}\geq1$, as some simple submodule of $S_{p}$ is isomorphic to $A_{p,j}$. Also, each $T_{k}$ is isomorphic to exactly one $A_{p,j}$, so $S_{p}=\bigoplus_{j}S_{p,j}$. By the definition of $S_{p,j}$ and as this sum is direct, every simple submodule of $S_{p}$ lies in exactly one $S_{p,j}$. This proves (1).\par
    By Lemma~\ref{lem:projective-submodules} with $R=\F_{p}G$, the module $S_{p,j}\cong A_{p,j}^{\,n_{p,j}}$ has exactly $$\frac{|D_{p,j}|^{n_{p,j}}-1}{|D_{p,j}|-1}$$ simple submodules. Summing over all $p$ and $j$ and adding the $b$ non-abelian ones gives (2). Finally, the $B_{i}$ and the $T_{k}$ (for all $p$) exhibit $\Soc(G)$ as a direct product $X_{1}\times\cdots\times X_{t}$ of $t=b+\sum_{p}\sum_{j}n_{p,j}$ minimal normal subgroups. As $X_{l}\cap X_{1}\cdots X_{l-1}=1$ and $\lambda_{G}(X_{l})=1$, \eqref{eq:modular-rank} of Lemma~\ref{lem:relative-length} gives $$\lambda_{G}(X_{1}\cdots X_{l})=\lambda_{G}(X_{1}\cdots X_{l-1})+1.$$ Hence $\sigma(G)=\lambda_{G}(\Soc(G))=t$, which is (3).
\end{proof}

We can now prove Theorem~B. The point is that the combinatorial condition of
Theorem~\ref{thm:reconstruction}(4) is equivalent to the numerical condition
$\sigma(G)\leq2$, and that Proposition~\ref{prop:socle-geometry} turns the
latter into a short list of socle types.

\begin{theorem}\label{thm:Mfull-classification}
Let $G$ be a non-trivial finite group and put $S=\Soc(G)$.  The following are
equivalent:
\begin{enumerate}
\item $G$ is $\M$-full;
\item $AB=S$ for every two distinct $A,B\in\MinNorm(G)$;
\item $\sigma(G)\leq2$;
\item one of the following mutually exclusive alternatives holds:
      \begin{enumerate}[label=\upshape(\alph*)]
      \item $G$ has exactly one minimal normal subgroup;
      \item $G$ has exactly two minimal normal subgroups;
      \item $G$ has at least three minimal normal subgroups, and there are a
            prime $p$ and an irreducible $\F_pG$-module $A$ such that
            $S\cong A\oplus A$ as $G$-modules.
      \end{enumerate}
\end{enumerate}
In case \textup{(4c)}, if $D=\End_{\F_pG}(A)$, then
\[
\MinNorm(G)\cong\mathds P^1(D)
\quad\text{and}\quad
|\MinNorm(G)|=|D|+1.
\]
\end{theorem}

\begin{proof} 
    (1)$\Leftrightarrow$(2) By Theorem~\ref{thm:reconstruction}(4).\par
    (2)$\Rightarrow$(3) Assume that $AB=S$, for any two distinct $A,B\in\MinNorm(G)$. Suppose, $\sigma(G)\geq 3$. Then, there is a chain
    \[
       1<N<M<S
    \]
    where $N,M\trnl G$. As $G$ is finite and hence normally atomic, there exists $A\in\MinNorm(G)$ such that $A\leq N$. As $M<S$ and $S$ is generated by the minimal normal subgroups of $G$, there exists $B\in\MinNorm(G)$ such that $B\nleq M$. Then $A\ne B$, as $A\leq N<M$. So by our hypothesis, $AB=S$. Now, $M\cap B$ is a normal subgroup of $G$ properly contained in $B$, and hence $M\cap B=1$. As $A\leq M$, by Lemma~\ref{lem:modular},
    \[
       M=AB\cap M=A(B\cap M)=A\leq N,
    \]
    which is a contradiction as $N<M$. Thus, $\sigma(G)\leq 2$.\par
    (3)$\Rightarrow$(4) Assume $\sigma(G)\leq 2$. As $G\ne\{1\}$, $S\ne\{1\}$ and $\sigma(G)\geq1$. We use the notation of Proposition~\ref{prop:socle-geometry}. First suppose $n_{p,j}=1$ for all $p$ and $j$. Then each term $$\frac{|D_{p,j}|^{n_{p,j}}-1}{|D_{p,j}|-1}$$ is $1$, and so by Proposition~\ref{prop:socle-geometry}(2),(3),
    \[
       |\MinNorm(G)|=b+\sum_{p}\sum_{j=1}^{m_{p}}n_{p,j}=\sigma(G)\in\{1,2\}.
    \]
    So, $G$ has exactly one or exactly two minimal normal subgroups. Now suppose $n_{p,j}\geq2$ for some $p$ and $j$. As $\sigma(G)=b+\sum_{p}\sum_{j}n_{p,j}\leq2$, $b=0$, $n_{p,j}=2$, and $S_{p,j}$ is the only homogeneous component. Thus, $S=S_{p}=S_{p,j}\cong A\oplus A$, where $A=A_{p,j}$ is an irreducible $\F_{p}G$-module. By Proposition~\ref{prop:socle-geometry}(2), with $D=\End_{\F_{p}G}(A)$,
    \[
       |\MinNorm(G)|=\frac{|D|^{2}-1}{|D|-1}=|D|+1\geq3
    \]
    as $|D|\geq2$. So (c) holds. Clearly, the alternatives (a), (b), (c) are mutually exclusive.\par
    (4)$\Rightarrow$(2) In case (a), there is nothing to prove. In case (b), if $M_{1},M_{2}$ are the two minimal normal subgroups of $G$, then $S=\langle M_{1},M_{2}\rangle=M_{1}M_{2}$. In case (c), $S\cong A\oplus A$ is an elementary abelian $p$-group and every subgroup of $S$ is an $\F_{p}$-subspace. So the normal subgroups of $G$ contained in $S$ are exactly the $\F_{p}G$-submodules of $S$. As every minimal normal subgroup of $G$ lies in $S$, the minimal normal subgroups of $G$ are exactly the simple submodules of $S\cong A^{2}$. Hence by Lemma~\ref{lem:projective-submodules}, \[\MinNorm(G)\cong\mathds{P}^{1}(D),\qquad |\MinNorm(G)|=|D|+1,\] and any two distinct $M_{1},M_{2}\in\MinNorm(G)$ satisfy $M_{1}\oplus M_{2}=S$, that is, $M_{1}M_{2}=S$.
\end{proof}
Whether a group is $\M$-full depends only on its normal subgroups, not on its size or on its other subgroups. For example, for every $n\ge 2$, the group $S_n$ is $\M$-full (Example~\ref{ex:D60-Sn}); thus $S_5$, which has $156$ subgroups, is $\M$-full because it has only three normal subgroups. The following example shows that each of the three alternatives in Theorem \ref{thm:Mfull-classification}(4) actually occurs, already inside one family of groups.

\begin{example}\label{ex:SL2}
	For $n\ge2$ let $\mathrm{SL}_2(\Z_n)$ be the group of $2\times2$
	matrices over $\Z_n$ with determinant $1$. 
 The three cases of Theorem~\ref{thm:Mfull-classification}(4) all occur.
	\begin{enumerate}
		\item For $n=3,5$, the group $\mathrm{SL}_2(\Z_n)$ is monolithic with minimal normal subgroup $\{I,-I\}$. Indeed, the normal subgroups of $\mathrm{SL}_2(\Z_3)$ are $1$, $\{I,-I\}$, the quaternion subgroup of order $8$ and $\mathrm{SL}_2(\Z_3)$, while those of $\mathrm{SL}_2(\Z_5)$ are $1$, $\{I,-I\}$ and $\mathrm{SL}_2(\Z_5)$, because $\mathrm{SL}_2(\Z_5)/\{I,-I\}\cong A_5$ is simple and $\mathrm{SL}_2(\Z_5)$ is perfect. This is case (4a).
		\item Clearly, $\mathrm{SL}_2(\Z_6)\cong\mathrm{SL}_2(\Z_2)\times\mathrm{SL}_2(\Z_3)\cong S_3\times\mathrm{SL}_2(\Z_3)$.
		  By (1) and Theorem~\ref{thm:product-atoms},
		\[
		\Soc(S_3\times\mathrm{SL}_2(\Z_3))=A_3\times\{I,-I\}.
		\]
		The only minimal normal subgroups are $A_3\times1$ and $1\times\{I,-I\}$: the only non-central minimal normal subgroup of a factor is $A_3\trnl S_3$, and, as $Z(S_3)=1$ and $Z(\mathrm{SL}_2(\Z_3))=\{I,-I\}$, the space $V_2$ is the line $1\times\{I,-I\}$ and $V_p=1$ for $p\neq2$. This is case (4b).
		\item  For $n=15$, we have $\mathrm{SL}_2(\Z_{15})\cong\mathrm{SL}_2(\Z_3)\times\mathrm{SL}_2(\Z_5)$.
		So by (1) and Theorem~\ref{thm:product-atoms},
		\[
		\Soc(\mathrm{SL}_2(\Z_{15}))\cong\{I,-I\}\times\{I,-I\}\cong C_2\times C_2 .
		\]
		So the minimal normal subgroups are exactly its
		three subgroups of order $2$. Inside $\mathrm{SL}_2(\Z_{15})$ the minimal normal subgroups are
		generated by $4I$, $11I$ and $14I=-I$ and they correspond to
		$(I,-I)$, $(-I,I)$ and $(-I,-I)$, since
		$4\equiv1\pmod3$, $4\equiv-1\pmod5$, and $11\equiv-1\pmod3$, $11\equiv1\pmod5$.
		Here $A=\Z_2$ is the trivial module and $D=\Z_2$, and
		$|\MinNorm(\mathrm{SL}_2(\Z_{15}))|=3=|D|+1$. This is case (4c).
	\end{enumerate}
\end{example}

For $p$-groups all minimal normal subgroups are central of order $p$, and
Theorem~\ref{thm:Mfull-classification} becomes a condition on the rank of the
centre.

\begin{corollary}\label{cor:pgroup-Mfull}
Let $G$ be a non-trivial finite $p$-group and put
\[
V=\Omega_1(Z(G))=\{z\in Z(G):z^p=1\}.
\]
Then the minimal normal subgroups of $G$ are exactly the one-dimensional
$\F_p$-subspaces of $V$, and $\Soc(G)=V$.  Moreover, $G$ is $\M$-full if and
only if $\dim_{\F_p}V\leq2$.
\end{corollary}

\begin{proof}
      As $Z(G)$ is abelian, $V$ is an $\F_{p}$-vector space. Let $N\in\MinNorm(G)$ and let $G\curvearrowright N$ by 
    \[
       g\cdot n\coloneqq gng^{-1}.
    \]
    Then, as $G$ is a $p$-group, $N$ is also a $p$-group and hence the class equation of this action gives:
    \[
       |N|=|N^G|+\sum_{i}|\mathcal{O}_{i}|,
    \]
    where $p|\,|\mathcal{O}_{i}|$ for each $i$. So it immediately follows that $p|\,|N^G|$. But, $$N^G=\{n\in N:gng^{-1}=n,\,\forall g\in G\}=N\cap Z(G).$$ Now, $1\ne N\cap Z(G)\trnl G$ and $N\cap Z(G)\leq N$. But, $N\in\MinNorm(G)$ implies that $N=N\cap Z(G)$ and hence, $N\leq Z(G)$. So, every subgroup of $N$ is normal in $G$. Now, if $|N|>p$, then there exists $L\leq N$ such that $1<L<N$ and $|L|=p$. But, $L\trnl G$ --- which is a contradiction as $N\in\MinNorm(G)$. So, $|N|=p$. Thus, every minimal normal subgroup of $G$ is cyclic of order $p$ and hence a one-dimensional subspace of $V$. So, $N\leq V$, and we get
    \[
       \MinNorm(G)\subseteq\{1\text{-dimensional }\F_{p}\text{-}\text{subspaces}\ \text{ of }V\}.
    \]
    Conversely, if $L\leq V$ is a one-dimensional subspace of $V$, then $|L|=p$. As $L\leq V\leq Z(G)$, we have $L\trnl G$. As a group of order $p$ has no non-trivial proper subgroups, it follows that $L\in\MinNorm(G)$. Hence,
    \[
       \{1\text{-dimensional }\F_{p}\text{-}\text{subspaces}\text{ of }V\}\subseteq\MinNorm(G).
    \]
    Thus, the minimal normal subgroups of $G$ are exactly the one-dimensional $\F_{p}$-subspaces of $V$.\par
    Now, as $M\in\MinNorm(G)$ implies that $M\leq V$, it follows at once, $\Soc(G)\leq V$. But, $V$ is an $\F_{p}$-vector space, so it is generated by its one-dimensional $\F_{p}$-subspaces. Hence, $V\leq\Soc(G)$. Thus, we get $\Soc(G)=V$.\par
    Now, observe that the interval $[1,V]_{\N(G)}$ is the interval of $\F_{p}$-subspaces of $V$. So,
    \[
       [1,V]_{\N(G)}=\text{SubSp}_{\F_{p}}(V),
    \]
    where $\text{SubSp}_{\F_{p}}(V)$ denotes the lattice of $\F_{p}$-subspaces of $V$. But, a maximal chain in $\text{SubSp}_{\F_{p}}(V)$ looks like
    \[
       V_{0}=\{0\}<V_{1}<V_{2}<\cdots<V_{n}=V,
    \]
    where $n=\dim_{\F_{p}}V$, i.e., there are $\dim_{\F_{p}}V$ number of inclusions. Hence, $\lambda_{G}(V)=\dim_{\F_{p}}V$. But, $\Soc(G)=V$. So, $[1,V]_{\N(G)}=[1,\Soc(G)]_{\N(G)}$, and thus
    \[
       \dim_{\F_{p}}V=\lambda_{G}(V)=\lambda_{G}(\Soc(G))=\sigma(G).
    \]
    But, by Theorem~\ref{thm:Mfull-classification}, we know that 
    \[
       G\text{ is }\M\text{-full} \iff \sigma(G)\leq 2.
    \]
    Thus, $G$ is $\M$-full if and only if $\dim_{\F_{p}}V\leq 2$.

\end{proof}

Since a finite nilpotent group is the direct product of its Sylow subgroups,
the same argument applies prime by prime.

\begin{corollary}\label{cor:nilpotent-Mfull}
Let $G$ be a non-trivial finite nilpotent group, let $P_p$ denote its Sylow
$p$-subgroup, and put
\[
 \rho(G)=\sum_{p\mid |G|}
 \dim_{\F_p}\{z\in Z(P_p):z^p=1\}.
\]
Then $G$ is $\M$-full if and only if $\rho(G)\leq2$.
\end{corollary}

\begin{proof}
    Let $G$ be a finite nilpotent group. Then,
    \[
       G=\prod_{p \mid|G|}P_{p},
    \]
    where each $P_{p}$ is the unique Sylow $p$-subgroup of $G$. Then, $$Z(G)=Z\left(\displaystyle{\prod_{p||G|}}P_{p}\right)=\displaystyle{\prod_{p||G|}}Z(P_{p}).$$ Now, as each $P_{p}$ is a finite $p$-group, $Z(P_{p})$ is a finite abelian $p$-group and hence, $\Omega_{1}(Z(P_{p}))$ is an elementary abelian $p$-group. Thus, for each $p|\,|G|$, $\Omega_{1}(Z(P_{p}))$ is an $\F_{p}$-vector space. Now, let $N\in\MinNorm(G)$. Then,
    \[
       N=\prod_{p||G|}(N\cap P_{p}),
    \]
    that is, every element of $N$ is the product of its $p$-parts, which are powers of it and so lie in $N\cap P_{p}$, where $N\cap P_{p}\trnl G$. As $N\ne 1$, there exists at least one $N\cap P_{p}\ne\{1\}$. But, $N\cap P_{p}\leq N$. So by minimality of $N$, we get $N\cap P_{p}=N$ and thus, $N\leq P_{p}$. As the other Sylow subgroups of $G$ centralize $P_{p}$, every normal subgroup of $P_{p}$ is normal in $G$, and so $N$ is a minimal normal subgroup of the finite $p$-group $P_{p}$. So, by Corollary~\ref{cor:pgroup-Mfull} we get
    \[
       N\leq\Omega_{1}(Z(P_{p})) \quad\text{and}\quad |N|=p.
    \]
    Hence, $N$ is a one-dimensional $\F_{p}$-subspace of $\Omega_{1}(Z(P_{p}))$.\par
    Conversely, suppose $L\leq\Omega_{1}(Z(P_{p}))$ is a one-dimensional subspace. Then, $|L|=p$. Now,
    \[
       L\leq\Omega_{1}(Z(P_{p}))\leq Z(P_{p})\leq Z(G)
    \]
    implies that $L\trnl G$ and as $|L|=p$, it follows $L\in\MinNorm(G)$. So, we have a correspondence,
    \[
       \MinNorm(G)=\coprod_{p||G|}\mathds{P}(\Omega_{1}(Z(P_{p}))),
    \]
    where, $\mathds{P}(\Omega_{1}(Z(P_{p})))$ is the collection of all one-dimensional $\F_{p}$-subspaces of $\Omega_{1}(Z(P_{p}))$. Hence, we get
    \[
       \Soc(G)=\prod_{p||G|}\Omega_{1}(Z(P_{p})).
    \]
    As $\Omega_{1}(Z(P_{p}))\leq Z(G)$, it follows that every subgroup of $\Omega_{1}(Z(P_{p}))$ is normal in $G$ and thus like before we get
    \[
       [1,\Omega_{1}(Z(P_{p}))]_{\N(G)}=\text{SubSp}_{\F_{p}}(\Omega_{1}(Z(P_{p}))).
    \]
    Hence, $$\lambda_{G}(\Omega_{1}(Z(P_{p})))=\dim_{\F_{p}}\Omega_{1}(Z(P_{p})).$$ Now, $\Soc(G)=\displaystyle{\prod_{p||G|}}\Omega_{1}(Z(P_{p}))$ and as each $\Omega_{1}(Z(P_{p}))$ has $p$-power order, these factors have pairwise coprime orders, and every subgroup of $\Soc(G)\leq Z(G)$ is normal in $G$, so we get
    \[
       [1,\Soc(G)]_{\N(G)} \cong\prod_{p||G|}[1,\Omega_{1}(Z(P_{p}))]_{\N(G)}.
    \]
    Therefore, we get
    \begin{align*}
        \sigma(G)&=\lambda_{G}(\Soc(G)) \\
        &=\sum_{p||G|}\lambda_{G}(\Omega_{1}(Z(P_{p}))) \\
        &=\sum_{p||G|}\dim_{\F_{p}}\Omega_{1}(Z(P_{p})) \\
        &=\sum_{p||G|}\dim_{\F_{p}}\{z\in Z(P_{p}):z^{p}=1\} \\
        &=\rho(G).
    \end{align*}
     By Theorem~\ref{thm:Mfull-classification}, we know that 
    \[
       G\text{ is }\M\text{-full} \iff \sigma(G)\leq 2.
    \]
    Hence, we conclude that $G$ is $\M$-full if and only if $\rho(G)\leq 2$.
\end{proof}

We close this section by showing that the equivalence $(1)\Leftrightarrow(3)$ of Theorem~\ref{thm:Mfull-classification} has a purely lattice-theoretic form, which holds for arbitrary groups.

Let $G$ be a fixed group. Then $\N(G)$, the set of normal subgroups of $G$, is a complete modular lattice by Lemma~\ref{lem:modular}, with meet the intersection and join the subgroup generated. Recall that a \emph{quantale} is a complete lattice $Q$ together with an associative multiplication $\cdot$ which distributes over arbitrary joins in both variables, that is,
\[
  a\cdot\bigvee_{i}b_{i}=\bigvee_{i}(a\cdot b_{i})
  \qquad\text{and}\qquad
  \Bigl(\bigvee_{i}a_{i}\Bigr)\cdot b=\bigvee_{i}(a_{i}\cdot b)
\]
for all $a,b,a_{i},b_{i}\in Q$. A quantale is \emph{commutative} if $\cdot$ is commutative, and \emph{integral} if the top element of $Q$ is an identity for $\cdot$. The $\mu$-elements of \cite{CavigliaGoswamiMesiti} are defined in integral commutative quantales, but their definition only uses the meet, and the notions below only use the meet and the join. So they make sense in $\N(G)$, although $\N(G)$ is in general not an integral commutative quantale: the commutator $[N,M]$ is not associative (for $G=S_{3}$, $[[G,G],A_{3}]=1\ne A_{3}=[G,[G,A_{3}]]$), the product $NM$ is just the join, for which $G$ is not an identity, and $\cap$ does not distribute over joins (for $G=C_{2}\times C_{2}$). In this language, the $\M$-subgroups of $G$ are precisely the $\mu$-elements of $\N(G)$. Following \cite{CavigliaGoswamiMesiti, GrzeszczukGoldie}, we use the following notion.

\begin{definition}
A family of non-trivial normal subgroups $N_1,\dots,N_k$ of $G$ is called \emph{independent}
(also called \emph{join-independent}) if, for every $i$,
\[
  N_i\cap \prod_{j\ne i}N_j=1.
\]
\end{definition}

As the lattice $\N(G)$ is modular, following \cite{GrzeszczukGoldie}, in the terminology of Goldie and dual Goldie dimensions, we make the following definition.

\begin{definition}
The \emph{uniform dimension} of $\N(G)$ is
\[
  \mathrm{udim}\, \N(G)\in \{0,1,2,\dots\}\cup\{\infty\},
\]
defined to be the supremum of the cardinalities of independent
families of non-trivial normal subgroups of $G$.
\end{definition}

We say that $\N(G)$ is \emph{uniform} if $\mathrm{udim}\,\N(G)\le 1$, that is, every non-trivial normal subgroup of $G$ is essential in $G$.

\begin{theorem}\label{thm:Mfull-udim}
The following conditions are equivalent for a group $G$.
\begin{enumerate}
  \item $G$ is $\M$-full, i.e., $\N^\M(G)=\N(G)$.
  \item There are no non-trivial normal subgroups $A,B,C$ of $G$ such
  that
  \[
    A\cap B=1
    \qquad\text{and}\qquad
    AB\cap C=1.
  \]
  \item $G$ contains no internal direct product $A\times B\times C$ of
  three non-trivial normal subgroups of $G$.
  \item $\mathrm{udim}\,\N(G)\le 2$.
\end{enumerate}
\end{theorem}

\begin{proof}
We begin with $(1)\Rightarrow(2)$, arguing contrapositively. Suppose that
$A,B,C$ are non-trivial normal subgroups of $G$ such that
\[
  A\cap B=1
  \qquad\text{and}\qquad
  AB\cap C=1.
\]
Put $E=AB$. Since $A$ and $B$ are normal, $E$ is again a normal
subgroup of $G$. Moreover,
\[
  A\le E,\qquad B\le E,\qquad E\cap C=1.
\]

Consider the normal subgroups $AC$ and $BC$. Since $C\le AC\cap BC$,
we have $AC\cap BC\ne 1$. Also,
\[
  A\le AC\cap E
  \qquad\text{and}\qquad
  B\le BC\cap E,
\]
so $AC\cap E\ne 1$ and $BC\cap E\ne 1$.

Finally, in fact $AC\cap BC=C$. Indeed, if $x\in AC\cap BC$, write
$x=ac=bc'$ with $a\in A$, $b\in B$, and $c,c'\in C$. Then
\[
  a^{-1}b=cc'^{-1}\in AB\cap C=1,
\]
whence $a=b\in A\cap B=1$, and therefore $x=c=c'\in C$. Consequently,
\[
  E\cap(AC\cap BC)=E\cap C=1.
\]
Thus $E$ is not an $\M$-subgroup, and hence $G$ is not $\M$-full. This proves $(1)\Rightarrow(2)$ by contraposition.\par For $(2)\Rightarrow(1)$, suppose that $E\trnl G$ is not an $\M$-subgroup. Then by Theorem~\ref{thm:dichotomy}, $E$ is neither essential nor relatively uniform in $G$. So there are non-trivial normal subgroups $A,B,C$ of $G$ such that $A,B\leq E$, $A\cap B=1$ and $E\cap C=1$. As $AB\leq E$, $AB\cap C\leq E\cap C=1$, and hence $(2)$ fails.\par
For $(2)\Leftrightarrow(3)$, observe that for non-trivial normal subgroups $A,B,C$ of $G$, the product $ABC$ is the internal direct product $A\times B\times C$ if and only if $A\cap B=1$ and $AB\cap C=1$.\par
For $(2)\Leftrightarrow(4)$, observe that a subfamily of an independent family is again independent. So $\mathrm{udim}\,\N(G)\geq3$ if and only if there are three independent non-trivial normal subgroups $A,B,C$ of $G$. Now if $A,B,C$ are independent, then $A\cap B\leq A\cap BC=1$ and $AB\cap C=1$. Conversely, let $A\cap B=1$ and $AB\cap C=1$. If $x=bc\in A\cap BC$ with $b\in B$ and $c\in C$, then $c=b^{-1}x\in AB\cap C=1$, and so $x=b\in A\cap B=1$. Thus $A\cap BC=1$, and similarly $B\cap AC=1$. So $A,B,C$ are independent. Hence $(2)$ holds if and only if $\mathrm{udim}\,\N(G)\leq2$.
\end{proof}

\begin{corollary}\label{cor:finiteudim}
  If $G$ is finite, then $\mathrm{udim}\,\N(G)=\sigma(G)$. Hence, a finite group is
  $\M$-full if and only if $\sigma(G)\le 2$.
\end{corollary}

\begin{proof}
If $G=1$, then both sides are $0$. So let $G\ne1$. By Proposition~\ref{prop:socle-geometry}, $\Soc(G)=B_{1}\times\cdots\times B_{b}\times\prod_{p}S_{p}$, where $S_{p}=\bigoplus_{j}S_{p,j}$ and each $S_{p,j}\cong A_{p,j}^{\,n_{p,j}}$ is a direct sum of $n_{p,j}$ simple submodules, which are minimal normal subgroups of $G$ by Proposition~\ref{prop:socle-geometry}$(1)$. So by Proposition~\ref{prop:socle-geometry}$(3)$, there are minimal normal subgroups $N_{1},\dots,N_{\sigma(G)}$ of $G$ such that
\[
  \Soc(G)=N_1\times\cdots\times N_{\sigma(G)}.
\]
The factors in a direct product are independent, that is,
$N_i\cap\prod_{j\neq i}N_j=1$ for every $i$. Thus $\{N_1,\dots,N_{\sigma(G)}\}$ is an
independent family of non-trivial normal subgroups of $G$, and the inequality
\[
  \mathrm{udim}\,\N(G)\ge\sigma(G)
\]
is immediate.

Conversely, let $N_1,\dots,N_k$ be independent non-trivial normal subgroups of
$G$, so that $N_i\cap\prod_{j\neq i}N_j=1$ for every $i$. For each $i$, we
choose a minimal normal subgroup $M_i$ of $G$ with $M_i\le N_i$; this is
possible because $G$ is finite and $N_i\neq 1$ (for instance, take $M_i$ to be
a non-trivial normal subgroup of $G$ of least order contained in $N_i$). Since
$M_1\cdots M_{i-1}\le N_1\cdots N_{i-1}\le\prod_{j\neq i}N_j$, we clearly have
\[
  M_i\cap M_1\cdots M_{i-1}\le N_i\cap\prod_{j\neq i}N_j=1.
\]
As $M_i\neq 1$, it follows that $M_i\not\le M_1\cdots M_{i-1}$, and hence
$M_1\cdots M_{i-1}<M_1\cdots M_i$. Moreover, each $M_i$, being a minimal normal
subgroup of $G$, is contained in $\Soc(G)$. Therefore
\[
  1<M_1<M_1M_2<\cdots<M_1M_2\cdots M_k\le\Soc(G)
\]
is a chain of length $k$ in $[1,\Soc(G)]_{\N(G)}$, and so $k\leq\lambda_{G}(\Soc(G))=\sigma(G)$.
As $\{N_1,\dots,N_k\}$ was an arbitrary independent family of non-trivial
normal subgroups of $G$, we obtain
\[
  \mathrm{udim}\,\N(G)\le\sigma(G).
\]

Hence $\mathrm{udim}\,\N(G)=\sigma(G)$. Since a group is $\M$-full if and only if $\mathrm{udim}\,\N(G)\le 2$ (\textit{cf.} Theorem~\ref{thm:Mfull-udim}), the last assertion follows.
\end{proof}

%-----------------------------------------------------------------------------------------------------
\section{Direct products and diagonal atoms}
\label{sec:products}

In a direct product the normal subgroups are not, in general, products of
normal subgroups of the factors, and the minimal normal subgroups need not be
coordinate subgroups. The deviation is caused entirely by central elements of
prime order, which produce diagonal minimal normal subgroups. The following
theorem, the first part of Theorem~C, makes this precise.

\begin{theorem}\label{thm:product-atoms}
Let $G=G_1\times\cdots\times G_t$ be a finite direct product.  Its minimal
normal subgroups are precisely
\begin{enumerate}
\item the coordinate subgroups $\widehat A$, where
      $A\in\MinNorm(G_i)$ and $A\nleq Z(G_i)$;
\item the one-dimensional $\F_p$-subspaces of $V_p$, as $p$ ranges over the
      prime divisors of $|G|$.
\end{enumerate}
Consequently,
\[
\Soc(G)=\Soc(G_1)\times\cdots\times\Soc(G_t).
\]
A line in $V_p$ is non-coordinate precisely when its projection is non-zero
in at least two factors.
\end{theorem}
\begin{proof}
    Suppose, $N\in\MinNorm(G)$. Then, we have two cases:\par
    Case 1: $[N,G_{i}^{*}]\ne1$ for some $i$. As, $N\trnl G$ and $G_{i}^{*}\trnl G$, so $[N,G_{i}^{*}]\leq N\cap G_{i}^{*}$. Hence, $N\cap G_{i}^{*}\ne 1$. But $N\cap G_{i}^{*}\trnl G$, and so by minimality of $N$, we get $N\cap G_{i}^{*}=N$, and hence $N\leq G_{i}^{*}$. Now, let $A\leq G_{i}$ such that $N=\hat{A}$. For any $g_{i}\in G_{i}$ and $a\in A$, we have
    \[
       \hat{g_{i}}\hat{a}\hat{g_{i}}^{-1}\in N=\hat{A}
    \]
    as $N\trnl G$. Thus, $g_{i}ag_{i}^{-1}\in A$. Hence, $A\trnl G_{i}$. Now, let $1\ne C\trnl G_{i}$ such that $C\leq A$. Then, $\hat{C}\trnl G$, as $C\trnl G_{i}$ and $\hat{C}$ is centralized by $G_{j}^{*}$ for $j\ne i$, and $\hat{C}\leq \hat{A}=N$. So, by minimality of $N$, we get $\hat{C}=N=\hat{A}$ and hence, $C=A$. Thus, $A\in\MinNorm(G_{i})$. As $[N,G_{i}^{*}]\ne1$, it implies that $[\hat{A},G_{i}^{*}]\ne1$ and hence $[A,G_{i}]\ne 1$. So, $A\nleq Z(G_{i})$.\par
    Case 2: $[N,G_{i}^{*}]=1$ for all $i$. As $G=G_{1}^{*}G_{2}^{*}\cdots G_{t}^{*}$ and $[N,G_{i}^{*}]=1$ for all $i$, it follows that $[N,G]=1$. Hence, $N\leq Z(G)$. So, $N$ is a finite abelian group. As $N$ is minimal normal, it follows that $|N|=p$ for some prime $p$. Hence, $N$ is a one-dimensional $\F_{p}$-vector space. Now, $Z(G)=\displaystyle{\prod_{i=1}^{t}}Z(G_{i})$ and $N\leq Z(G)$, $|N|=p$. So, it follows that
    \[
       N\leq\prod_{i=1}^{t} V_{i,p}=V_{p}.
    \]
    So, $N$ is a one-dimensional $\F_{p}$-subspace of $V_{p}$.\par
    Conversely, suppose that $A\in\MinNorm(G_{i})$ such that $A\nleq Z(G_{i})$. Then, $A\leq G_{i}$ and hence $\hat{A}\leq G_{i}^{*}$. Now, let $g\in G$ and $\hat{a}\in\hat{A}$. Then, as $A\trnl G_{i}$, we get $g\hat{a}g^{-1}\in\hat{A}$. Thus, $\hat{A}\trnl G$. Now, let $1\ne N\trnl G$ be a normal subgroup such that $N\leq\hat{A}$. Then $N=\hat{C}$ for some $1\ne C\leq A$, and $C\trnl G_{i}$ as $N\trnl G$. But as $A\in\MinNorm(G_{i})$, we must have $C=A$ and hence, $N=\hat{A}$. Thus, $\hat{A}\in\MinNorm(G)$.\par
    Now suppose $L\leq V_{p}$ is a one-dimensional $\F_{p}$-subspace. Then, $|L|=p$ and $L\leq V_{p}\leq Z(G)$. So, $L\trnl G$, and as a group of prime order, it has no proper non-trivial subgroups, it follows that $L\in\MinNorm(G)$. Thus, the minimal normal subgroups of $G$ are precisely the coordinate subgroups $\hat{A}$, where $A\in\MinNorm(G_{i})$ and $A\nleq Z(G_{i})$ along with the one-dimensional $\F_{p}$-subspaces of $V_{p}$, where $p$ varies over the prime divisors of $|G|$.\par
    Now, for each $i$, we have
    \[
       \Soc(G_{i})=\langle A :A\in\MinNorm(G_{i})\rangle.
    \]
    So, the coordinate copy of $\Soc(G_{i})$, denoted $\Soc^*(G_{i})$, is generated by the coordinate minimal normal subgroups, $\hat{A}$ where $\hat{A}\in\MinNorm(G_{i}^*)$. But, we just now proved that every $A\in\MinNorm(G_{i})$ gives us $\hat{A}\in\MinNorm(G)$, and the argument above does not use $A\nleq Z(G_{i})$. Therefore,
    \[
       \Soc^*(G_{1})\Soc^*(G_{2})\cdots\Soc^*(G_{t})\leq\Soc(G).
    \]
    As $\Soc^*(G_{i})$ is the canonical subgroup, we hence get
    \[
       \Soc(G_{1})\times\cdots\times\Soc(G_{t})\leq\Soc(G).
    \]
    Now, for any $N\in\MinNorm(G)$, either $N=\hat{A}$ for some $A\in\MinNorm(G_{i}),A\nleq Z(G_{i})$ or $N\leq V_{p}$ is a one-dimensional $\F_{p}$-subspace of $V_{p}$. Now, if $N=\hat{A}$, then $N\in\MinNorm(G_{i}^*)$. So, if we canonically identify, $G_{i}$ with $G_{i}^*$, we then get $N\in\MinNorm(G_{i})$, and hence $N\leq\Soc^{*}(G_{i})\leq\prod_{j=1}^{t}\Soc(G_{j})$. Now suppose, $N\leq V_{p}$ is a one-dimensional $\F_{p}$-subspace. Recall that, 
    \[
       V_{p}=\prod_{i=1}^{t}V_{i,p}.
    \]
    So, for any $1\ne x\in V_{i,p}$, $\langle x\rangle$ has order $p$. Also, as $\langle x\rangle\leq V_{i,p}\leq Z(G_{i})$, we get $\langle x\rangle\trnl G_{i}$. So, $\langle x\rangle\in\MinNorm(G_{i})$. Hence, we get $\langle x\rangle\leq\Soc(G_{i})$ and thus $x\in\Soc(G_{i})$. So, $V_{i,p}\leq\Soc(G_{i})$. Hence,
    \[
       N\leq V_{p}\leq\prod_{i=1}^{t}\Soc(G_{i}).
    \]
    As $N$ was an arbitrary minimal normal subgroup of $G$, it follows then
    \[
       \Soc(G)\leq\Soc(G_{1})\times\cdots\times\Soc(G_{t}).
    \]
    Thus, $\Soc(G)=\Soc(G_{1})\times\cdots\times\Soc(G_{t})$. Now, let $L$ be a line in $V_{p}$. Then, $L$ is a one-dimensional $\F_{p}$-subspace of $V_{p}$. As $V_{p}=\displaystyle{\prod_{i=1}^{t}}V_{i,p}$, we have the projections:
    \begin{align*}
        \pi_{i}:V_{p}&\rightarrow V_{i,p} \\
        (v_{1},\cdots,v_{t})&\longmapsto v_{i}.
    \end{align*}
    Now, suppose, $0\ne v=(v_{1},\cdots,v_{t})\in L$. As $L=\F_{p}v$, so
    \[
       L\leq\widehat{V_{i,p}} \iff \pi_{i}(v)\ne0\text{ and }\pi_{j}(v)=0\text{ for all }j\ne i.
    \]
    Thus, $L$ is coordinate in $V_{p}$ if and only if exactly one of its projections is non-zero. Hence, as $v\ne0$, we get that $L$ is non-coordinate in $V_{p}$ if and only if at least two distinct projections are non-zero.
\end{proof}
As a consequence, the number of minimal normal subgroups below a given normal
subgroup of a direct product, and hence its $\M$-status, can be computed from
the factors and from linear algebra over the prime fields.

\begin{corollary}\label{cor:product-count}
Let $G_1,\ldots,G_t$ be finite groups, put $G=G_1\times\cdots\times G_t$, and let $H\trnl G$.  Define
\[
a(H)=\#\{(i,A):A\in\MinNorm(G_i),\ A\nleq Z(G_i),\ \widehat A\leq H\}
\]
and
\[
d_p(H)=\dim_{\F_p}(H\cap V_p).
\]
Then the number of minimal normal subgroups of $G$ contained in $H$ is
\[
 a(H)+\sum_p\frac{p^{d_p(H)}-1}{p-1}.
\]
In particular, $H$ is an $\M$-subgroup if and only if $H=1$, the displayed
number is $1$, or
\(\Soc(G_1)\times\cdots\times\Soc(G_t)\leq H\).
\end{corollary}

\begin{proof}
    Suppose, $N\in\MinNorm(G)$ such that $N\leq H$. Then by Theorem~\ref{thm:product-atoms}, $N=\hat{A}$ for some $A\in\MinNorm(G_{i})$ such that $A\nleq Z(G_{i})$ or $N\leq V_{p}$ is a one-dimensional $\F_{p}$-subspace. Now, if $N=\hat{A}$, then $(i,A)\in$ $\{(i,A):A\in\MinNorm(G_{i}),\ A\nleq Z(G_{i}),\ \hat{A}\leq H\}$.\par
    Now, if $N\leq V_{p}$ is a one-dimensional $\F_{p}$-subspace then, $N\leq H\cap V_{p}$ is a one-dimensional $\F_{p}$-subspace as $H\cap V_{p}$ is an $\F_{p}$-subspace of $V_{p}$. Conversely, by Theorem~\ref{thm:product-atoms}, every such $\hat{A}$ and every one-dimensional $\F_{p}$-subspace of $H\cap V_{p}$ is a minimal normal subgroup of $G$ contained in $H$. Moreover, distinct pairs $(i,A)$ give distinct $\hat{A}$, as $\hat{A}\leq G_{i}^{*}$ and $G_{i}^{*}\cap G_{j}^{*}=1$ for $i\ne j$. Also, no $\hat{A}$ lies in any $V_{p}$, as $\hat{A}\nleq Z(G)$ while $V_{p}\leq Z(G)$, and the subspaces of $V_{p}$ for different primes $p$ are different. So, no minimal normal subgroup is counted twice. Now, the size of the set $\{(i,A):A\in\MinNorm(G_{i}),\ A\nleq Z(G_{i}),\ \hat{A}\leq H\}$ is given as:
    \[
      a(H)=\#\{(i,A):A\in\MinNorm(G_{i}),\ A\nleq Z(G_{i}),\ \hat{A}\leq H\}.
    \]
    Also, the number of one-dimensional $\F_{p}$-subspaces of $H\cap V_{p}$ is given by,
    \[
       \frac{p^{\dim_{\F_{p}}(H\cap V_{p})}-1}{p-1}=\frac{p^{d_{p}(H)}-1}{p-1}.
    \]
    Hence, the total number of minimal normal subgroups of $G$ contained in $H$ is
    \[
       a(H)+\sum_{p}\frac{p^{d_{p}(H)}-1}{p-1},
    \]
    where $p$ ranges over all prime divisors of $|G|$.\par
    Now, as $G$ is finite, by Corollary~\ref{cor:finite-classification}, we have
    \[
       H\text{ is an }\M\text{-subgroup} \iff H=1,\ |\supp(H)|=1,\ \Soc(G)\leq H.
    \]
    But, $$|\supp(H)|=\#\{N\in\MinNorm(G):N\leq H\}=a(H)+\displaystyle{\sum_{p}}\frac{p^{d_{p}(H)}-1}{p-1}.$$ Also, by Theorem~\ref{thm:product-atoms} we have
    $
       \Soc(G)=\Soc(G_{1})\times\cdots\times\Soc(G_{t}).
    $
    Hence, $H$ is an $\M$-subgroup if and only if $H=1$, $a(H)+\displaystyle{\sum_{p}}\frac{p^{d_{p}(H)}-1}{p-1}=1$ or $\Soc(G_{1})\times\cdots\times\Soc(G_{t})\leq H$.
\end{proof}

Combining Theorem~\ref{thm:product-atoms} with Theorem~\ref{thm:Mfull-classification}
we obtain the second assertion of Theorem~C: a non-trivial direct product is
$\M$-full only in the smallest possible case.

\begin{corollary}\label{cor:product-monolithic}
Let $G_1,\ldots,G_t$ be non-trivial finite groups.  If $t\geq2$, then
$G_1\times\cdots\times G_t$ is $\M$-full if and only if $t=2$ and both
$G_1$ and $G_2$ are monolithic.
\end{corollary}

\begin{proof}
    Suppose, $G=G_{1}\times\cdots\times G_{t}$ is $\M$-full. Then, by Theorem~\ref{thm:Mfull-classification}, $\sigma(G)\leq 2$. Also, by Theorem~\ref{thm:product-atoms}, we have
    \[
       \Soc(G)=\Soc(G_{1})\times\cdots\times\Soc(G_{t}).
    \]
    Now, for all $i\in\{1,\cdots,t\}$, let $A_{i}\in\MinNorm(G_{i})$. Then, $\hat{A_{i}}\leq G_{i}^*$ and $\hat{A_{i}}\trnl G$. Now, suppose without loss of generality $i<j$. If $\hat{A_{i}}=\{(1,\cdots,a_{i},\cdots,1):a_{i}\in A_{i}\}$ and $\hat{A_{j}}=\{(1,\cdots,a_{j},\cdots,1):a_{j}\in A_{j}\}$, then,
    \[
       (1,\cdots,a_{i},\cdots,1)(1,\cdots,a_{j},\cdots,1)=(1,\cdots,a_{i},\cdots,a_{j},\cdots,1)=(1,\cdots,a_{j},\cdots,1)(1,\cdots,a_{i},\cdots,1).
    \]
    So, each distinct coordinate copy commutes. Hence, the subgroups, $\hat{A_{1}}$, $\hat{A_{1}}\hat{A_{2}}$, $\cdots$, $\hat{A_{1}}\hat{A_{2}}\cdots\hat{A_{t}}$ give us a chain,
    \[
       1<\hat{A_{1}}<\hat{A_{1}}\hat{A_{2}}<\cdots<\hat{A_{1}}\hat{A_{2}}\cdots\hat{A_{t}}\leq\Soc(G).
    \]
    This is because $A_{i}\leq\Soc(G_{i})$. So identifying with their coordinate copies we get, $\hat{A_{i}}\leq\Soc^*(G_{i})$. So, $\hat{A_{i}}\leq\Soc(G)$. Each inclusion is strict, as $1\ne\hat{A_{k+1}}\leq G_{k+1}^{*}$, while $\hat{A_{1}}\cdots\hat{A_{k}}\leq G_{1}^{*}\cdots G_{k}^{*}$ meets $G_{k+1}^{*}$ trivially. Thus, $\sigma(G)\geq t$. But, $\sigma(G)\leq 2$ and so $t\leq 2$. As $t\geq 2$, it follows $t=2$. Hence, $G=G_{1}\times G_{2}$ and so by Theorem~\ref{thm:product-atoms} $\Soc(G)=\Soc(G_{1})\times\Soc(G_{2})$. Now, suppose, $G_{1}$ is not monolithic. Then, there exist distinct $A,A'\in\MinNorm(G_{1})$. As $G_{2}\ne 1$, let $B\in\MinNorm(G_{2})$. Now, $\hat{A},\hat{A'},\hat{B}$ are three distinct normal subgroups of $G$. As $A\cap A'=1$ by Lemma~\ref{lem:distinct-atoms}, and $1\ne\hat{B}\leq G_{2}^{*}$ while $\hat{A}\hat{A'}\leq G_{1}^{*}$, the inclusions below are strict. So we have a chain,
    \[
       1<\hat{A}<\hat{A}\hat{A'}<\hat{A}\hat{A'}\hat{B}\leq \Soc(G).
    \]
    This means $\sigma(G)\geq 3$ --- which is a contradiction, as $\sigma(G)\leq 2$. Thus, $G_{1}$ must be monolithic. Similarly, $G_{2}$ is also monolithic.\par
    Conversely, suppose $G=G_{1}\times G_{2}$ and $G_{1},G_{2}$ are monolithic. So by Theorem~\ref{thm:product-atoms} we have $\Soc(G)=\Soc(G_{1})\times\Soc(G_{2}).$ Now we claim the following
    \[\text{Claim: If }\Soc(G)=\Soc(G_{1})\times \Soc(G_{2})\text{ then }\sigma(G)=\sigma(G_{1})+\sigma(G_{2}).
    \]
    Now, by Proposition~\ref{prop:socle-geometry} we have:
    \[
       \Soc(G_{1})=B_{1,1}\times\cdots\times B_{1,r_{1}}\times\prod_{p}S_{1,p} \quad\text{and}\quad \Soc(G_{2})=B_{2,1}\times\cdots\times B_{2,r_{2}}\times\prod_{p}S_{2,p},
    \]
    where, $B_{1,1},\cdots,B_{1,r_{1}}$ and $B_{2,1},\cdots,B_{2,r_{2}}$ are non-abelian minimal normal subgroups of $G_{1}$ and $G_{2}$ respectively, and $S_{i,p}$ are the subgroups generated by the abelian minimal normal $p$-subgroups. Then, we get:
    \[
       \Soc(G)=B_{1,1}\times\cdots\times B_{1,r_{1}}\times B_{2,1}\times\cdots\times B_{2,r_{2}}\times\prod_{p}(S_{1,p}\times S_{2,p}).
    \]
    Suppose, $S_{p}=S_{1,p}\times  S_{2,p}$. By Theorem~\ref{thm:product-atoms}, the non-abelian minimal normal subgroups of $G$ are the $B_{i,k}$, and the abelian minimal normal $p$-subgroups of $G$ lie in $S_{p}$ and generate it. So, $S_{p}$ is the subgroup of Proposition~\ref{prop:socle-geometry} for $G$. Then, as $\F_{p}G$-modules, $S_{p}=S_{1,p}\oplus S_{2,p}$. By Proposition~\ref{prop:socle-geometry}, each $S_{i,p}$ is a semi-simple $\F_{p}G_{i}$-module. So, we can consider an isotypic decomposition for both of them as:
    \[
       S_{1,p}\cong\bigoplus_{\alpha\in I_{1}}X_{\alpha}^{n_\alpha} \quad\text{and}\quad S_{2,p}\cong\bigoplus_{\beta\in I_{2}}Y_{\beta}^{m_{\beta}},
    \]
    here and below, these are isomorphisms of modules, and $I_{1},I_{2}$ (and the exponents $r,s$ introduced below) depend on $p$, where $X_{\alpha}$'s are pairwise non-isomorphic simple $\F_{p}G_{1}$-modules and $Y_{\beta}$'s are pairwise non-isomorphic simple $\F_{p}G_{2}$-modules. Now, when each $X_{\alpha}$ and each $Y_{\beta}$ is considered as an $\F_{p}G$-module, $G_{i}$ acts on them via the projections $\pi_{1}$ and $\pi_{2}$. So, if we consider $x\in X_{\alpha}$ and $(1,g_{2})\in G_{2}^{*}$, then $(1,g_{2})\cdot x=\pi_{1}(1,g_{2})\cdot x=x$. Similarly, for $y\in Y_{\beta}$ and $(g_{1},1)\in G_{1}^{*}$, $(g_{1},1)\cdot y=y$. So, $G_{2}$ acts trivially on each $X_{\alpha}$ and $G_{1}$ acts trivially on each $Y_{\beta}$. Now, suppose, $X_{\alpha}\cong Y_{\beta}$ as $\F_{p}G$-modules. Suppose, $f:X_{\alpha}\xrightarrow{\cong}Y_{\beta}$. Then, for any $g\in G_{1}$ and $x\in X_{\alpha}$, $f(g\cdot x)=g\cdot f(x)$. But $f(x)\in Y_{\beta}$ and $G_{1}$ acts trivially on $Y_{\beta}$, so 
    \[
       f(g\cdot x)=g\cdot f(x)=f(x).
    \]
    As $f$ is an isomorphism, and hence is injective, it follows that $g\cdot x=x$. Thus, $G_{1}$ also acts trivially on $X_{\alpha}$. So, $G$ acts trivially on $X_{\alpha}$ and hence $X_{\alpha}$ is the trivial $\F_{p}G$-module. So, $Y_{\beta}$ is also trivial $\F_{p}G$-module. So, $X_{\alpha}\cong Y_{\beta}$ if and only if they are the same trivial $\F_{p}G$-module. So, we get
    \[
       S_{1,p}=\left(\bigoplus_{\alpha\in I'_{1}}X_{\alpha}^{n_{\alpha}}\right)\oplus T^{r} \quad\text{and}\quad S_{2,p}=\left(\bigoplus_{\beta\in I'_{2}}Y_{\beta}^{m_{\beta}}\right)\oplus T^{s}
    \]
    as $\F_{p}G$-modules, where $I'_{1}$ and $I'_{2}$ are the index sets for which $X_{\alpha}$ and $Y_{\beta}$ are non-trivial respectively. So, we get
    \[
       S_{p}=S_{1,p}\oplus S_{2,p}=\left(\bigoplus_{\alpha\in I'_{1}}X_{\alpha}^{n_{\alpha}}\right)\oplus\left(\bigoplus_{\beta\in I'_{2}}Y_{\beta}^{m_{\beta}}\right)\oplus T^{r+s}.
    \]
    Now, by Proposition~\ref{prop:socle-geometry}, we have
    
    \[
       \sigma(G_{1})=r_{1}+\sum_{p}\Bigl(\sum_{\alpha\in I'_{1}}n_{\alpha}+r\Bigr)\quad\text{and}\quad\sigma(G_{2})=r_{2}+\sum_{p}\Bigl(\sum_{\beta\in I'_{2}}m_{\beta}+s\Bigr).
    \]
    But, $\displaystyle{\sum_{\alpha\in I'_{1}}n_{\alpha}+r=\sum_{\alpha\in I_{1}}n_{\alpha}}$ and $\displaystyle{\sum_{\beta\in I'_{2}}m_{\beta}+s=\sum_{\beta\in I_{2}}m_{\beta}}$. So,
    \[
       \sigma(G_{1})=r_{1}+\sum_{p}\sum_{\alpha\in I_{1}}n_{\alpha} \quad\text{and}\quad \sigma(G_{2})=r_{2}+\sum_{p}\sum_{\beta\in I_{2}}m_{\beta}.
    \]
    Finally then we get
    \begin{align*}
        \sigma(G_{1})+\sigma(G_{2})&=r_{1}+\sum_{p}\sum_{\alpha\in I_{1}}n_{\alpha}+r_{2}+\sum_{p}\sum_{\beta\in I_{2}}m_{\beta} \\
        &=r_{1}+r_{2}+\sum_{p}\left(\sum_{\alpha\in I_{1}}n_{\alpha}+\sum_{\beta\in I_{2}}m_{\beta}\right) \\
        &=\sigma(G).
    \end{align*}
    The last equality again holds by Proposition~\ref{prop:socle-geometry}.\par
    Now, as $G_{1}$ and $G_{2}$ are monolithic, it follows that $\sigma(G_{1})=1=\sigma(G_{2})$. Hence, by our claim, it follows that $\sigma(G)=2$. So, by Theorem~\ref{thm:Mfull-classification}, it follows then, $G=G_{1}\times G_{2}$ is $\M$-full.
\end{proof}

Corollary~\ref{cor:product-monolithic} allows exactly two monolithic factors. The next example shows that two cannot be replaced by three.
 
\begin{example}\label{ex:S33}
The group $S_3\times S_3$ is $\M$-full by Corollary~\ref{cor:product-monolithic}, since $S_3$
is monolithic with minimal normal subgroup $A_3$. But $S_3\times S_3\times S_3$ is not. Let
$A_3^{(i)}$ denote $A_3$ placed in the $i$-th coordinate of $S_3\times S_3\times S_3$. Since $A_3$ is
not central in $S_3$, Theorem~\ref{thm:product-atoms} shows that $A_3^{(1)},A_3^{(2)},A_3^{(3)}$ are exactly the minimal normal subgroups of $S_3\times S_3\times S_3$. The normal subgroup $A_3\times A_3\times1$ contains $A_3^{(1)}$ and $A_3^{(2)}$ but not $A_3^{(3)}$, so it is not an $\M$-subgroup by Corollary~\ref{cor:finite-classification}.
\end{example}

 The same phenomenon occurs for groups of continuous maps into a finite group, where the factors correspond to the pieces of a partition of the underlying space into open connected subsets.
\begin{example}\label{ex:CXG}
Let $X$ be a non-empty topological space and $G$ a non-trivial finite
group with the discrete topology. Write $C(X,G)$ for the group of
continuous maps $X\to G$, with pointwise multiplication. Then
$C(X,G)$ is $\M$-full if and only if either
\begin{enumerate}[label=\upshape(\roman*)]
\item $X$ is connected and $G$ is $\M$-full; or
\item $X$ is the disjoint union of two non-empty open connected subsets
and $G$ is monolithic.
\end{enumerate}
Indeed, in case (i) every continuous map is constant, so
$C(X,G)\cong G$. In case (ii), restriction to the two connected
components gives $C(X,G)\cong G\times G$, which is $\M$-full by
Corollary~\ref{cor:product-monolithic}.

Conversely, suppose $C(X,G)$ is $\M$-full and $X$ is disconnected.
If $X$ admitted a partition into three non-empty clopen subsets
$U_1,U_2,U_3$, let $K_i$ be the normal subgroup of continuous maps
that are trivial on $U_i$. Choose a nonidentity $g\in G$ and define
$e_i$ to be $g$ on $U_i$ and $1$ elsewhere. Then
\[
 e_3\in K_1\cap K_2\setminus\{1\},\qquad
 e_1\in K_3\cap K_2\setminus\{1\},\qquad
 e_2\in K_3\cap K_1\setminus\{1\},
\]
whereas $K_1\cap K_2\cap K_3=1$. Thus $K_3$ would fail the
$\M$-condition. Hence no such partition exists. Any decomposition
$X=U\sqcup V$ into two non-empty clopen subsets must therefore
have $U$ and $V$ connected, for otherwise one could split a
disconnected part and obtain a partition into three. Consequently
$C(X,G)\cong G\times G$, and Corollary~\ref{cor:product-monolithic}
forces $G$ to be monolithic. The connected case reduces to $G$
itself.
\end{example}

We end this section with the two extreme kinds of characteristically simple
groups, where the $\M$-lattice can be written down completely; they will
reappear as the crown-based powers with trivial top in
Section~\ref{sec:crowns}.

\begin{proposition}\label{prop:simple-powers}
Let $k\geq1$. Then the following statements hold. 
\begin{enumerate}
\item If $T$ is a non-abelian finite simple group, then every normal subgroup
      of $T^k$ is a product of coordinate factors.  Hence
      \[
      \N^{\M}(T^k)=
      \begin{cases}
      \{1,T\},&k=1,\\
      \{1,T_1,\ldots,T_k,T^k\},&k\geq2,
      \end{cases}
      \]
      where $T_i$ is the $i$-th coordinate factor.
\item For the elementary abelian group $C_p^k$, the $\M$-subgroups are
      $0$, all one-dimensional subspaces, and the whole space.  For $k\leq2$
      this is the entire subgroup lattice; for $k\geq3$ it is a proper
      meet-subsemilattice.
\end{enumerate}
\end{proposition}

\begin{proof}
    For $(1)$, let $T^{k}=T_{1}\times\cdots\times T_{k}$, where $T_{i}$ is the $i$-th coordinate subgroup. Then, $T_{i}\trnl T^{k}$. Now suppose $N\in\N(T^{k})$. Then, as $N\trnl T^{k}$ and $T_{i}\trnl T^{k}$ we have $N\cap T_{i}\trnl T^{k}$. In particular, $N\cap T_{i}\trnl T_{i}$. But, $T_{i}\cong T$. So
    \[
       N\cap T_{i}=1 \quad\text{or}\quad N\cap T_{i}=T_{i}.
    \]
    Define, $I=\{i : N\cap T_{i}=T_{i}\}$. Then, $\displaystyle{\prod_{i\in I}}T_{i}\leq N$. Let $n=(t_{1},\cdots,t_{k})\in N$ and suppose $t_{i}\ne1$. As $T$ is non-abelian and simple, we get $Z(T)=1$. Then, there exists $x_{i}\in T_{i}$ such that $[t_{i},x_{i}]\ne1$. As $N\trnl T^{k}$, it follows that $[n,x_{i}]\in N$. So, $[n,x_{i}]=(1,\cdots,[t_{i},x_{i}],\cdots,1)\in T_{i}$ and thus $1\ne[n,x_{i}]\in N\cap T_{i}$. Now, since $T_{i}\cong T$ is simple, so $N\cap T_{i}=T_{i}$ and so $i\in I$. Therefore,
    \[
       n\in\prod_{i\in I}T_{i} \quad\text{and}\quad N\leq\prod_{i\in I}T_{i}.
    \]
    So, $N=\displaystyle{\prod_{i\in I}T_{i}}$.\par
    Conversely, for any $I\subseteq\{1,\cdots,k\}$ we have $\displaystyle{\prod_{i\in I}T_{i}}\trnl T^{k}$. So,
    \[
       \N(T^{k})=\bigl\{\prod_{i\in I}T_{i}:I\subseteq\{1,\cdots,k\}\bigr\}.
    \]
    Now, as $T_{i}\cong T$ is non-abelian and simple, $\MinNorm(T^{k})=\{T_{1},\cdots,T_{k}\}$. So when $k=1$, $T^{1}=T$ and hence $\N(T^{1})=\N(T)=\{1,T\}$. For any $E\trnl T^{1}=T$, $E=1$ or $E=T$. But by Corollary~\ref{cor:finite-classification},
    \[
       E\in\N^\M(T^{1})\iff E=1\quad\text{or}\quad|\supp(E)|=1\quad\text{or}\quad\Soc(T^{1})\leq E
.    \]
    So, if $E=1$, then $E\in\N^\M(T^{1})$. If $E=T$, then as $\MinNorm(T^{1})=\MinNorm(T)=\{T\}$, $\Soc(T^{1})=T\leq E$. Thus, $E\in\N^\M(T^{1})$. Hence,
    \[
       \N^\M(T^{1})=\{1,T\}.
    \]
    For $k\geq2$, let $E\trnl T^{k}$. Then, as we showed either $E=1$ or $E=\displaystyle{\prod_{i\in I}T_{i}}$ for some $I\subseteq\{1,\cdots,k\}$. If $E=1$, then trivially $E\in\N^\M(T^{k})$. If $E=\displaystyle{\prod_{i\in I}T_{i}}$, then we have two cases:
    
    \text{Case }1: $|I|=1$. Then, $E=T_{i}$ and hence $\supp(E)=\{T_{i}\}$. So, $E\in\N^\M(T^{k})$, by Corollary~\ref{cor:finite-classification}.

       \text{Case }2: $|I|\geq2$. As $\Soc(T^{k})=\Soc(T_{1})\times\cdots\times\Soc(T_{k})$, and each $T_{i}\cong T$, we have,
        \[
           \Soc(T^{k})=T_{1}\times\cdots\times T_{k}=T^{k}.
        \]
        Since, $|I|\geq2$, so $E$ contains at least two distinct minimal normal subgroups, say $T_{i}$ and $T_{j}$. But, if $I\ne\{1,\cdots,k\}$, then
        \[
           E=\prod_{i\in I}T_{i}<T_{1}T_{2}\cdots T_{k}=\Soc(T^{k})
        \]
        and $|\supp(E)|\geq2$. Therefore, $E$ is not an $\M$-subgroup. But when $I=\{1,\cdots,k\}$, then $E=T_{1}\cdots T_{k}=T^{k}=\Soc(T^{k})$. Hence, $E\in\N^\M(T^{k})$.

        So for $k\geq2$,
        \[
           \N^\M(T^{k})=\{1,T_{1},\cdots,T_{k},T^{k}\}.
        \]

        For $(2)$, let $G=C_{p}^{k}$. Then, $G$ is an elementary abelian $p$-group. So, $G$ is an $\F_{p}$-vector space. Consequently, every subgroup of $G$ is normal and hence $E\trnl G$ implies that $E$ is an $\F_{p}$-subspace. By Corollary~\ref{cor:pgroup-Mfull}, $\Soc(G)=\Omega_{1}(Z(G))$. But $G$ is abelian and is a non-trivial finite $p$-group. So $Z(G)=G$ and hence, $\Soc(G)=\Omega_{1}(G)=G$. Now, if $E\trnl G$ and $E=0$, then clearly $E\in\N^\M(G)$. If $E\ne0$, then let $\dim_{\F_{p}}E=d$. By Corollary~\ref{cor:pgroup-Mfull}, the minimal normal subgroups of $G$ are all the one-dimensional $\F_{p}$-subspaces of $G$. So, by Corollary~\ref{cor:product-count}, the number of minimal normal subgroups contained in $E$ is
        \[
           \frac{p^{d}-1}{p-1}.
        \]
        As $E\ne0$, 
        \[
           E\in\N^\M(G)\iff|\supp(E)|=1\quad\text{or}\quad\Soc(G)\leq E.
        \]
        But, $\Soc(G)\leq E$ if and only if $E=G$ and $|\supp(E)|=1$ if and only if $\frac{p^{d}-1}{p-1}=1$ if and only if $d=1$ if and only if $\dim_{\F_{p}}E=1$. So,
        \[
           E\in\N^\M(G)\iff E=0\quad\text{or}\quad\dim_{\F_{p}}E=1\quad\text{or}\quad E=G.
        \]
        Now, as $G=C_{p}^{k}$ is abelian, $\N(G)=\Sub(G)$. If $k=1$, then $G=C_{p}$ and hence $\N(G)=\{0,C_{p}\}$. So clearly, $\N^\M(G)=\N(G)$. So, $G$ is $\M$-full and $\Sub(G)=\N^\M(G)$. If $k=2$, then $G=C_{p}^{2}$. But $C_{p}$ is monolithic and so by Corollary~\ref{cor:product-monolithic}, $\N^\M(G)=\N(G)=\Sub(G)$. So for $k\leq 2$, the $\M$-subgroup lattice is precisely the entire subgroup lattice of $G$.\par
        If $k\geq3$, then we already have 
        \[
           \N^\M(G)=\{0\}\sqcup\{L\leq C_{p}^{k}:\dim_{\F_{p}}L=1\}\sqcup\{C_{p}^{k}\}.
        \]
        As $k\geq3$, there exists a $2$-dimensional $\F_{p}$-subspace $E$ of $C_{p}^{k}$. So, $E\notin\N^\M(G)$ but $E\trnl G$ and hence $E\in\N(G)=\Sub(G)$. So, $\N^\M(G)\subsetneq\Sub(G)$. By Proposition~\ref{prop:meet}, for any two $E,F\in\N^\M(G)$
        \[
           E\wedge F=E\cap F\in\N^\M(G).
        \]
        So, $\N^\M(G)$ is a proper meet-subsemilattice of $\Sub(G)$.
\end{proof}

\section{Crown-based powers: normal lattices and quotients}
\label{sec:crowns}

Let $L$ be a finite monolithic primitive group, let
$A=\Soc(L)$, and put $H=L/A$.  For $k\geq1$, the \emph{crown-based power}
of $L$ of size $k$ is
\[
L_k=\{(x_1,\ldots,x_k)\in L^k:x_1A=\cdots=x_kA\}.
\]
Write $V=A^k$.  Then $V\trnl L_k$ and $L_k/V\cong H$.  Crown-based powers
were introduced in \cite{DallaVoltaLucchini} to describe the finite groups
that need more generators than any proper quotient, and through the
isomorphism $G/R_G(B)\cong(L_B)_{\delta_G(B)}$ of Detomi and
Lucchini~\cite{DetomiLucchiniCrowns} they describe the crowns of an arbitrary
finite group (Corollary~\ref{cor:arbitrary-crown}). This crown-quotient
construction is also used in the study of maximal-subgroup intersection
numbers~\cite[Section~5.1]{BurnessGaronziLucchini};
see~\cite{LucchiniCrowns} for a survey. It is known that
$\Soc(L_k)=A^k$ is a direct product of $k$ minimal normal
subgroups~\cite{DallaVoltaLucchini}. Furthermore,
Lucchini, Moscatiello and Spiga~\cite[Lemma~2.2]{LucchiniMoscatielloSpiga}
proved that every normal subgroup of $L_k$ is comparable with $V$.
We recall the proof of this fact and combine it with the standard
module-theoretic description of $V$ to obtain an explicit presentation
of the whole normal-subgroup lattice, which is the starting point for
our results on $\M$-lattices, rigidity and quotients.

\begin{theorem}\label{thm:crown-normal-lattice}
Let $L$, $A$, $H$, $L_k$, and $V=A^k$ be as above.
\begin{enumerate}
\item \textup{\cite[Lemma~2.2]{LucchiniMoscatielloSpiga}.}
      Every normal subgroup of $L_k$ is comparable with $V$:
      \[
      N\trnl L_k\quad\Longrightarrow\quad N\leq V\ \text{or}\ V\leq N.
      \]
      Consequently,
      \[
      \N(L_k)=[1,V]_{\N(L_k)}\cup[V,L_k]_{\N(L_k)},
      \qquad
      [1,V]_{\N(L_k)}\cap[V,L_k]_{\N(L_k)}=\{V\}.
      \]
\item The quotient map induces a lattice isomorphism
      \[
      [V,L_k]_{\N(L_k)}\longrightarrow \N(H),
      \qquad N\longmapsto N/V.
      \]
\item If $A$ is non-abelian, then
      \[
      I\longmapsto A_I:=\prod_{i\in I}A_i
      \]
      is a lattice isomorphism from the Boolean lattice of subsets of
      $\{1,\ldots,k\}$ onto $[1,V]_{\N(L_k)}$.  Here $A_i$ is the $i$-th
      coordinate copy of $A$ and $A_\varnothing=1$.
\item If $A$ is abelian, then $A$ is an irreducible faithful
      $\F_pH$-module for some prime $p$.  With
      \[
      D=\End_{\F_pH}(A),
      \]
      the map
      \[
      K\longmapsto\Hom_{\F_pH}(A,K)
      \]
      identifies $[1,V]_{\N(L_k)}$ with the lattice of $D$-subspaces of
      $D^k$.
\end{enumerate}
In particular, $\Soc(L_k)=V$.  Its minimal normal subgroups are the
coordinate factors $A_1,\ldots,A_k$ when $A$ is non-abelian, and are naturally
parametrized by $\mathds P^{k-1}(D)$ when $A$ is abelian.
\end{theorem}

\begin{proof}
(1) For completeness we give the argument behind the cited comparability
lemma. Let $C=C_L(A)=\{x\in L:xax^{-1}=a\text{ for every }a\in A\}$. Then $C\trnl L$. We first show that $C\leq A$. If $A$ is non-abelian, then $C\cap A=Z(A)=1$, as $Z(A)$ is characteristic in $A$, hence normal in $L$, and $Z(A)\ne A$. As $L$ is monolithic, every non-trivial normal subgroup of $L$ contains $A$, and so $C=1$. Now suppose $A$ is abelian, so that $A\leq C$. Let $M$ be a maximal subgroup of $L$ with $\operatorname{core}_{L}(M)=1$. If $A\leq M$, then $A\leq\operatorname{core}_{L}(M)=1$, which is a contradiction. So $A\nleq M$, and as $M$ is maximal, $L=AM$. Now $C\cap M$ is normalized by $M$ and centralized by $A$. So $C\cap M\trnl AM=L$, and hence $C\cap M\leq\operatorname{core}_{L}(M)=1$. If $x\in C$, then $x=am$ for some $a\in A$ and $m\in M$, and $m=a^{-1}x\in C\cap M=1$. Thus, $x=a\in A$. So in both cases $C\leq A$, that is,
        \[
           x\in L\setminus A\implies [x,a]\ne1\text{ for some }a\in A.
        \]
        Now, let $N\trnl L_{k}$ with $N\nleq V$, and let $n=(n_{1},\cdots,n_{k})\in N\setminus V$. As $n_{1}A=\cdots=n_{k}A$ and $n\notin A^{k}$, $n_{j}\notin A$ for every $j$. Fix $j$. Then $[n_{j},a]\ne1$ for some $a\in A$. Let $\hat{a}_{j}=(1,\cdots,a,\cdots,1)\in A_{j}$, with $a$ in the $j$-th place. As $N\trnl L_{k}$, so
        \[
           1\ne[n,\hat{a}_{j}]=(1,\cdots,[n_{j},a],\cdots,1)\in N\cap A_{j}.
        \]
        Let $K\leq A$ be the set of $j$-th coordinates of the elements of $N\cap A_{j}$. As $N\cap A_{j}\trnl L_{k}$ and $(l,\cdots,l)\in L_{k}$ for every $l\in L$, $K\trnl L$. As $K\ne1$ and $A$ is minimal normal in $L$, $K=A$, that is, $A_{j}\leq N$. As $j$ was arbitrary, $V\leq N$. The rest of (1) follows at once.

(2) The map $L_{k}\rightarrow H$, $(x_{1},\cdots,x_{k})\mapsto x_{1}A$ is a surjective homomorphism with kernel $V$. So, by the correspondence theorem, $N\mapsto N/V$ is a bijection from $[V,L_{k}]_{\N(L_{k})}$ onto $\N(L_{k}/V)\cong\N(H)$, which preserves inclusion in both directions. Hence, it is a lattice isomorphism.

(3) Clearly, each $A_{I}\trnl L_{k}$, and $A_{I}\leq A_{J}$ if and only if $I\subseteq J$. Now let $N\trnl L_{k}$ with $N\leq V$, and put $I=\{i:A_{i}\leq N\}$. Then $A_{I}\leq N$. Let $n=(a_{1},\cdots,a_{k})\in N$ and $i\notin I$. If $a_{i}\ne1$, then as $Z(A)=1$, $[a_{i},a]\ne1$ for some $a\in A$. Then, as in (1), $1\ne[n,\hat{a}_{i}]\in N\cap A_{i}$ and so $A_{i}\leq N$, which is a contradiction. So $a_{i}=1$ for all $i\notin I$, and hence $N\leq A_{I}$. Thus, $N=A_{I}$. So $I\mapsto A_{I}$ is a bijection which preserves inclusion in both directions, and hence it is a lattice isomorphism.

(4) As $A$ is abelian, the Sylow subgroups of $A$ and $\{a\in A:a^{p}=1\}$ are characteristic in $A$, and hence normal in $L$. So, by minimality of $A$, $A$ is an elementary abelian $p$-group for some prime $p$. Now $L$ acts on $A$ by conjugation, and $A$ acts trivially as $A$ is abelian. So $A$ is an $\F_{p}H$-module. Its submodules are exactly the normal subgroups of $L$ contained in $A$, so $A$ is irreducible. It is faithful, as $lA\in H$ acts trivially only if $l\in C\leq A$. Moreover, for $(l_{i})\in L_{k}$ and $(a_{i})\in V$ we have $l_{i}a_{i}l_{i}^{-1}=l_{1}a_{i}l_{1}^{-1}$, as $l_{i}\in l_{1}A$ and $A$ is abelian. So, $V\cong A^{k}$ with $H$ acting diagonally, and $[1,V]_{\N(L_{k})}$ is the lattice of $\F_{p}H$-submodules of $A^{k}$.\par
        For $x\in D^{k}$, let $\Phi_{x}:A\rightarrow A^{k}$ be as in the proof of Lemma~\ref{lem:projective-submodules}. Then $x\mapsto\Phi_{x}$ is a bijection from $D^{k}$ onto $\Hom_{\F_{p}H}(A,A^{k})$, $\Phi_{x+y}=\Phi_{x}+\Phi_{y}$ and $\Phi_{xd}=\Phi_{x}\circ d$ for $d\in D$. For a submodule $K\leq A^{k}$ and a subspace $W\leq D^{k}$ put
        \[
           W(K)=\{x\in D^{k}:\im\Phi_{x}\leq K\}\quad\text{and}\quad K(W)=\sum_{x\in W}\im\Phi_{x}.
        \]
        So $W(K)$ corresponds to $\Hom_{\F_{p}H}(A,K)$, and it is a subspace of $D^{k}$, as $\im\Phi_{x+y}\leq\im\Phi_{x}+\im\Phi_{y}$ and $\im\Phi_{xd}\leq\im\Phi_{x}$. Clearly, both maps preserve inclusion. As $A^{k}$ is semisimple, so is $K$, and hence $K$ is a sum of simple submodules. By the proof of Lemma~\ref{lem:projective-submodules}, each of them is $\im\Phi_{x}$ for some $x\in D^{k}$, and then $x\in W(K)$. Hence, $K(W(K))=K$.\par
        Conversely, clearly $W\subseteq W(K(W))$. Let $b_{1},\cdots,b_{r}$ be a basis of $W$ and define $F:A^{r}\rightarrow K(W)$ by $F(a_{1},\cdots,a_{r})=\sum_{i}\Phi_{b_{i}}(a_{i})$. Then $F$ is a surjective homomorphism, as $\Phi_{x}=\sum_{i}\Phi_{b_{i}}\circ d_{i}$ for $x=b_{1}d_{1}+\cdots+b_{r}d_{r}\in W$. As $A^{r}$ is semisimple, there exists a homomorphism $s:K(W)\rightarrow A^{r}$ with $F\circ s=\id$. Now let $x\in W(K(W))$. Then $s\circ\Phi_{x}\in\Hom_{\F_{p}H}(A,A^{r})$, so $s\circ\Phi_{x}=(c_{1},\cdots,c_{r})$ for some $c_{i}\in D$. Hence,
        \[
           \Phi_{x}=F\circ s\circ\Phi_{x}=\sum_{i}\Phi_{b_{i}}\circ c_{i}=\Phi_{b_{1}c_{1}+\cdots+b_{r}c_{r}}
        \]
        and so $x=b_{1}c_{1}+\cdots+b_{r}c_{r}\in W$. Thus, $W(K(W))=W$. Therefore, $K\mapsto W(K)$ and $W\mapsto K(W)$ are mutually inverse bijections which preserve inclusion, and so $K\mapsto\Hom_{\F_{p}H}(A,K)$ is a lattice isomorphism onto the lattice of $D$-subspaces of $D^{k}$.

    Finally, if $N\in\MinNorm(L_{k})$, then by (1), $N\leq V$ or $V\leq N$, and in the latter case $N=V$. So $\MinNorm(L_{k})$ is the set of atoms of $[1,V]_{\N(L_{k})}$. These are the $A_{i}$ when $A$ is non-abelian, by (3), and they correspond to the one-dimensional subspaces of $D^{k}$, that is, to $\mathds{P}^{k-1}(D)$, when $A$ is abelian, by (4). In both cases each $A_{i}$ is an atom (in the abelian case $A_{i}\cong A$ is a simple submodule), and $V=A_{1}\cdots A_{k}$. Hence, $\Soc(L_{k})=V$.
\end{proof}

Theorem~\ref{thm:crown-normal-lattice} turns the lower part of the normal
lattice into an explicit finite geometry.  Here $\qbinom{k}{r}_q$ denotes the number of $r$-dimensional subspaces
of $\F_q^k$ when $0\leq r\leq k$, and is zero otherwise.

\begin{corollary}\label{cor:crown-M-lattice}
Assume $k\geq2$, and, in the abelian case, put $q=|D|$.  Set
\[
a_k=
\begin{cases}
k,&A\text{ non-abelian},\\[2mm]
\qbinom{k}{1}_q=\dfrac{q^k-1}{q-1},&A\text{ abelian}.
\end{cases}
\]
Then
\[
\N^{\M}(L_k)
 =\{1\}\,\sqcup\,\MinNorm(L_k)\,\sqcup\,[V,L_k]_{\N(L_k)}.
\]
Equivalently, the $\M$-lattice is obtained by placing $a_k$ incomparable
atoms between $1$ and $V$ and then attaching a copy of $\N(H)$ above $V$.
In particular,
\[
|\N^{\M}(L_k)|=|\N(H)|+a_k+1.
\]
Moreover,
\[
|\N(L_k)|=
\begin{cases}
|\N(H)|+2^k-1,&A\text{ non-abelian},\\[1mm]
|\N(H)|+\displaystyle\sum_{r=0}^k\qbinom{k}{r}_q-1,
   &A\text{ abelian},
\end{cases}
\]
and the number of normal subgroups discarded by the $\M$-closure system is
\[
|\N(L_k)\setminus\N^{\M}(L_k)|=
\begin{cases}
2^k-k-2,&A\text{ non-abelian},\\[1mm]
\displaystyle\sum_{r=2}^{k-1}\qbinom{k}{r}_q,&A\text{ abelian}.
\end{cases}
\]
For every $K\trnl L_k$ with $K\leq V$, the closure is
\[
\clM(K)=
\begin{cases}
K,&\lambda_{L_k}(K)\in\{0,1,k\},\\
V,&2\leq\lambda_{L_k}(K)\leq k-1.
\end{cases}
\]
The displayed loss vanishes exactly when $k=2$.
\end{corollary}

\begin{proof}
    Let $E\trnl L_{k}$. By Theorem~\ref{thm:crown-normal-lattice}, either $E\leq V$ or $V\leq E$. If $V\leq E$, then as $V=\Soc(L_{k})$, $E\in\N^\M(L_{k})$ by Corollary~\ref{cor:finite-classification}. Now suppose $E\leq V$. By Theorem~\ref{thm:crown-normal-lattice}, $[1,V]_{\N(L_{k})}$ is isomorphic to the Boolean lattice of subsets of $\{1,\cdots,k\}$ if $A$ is non-abelian, and to the lattice of $D$-subspaces of $D^{k}$ if $A$ is abelian. So $\lambda_{L_{k}}(E)=|I|$ if $E=A_{I}$, and $\lambda_{L_{k}}(E)=r$ if $E$ corresponds to an $r$-dimensional subspace. Moreover, the minimal normal subgroups contained in $E$ are the atoms below $E$, so
    \[
       |\supp(E)|=\begin{cases}
           \lambda_{L_{k}}(E), & A\text{ non-abelian},\\
           \qbinom{r}{1}_{q}, & A\text{ abelian}.
       \end{cases}
    \]
    Hence, $|\supp(E)|=1$ if and only if $\lambda_{L_{k}}(E)=1$, and $V\leq E$ if and only if $\lambda_{L_{k}}(E)=k$. So, by Corollary~\ref{cor:finite-classification},
    \[
       E\in\N^\M(L_{k})\iff\lambda_{L_{k}}(E)\in\{0,1,k\}\iff E=1,\ E\in\MinNorm(L_{k})\ \text{or}\ E=V.
    \]
   As $k\geq2$, $V$ is not minimal normal, and so the union below is disjoint. Thus,
    \[
       \N^\M(L_{k})=\{1\}\sqcup\MinNorm(L_{k})\sqcup[V,L_{k}]_{\N(L_{k})}.
    \]
    Now, by Theorem~\ref{thm:crown-normal-lattice}, $|[V,L_{k}]_{\N(L_{k})}|=|\N(H)|$, and $|[1,V]_{\N(L_{k})}|$ is $2^{k}$ or $\sum_{r=0}^{k}\qbinom{k}{r}_{q}$. Also, $L_{k}$ has exactly $a_{k}$ minimal normal subgroups, namely the $k$ factors $A_{i}$, or the ones corresponding to the $\qbinom{k}{1}_{q}$ one-dimensional subspaces of $D^{k}$. As the two intervals meet only in $V$, so
    \[
       |\N(L_{k})|=|[1,V]_{\N(L_{k})}|+|\N(H)|-1\quad\text{and}\quad|\N^\M(L_{k})|=|\N(H)|+a_{k}+1,
    \]
    which are the formulas in the statement. Subtracting, and using $\qbinom{k}{0}_{q}=\qbinom{k}{k}_{q}=1$, we get
    \[
       |\N(L_{k})\setminus\N^\M(L_{k})|=\begin{cases}
           2^{k}-1-k-1=2^{k}-k-2, & A\text{ non-abelian}\\
           \displaystyle{\sum_{r=0}^{k}\qbinom{k}{r}_{q}-1-\qbinom{k}{1}_{q}-1=\sum_{r=2}^{k-1}\qbinom{k}{r}_{q}}, & A\text{ abelian}.
       \end{cases}
    \]
    Now let $K\trnl L_{k}$ with $K\leq V$. If $K$ is an $\M$-subgroup, then $\clM(K)=K$. Otherwise $K\ne1$ and $|\supp(K)|\ne1$, so $|\supp(K)|\geq2$, and by Proposition~\ref{prop:Mclosure}, $\clM(K)=K\Soc(L_{k})=KV=V$. As $K$ is an $\M$-subgroup if and only if $\lambda_{L_{k}}(K)\in\{0,1,k\}$, this gives the formula for $\clM(K)$.\par
    Finally, if $k=2$, then $2^{k}-k-2=0$ and $\sum_{r=2}^{1}\qbinom{2}{r}_{q}=0$. If $k\geq3$, then $2^{k}-k-2\geq3$ and $$\sum_{r=2}^{k-1}\qbinom{k}{r}_{q}\geq\qbinom{k}{2}_{q}\geq1.$$ So the loss vanishes exactly when $k=2$.
\end{proof}

In particular, we recover the crown-based part of Theorem~C.

\begin{corollary}\label{cor:crown-M-full}
For every $k\geq1$, the crown-based power $L_k$ is $\M$-full if and only if
$k\leq2$.
\end{corollary}

\begin{proof}
If $k=1$, then $L^{k}=L$ and thus $L_{1}=L$. So as $L$ is finite and monolithic, it follows by Theorem~\ref{thm:Mfull-classification}(4a) that $L_{1}$ is $\M$-full. For $k\geq2$, by Corollary~\ref{cor:crown-M-lattice} we have $\N(L_{k})\setminus\N^\M(L_{k})=\emptyset$ exactly when $k=2$. So, $L_{k}$ is $\M$-full when $k=2$ and since it is also $\M$-full when $k=1$, so we get, $L_{k}$ is $\M$-full when $k\leq2$. Also, for $k\geq3$ we have $\N(L_{k})\setminus\N^\M(L_{k})\ne\emptyset$, so $L_{k}$ is not $\M$-full. Hence, $L_{k}$ is $\M$-full if and only if $k\leq2$.
\end{proof}

The complete lattice description also determines all maximal $\M$-chains,
not only their maximum length.

\begin{proposition}\label{prop:crown-M-chains}
Assume $k\geq2$.  The lattice $\N^{\M}(L_k)$ is graded and modular.  Its rank
function is
\[
\operatorname{rk}_{\M}(E)=
\begin{cases}
0,&E=1,\\
1,&E\in\MinNorm(L_k),\\
2+\lambda_H(E/V),&V\leq E.
\end{cases}
\]
Every maximal $\M$-chain has the form
\[
1<W<V=E_0<E_1<\cdots<E_s=L_k,
\]
where $W\in\MinNorm(L_k)$ and
\[
1=E_0/V<E_1/V<\cdots<E_s/V=H
\]
is a chief series of $H$.  Hence every maximal $\M$-chain has length
$\elln(H)+2$.  If $c(H)$ denotes the number of chief series of $H$, with
$c(1)=1$, then $L_k$ has exactly $a_kc(H)$ maximal $\M$-chains.

Furthermore, $\N^{\M}(L_k)$ is distributive if and only if $A$ is
non-abelian, $k=2$, and $\N(H)$ is distributive.
\end{proposition}

\begin{proof}
    By Corollary~\ref{cor:crown-M-lattice} we know that
    \[
       \N^\M(L_{k})=\{1\}\sqcup\MinNorm(L_{k})\sqcup[V,L_{k}]_{\N(L_{k})}.
    \]
    So we define, $\N_{1}=\{1\}\sqcup\MinNorm(L_{k})\cup\{V\}$ and $\N_{2}=[V,L_{k}]_{\N(L_{k})}$. Then, $\N^\M(L_{k})=\N_{1}\cup\N_{2}$ and $\N_{1}\cap\N_{2}=\{V\}$. Moreover, for any $X\in\N_{1}\setminus\{V\}$ and $Y\in\N_{2}\setminus\{V\}$, we have
    \[
       X<V<Y.
    \]
    Thus, $\N^\M(L_{k})=\N_{1}\oplus\N_{2}$, the vertical sum of $\N_{1}$ and $\N_{2}$, that is, $\N_{1}$ is placed below $\N_{2}$ and the top $V$ of $\N_{1}$ is identified with the bottom $V$ of $\N_{2}$. Now for any $M\ne N\in\MinNorm(L_{k})$ we have $M\cap N=1$ and $M\vee_{\M}N=\clM(MN)=V$ by Corollary~\ref{cor:crown-M-lattice}, although the product $MN$ itself is smaller than $V$ when $k\geq3$. So, $\N_{1}$ is a lattice of height $2$ and hence $\N_{1}$ is modular. Further by Theorem~\ref{thm:crown-normal-lattice}, $[V,L_{k}]_{\N(L_{k})}\cong\N(H)$ and by Lemma~\ref{lem:modular}, $\N(H)$ is modular. Also, $\N_{2}$ is a sublattice of $\N^\M(L_{k})$, as $\clM(EF)=EF$ for $E,F\in\N_{2}$. As we are gluing the bottom of $\N_{2}$ with the top of $\N_{1}$ and both $\N_{1},\N_{2}$ are modular, it follows that $\N^\M(L_{k})$ is also modular.\par
    Now consider $V\leq E\trnl L_{k}$. As $[V,L_{k}]_{\N(L_{k})}\cong\N(H)$, if we restrict this isomorphism, we get $[V,E]_{\N(L_{k})}\cong[1,E/V]_{\N(H)}$. So a chain,
    \[
       V=E_{0}<E_{1}<\cdots<E_{r}=E.
    \]
    corresponds to 
    \[
       1=E_{0}/V<E_{1}/V<\cdots<E_{r}/V=E/V.
    \]
    Thus, $\lambda_{H}(E/V)$ is the length of a maximal chain from $V$ to $E$. Now, for any $E\in\N^\M(L_{k})$, if $E=1$ then $\operatorname{rk}_{\M}(E)=0$. If $E\in\MinNorm(L_{k})$, then a maximal chain looks like $1<E$. Hence, $\operatorname{rk}_{\M}(E)=1$. If $E=V$, then a maximal chain is $1<M<E=V$, where $M\in\MinNorm(L_{k})$ and so $\operatorname{rk}_{\M}(E)=2$. If $V<E$, then a maximal chain from $1$ to $E$ looks like
    \[
       1<M<V=E_{0}<E_{1}<E_{2}<\cdots<E_{r}=E, 
    \]
    and hence, $\operatorname{rk}_{\M}(E)=2+\lambda_{H}(E/V)$. Thus the rank function increases by one across every cover and hence $\N^\M(L_{k})$ is graded. Now, observe that a maximal $\M$-chain must first pass through the chain $1<W<V$, where $W\in\MinNorm(L_{k})$. After that, it follows a maximal chain in $[V,L_{k}]_{\N(L_{k})}$
    \[
       V=E_{0}<E_{1}<\cdots<E_{s}=L_{k},
    \]
    which as we showed corresponds to a chief series,
    \[
       1=E_{0}/V<E_{1}/V<\cdots<E_{s}/V=L_{k}/V\cong H
    \]
    in $\N(H)$, due to the isomorphism $[V,L_{k}]_{\N(L_{k})}\cong\N(H)$ in Theorem~\ref{thm:crown-normal-lattice}. Thus, if $\elln(H)$ denotes the length of a chief series in $H$, the length of a maximal $\M$-chain in $L_{k}$ is $2+\elln(H)$. Furthermore, if $a_{k}$ denotes the number of atoms of $\N^\M(L_{k})$ as in Corollary~\ref{cor:crown-M-lattice}, then due to our two independent choices: first a $W\in\MinNorm(L_{k})$ and then a chief series of $H$, the number of maximal $\M$-chains that $L_{k}$ has is $a_{k}c(H)$.\par
    Now suppose $\N^\M(L_{k})$ is distributive and $a_{k}=|\MinNorm(L_{k})|$. If $N,M,K$ are three atoms such that their pairwise join is $V$, then, $M\wedge N=1$ and $M\wedge K=1$ and thus,
    \[
       M=M\wedge V=M\wedge(N\vee K)=(M\wedge N)\vee(M\wedge K)=1,
    \]
    which is a contradiction as $M$ is an atom. So, this tells us that $a_{k}\leq2$. Now if $A$ is abelian, then by Corollary~\ref{cor:crown-M-lattice} we have
    \[
       a_{k}=\qbinom{k}{1}_{q}=\frac{q^{k}-1}{q-1},
    \]
    where $q=|D|$. As $k\geq2$, $q^{k}-1\geq q^{2}-1$ and hence $a_{k}\geq q+1\geq3$ as $q\geq2$. So, $A$ cannot be abelian. Then, $A$ must be non-abelian and $a_{k}=k$. But $k\geq2$ and $k=a_{k}\leq2$ gives us $k=2$. Also as every normal subgroup containing $V$ is an $\M$-subgroup so $$[V,L_{k}]_{\N^\M(L_{k})}=[V,L_{k}]_{\N(L_{k})}\cong\N(H).$$ As $\N^\M(L_{k})$ is distributive, $[V,L_{k}]_{\N^\M(L_{k})}$ is also distributive and hence $\N(H)$ is distributive.\par
    Conversely suppose $A$ is non-abelian, $k=2$ and $\N(H)$ is distributive. Then, both the interval below $V$ and the interval above $V$ are distributive. But joining the top of the lower interval and the bottom of the above interval preserves distributivity and gives us $\N^\M(L_{k})$. Hence, $\N^\M(L_{k})$ is distributive.
\end{proof}

The full normal lattice and the $\M$-lattice retain different data
about a crown-based power. The following isomorphism criteria make this
comparison precise; notably, the $\M$-lattice need not determine either
the crown multiplicity or the field appearing in the abelian socle.

\begin{theorem}\label{thm:crown-lattice-rigidity}
Let $L$ and $L'$ be finite monolithic primitive groups with socles $A$ and
$B$, respectively.  Put
\[
H=L/A,\qquad H'=L'/B,
\]
and let $L_k$ and $L'_m$ be crown-based powers with $k,m\geq2$.  Define
\[
a(L,k)=|\MinNorm(L_k)|=
\begin{cases}
k,&A\text{ non-abelian},\\[1mm]
\dfrac{q^k-1}{q-1},&A\text{ abelian},
\end{cases}
\]
where $q=|\End_{\F_pH}(A)|$ in the abelian case, and define $a(L',m)$
analogously.
\begin{enumerate}
\item The normal lattices $\N(L_k)$ and $\N(L'_m)$ are isomorphic if and only
      if all of the following hold:
      \begin{enumerate}[label=\upshape(\alph*)]
      \item $k=m$;
      \item $A$ and $B$ are either both abelian or both non-abelian;
      \item $\N(H)\cong\N(H')$;
      \item in the abelian case,
            \[
            |\End_{\F_pH}(A)|=|\End_{\F_{p'}H'}(B)|,
            \]
            where $B$ is an irreducible $\F_{p'}H'$-module.
      \end{enumerate}
\item The $\M$-lattices $\N^{\M}(L_k)$ and $\N^{\M}(L'_m)$ are isomorphic if
      and only if
      \[
      a(L,k)=a(L',m)
      \qquad\text{and}\qquad
      \N(H)\cong\N(H').
      \]
\end{enumerate}
Thus the full normal lattice determines the crown multiplicity, the abelian
or non-abelian nature of the socle, and the order of the endomorphism field in
the abelian case.  None of these data is determined in general by the
$\M$-lattice.  Indeed, for every non-abelian finite simple group $T$,
\[
\N^{\M}(T^3)\cong\N^{\M}(C_2^2),
\qquad
\N(T^3)\not\cong\N(C_2^2),
\]
and one also has
\[
\N^{\M}(C_2^5)\cong\N^{\M}(C_5^3),
\qquad
\N(C_2^5)\not\cong\N(C_5^3).
\]
\end{theorem}

\begin{proof}
    For $(1)$, first suppose that $\Phi:\N(L_{k})\to\N(L_{m}')$ is an isomorphism. Let $V=A^{k}$ and $V'=B^{m}$. Then, 
    \[
       V=\bigvee\MinNorm(L_{k}) \quad\text{and}\quad V'=\bigvee\MinNorm(L_{m}').
    \]
    So we must have $\Phi(V)=V'$. So by Theorem~\ref{thm:crown-normal-lattice}, we get $[1,V]_{\N(L_{k})}\cong[1,V']_{\N(L_{m}')}$ and $k=\lambda_{L_{k}}(V)=\lambda_{L_{m}'}(V')=m$. Now, if $A$ is non-abelian then $[1,V]_{\N(L_{k})}$ is distributive as it is isomorphic to the Boolean lattice $\P(\{1,\cdots,k\})$. If $A$ is abelian, then $[1,V]_{\N(L_{k})}$ is not distributive, as $\Sub_{D}(D^{k})$ contains three distinct lines $\ell_{1},\ell_{2},\ell_{3}$ lying in one plane $P$ (as $k\geq2$ and $|D|+1\geq3$), and then
    \[
       \ell_{1}\wedge(\ell_{2}\vee \ell_{3})=\ell_{1}\wedge P=\ell_{1}\ne0=(\ell_{1}\wedge \ell_{2})\vee(\ell_{1}\wedge \ell_{3}).
    \]
    As $[1,V]_{\N(L_{k})}\cong[1,V']_{\N(L_{m}')}$, either both $A$ and $B$ are abelian or both are non-abelian. Since we also have $[V,L_{k}]_{\N(L_{k})}\cong\N(H)$ and $[V',L_{m}']_{\N(L_{m}')}\cong\N(H')$, so by our assumption it follows that $\N(H)\cong\N(H')$. Now if both $A$ and $B$ are abelian, then $[1,V]_{\N(L_{k})}\cong[1,V']_{\N(L_{m}')}$ so they contain the same number of atoms and hence,
    \[
       1+|D|+\cdots+|D|^{k-1}=\qbinom{k}{1}_{|D|}=\qbinom{m}{1}_{|D'|}=1+|D'|+\cdots+|D'|^{m-1},
    \]
   where $D=\End_{\F_{p}H}(A)$ and $D'=\End_{\F_{p'}H'}(B)$. As $k=m\geq2$, the sum $1+q+\cdots+q^{k-1}$ is strictly increasing in $q$, and hence $|D|=|D'|$. So, $|\End_{\F_{p}H}(A)|=|\End_{\F_{p'}H'}(B)|$.\par
    Conversely, assume that conditions $(a)$ through $(d)$ hold. Then clearly $[V,L_{k}]_{\N(L_{k})}\cong[V',L_{m}']_{\N(L_{m}')}$. Now if both $A$ and $B$ are non-abelian then
    \[
       [1,V]_{\N(L_{k})}\cong\P(\{1,\cdots,k\}) \quad\text{and}\quad [1,V']_{\N(L_{m}')}\cong\P(\{1,\cdots,m\}).
    \]
    But $k=m$ and hence $[1,V]_{\N(L_{k})}\cong[1,V']_{\N(L_{m}')}$. Now if both $A$ and $B$ are abelian then,
    \[
       [1,V]_{\N(L_{k})}\cong\Sub_{D}(D^{k}) \quad\text{and}\quad [1,V']_{\N(L_{m}')}\cong\Sub_{D'}(D'^{k}).
    \]
    But $k=m$ and $|D|=|D'|=q$ (say). By Wedderburn's little theorem, the finite division rings $D$ and $D'$ are fields. As, up to isomorphism, there is a unique field with $q$ elements, it follows that 
    \[
       D\cong D'.
    \]
    Let $\varphi:D\to D'$ be a field isomorphism. As $k=m$, this field isomorphism induces a bijection $\tilde{\varphi}:D^{k}\to D'^{k}$ given by $(d_{1},\cdots,d_{k})\longmapsto(\varphi(d_{1}),\cdots,\varphi(d_{k}))$. Now let $U\leq D^{k}$ be a $D$-subspace and $u=(u_{1},\cdots,u_{k})\in U$. Then, if $\varphi(d)=d'$ we have,
    \[
       \tilde{\varphi}(ud)=(\varphi(u_{1})\varphi(d),\cdots,\varphi(u_{k})\varphi(d))=(\varphi(u_{1})d',\cdots,\varphi(u_{k})d')=\tilde{\varphi}(u)d'.
    \]
    As $ud\in U$ it follows that $\tilde{\varphi}(u)d'\in\tilde{\varphi}(U)$. As $\varphi$ preserves addition, $\tilde{\varphi}$ also preserves addition and hence $\tilde{\varphi}(U)$ is a $D'$-subspace of $D'^{k}$. So, we obtain a map $\psi\colon\Sub_{D}(D^{k})\to\Sub_{D'}(D'^{k})$ defined as $U\longmapsto\tilde{\varphi}(U)$. This is an isomorphism as sums and intersections are preserved, and its inverse is $\psi^{-1}\colon U'\longmapsto\tilde{\varphi}^{-1}(U')$. Thus, $$[1,V]_{\N(L_{k})}\cong[1,V']_{\N(L_{m}')}.$$
    Now, both the lower and upper intervals of $\N(L_{k})$ and $\N(L_{m}')$ are isomorphic. Both isomorphisms send $V$ to $V'$, and by Theorem~\ref{thm:crown-normal-lattice}(1) every normal subgroup of $L_{k}$ (resp.\ $L_{m}'$) is comparable with $V$ (resp.\ $V'$). So the two isomorphisms glue to a bijection $\N(L_{k})\to\N(L_{m}')$ which preserves inclusion in both directions, and hence it is a lattice isomorphism.\par
    For $(2)$, first suppose that $\N^\M(L_{k})\cong\N^\M(L_{m}')$. As $V$ is the join of the atoms of $\N^\M(L_{k})$, and similarly for $V'$, the isomorphism maps $V$ to $V'$. Then by Corollary~\ref{cor:crown-M-lattice} we have $$[V,L_{k}]_{\N^\M(L_{k})}\cong[V',L_{m}']_{\N^\M(L_{m}')}$$ and hence $\N(H)\cong\N(H')$. Also, the intervals below $V$ and $V'$ are isomorphic, implying they have the same number of atoms. So, $a(L,k)=a(L',m)$. Conversely if $\N(H)\cong\N(H')$ and $a(L,k)=a(L',m)$ then $[V,L_{k}]_{\N^\M(L_{k})}\cong[V',L_{m}']_{\N^\M(L_{m}')}$. Also, as the interval below $V$ and $V'$ is determined completely by the number of atoms, it follows that $[1,V]_{\N^\M(L_{k})}\cong[1,V']_{\N^\M(L_{m}')}$. So if we glue these two isomorphisms at the common socle, we get an extended isomorphism to $\N^\M(L_{k})\cong\N^\M(L_{m}')$.\par
    Now let $T$ be a non-abelian finite simple group. Then $T$ is monolithic with $\Soc(T)=T$ and $H=T/T=1$. So the crown-based power of $T$ of size $3$ is $T^{3}$. As $T$ is non-abelian, we get the normal subgroup interval below the socle is $[1,T^{3}]_{\N(T^{3})}\cong\P(\{1,2,3\})$. So $T^{3}$ has exactly three atoms. Now observe that $C_{2}^{2}\cong\F_{2}^{2}$. So its subgroups are precisely its subspaces, and hence the three one-dimensional subspaces are precisely the three atoms. So, by (2), as $a(T,3)=3=a(C_{2},2)$ and $H=1$ in both cases,
    \[
       \N^\M(T^{3})\cong\N^\M(C_{2}^{2}).
    \]
    However, the full normal lattice of $T^{3}$ contains intermediate normal subgroups: $T_{1}T_{2},T_{2}T_{3},T_{1}T_{3}$ that are not $\M$-subgroups. As $[1,V]_{\N(L_{k})}\cong\P(\{1,2,3\})$, $|\N(T^{3})|=8$. But in $C_{2}^{2}$, as the dimension is $2$, there is no subspace strictly between a line and the whole space. Hence, $|\N(C_{2}^{2})|=1+3+1=5$. So, 
    \[
       \N(T^{3})\ncong\N(C_{2}^{2}).
    \]
    Now consider the groups $C_{2}$ and $C_{5}$. Observe that $\Soc(C_{2})=C_{2}$ and $\Soc(C_{5})=C_{5}$. So, it follows that 
    \[
       H=C_{2}/\Soc(C_{2})=1 \quad\text{and}\quad H'=C_{5}/\Soc(C_{5})=1.
    \]
    Moreover, $(C_{2})_{5}=C_{2}^{5}$ and $(C_{5})_{3}=C_{5}^{3}$. So both the socles of $C_{2}^{5}$ and $C_{5}^{3}$ are abelian. Now, the $1$-dimensional subspaces of $\F_{2}^{5}$ are:
    \[
       \qbinom{5}{1}_{2}=\frac{2^{5}-1}{2-1}=31
    \]
    and the $1$-dimensional subspaces of $\F_{5}^{3}$ are 
    \[
       \qbinom{3}{1}_{5}=\frac{5^{3}-1}{5-1}=31.
    \]
    Thus the $\M$-lattice of both $C_{2}^{5}$ and $C_{5}^{3}$ has height $2$ with the same number of atoms. Hence, 
    \[
       \N^\M(C_{2}^{5})\cong\N^\M(C_{5}^{3}).
    \]
    However, the normal subgroup interval below the socle of $C_{2}^{5}$ contains proper non-zero subspaces of dimensions $1$, $2$, $3$ and $4$, whereas the normal subgroup interval below the socle of $C_{5}^{3}$ contains proper non-zero subspaces of dimensions $1$ and $2$ only. Thus, $\N(C_{2}^{5})$ has a height of $5$ and $\N(C_{5}^{3})$ has a height of $3$. Hence $\N(C_{2}^{5})\ncong\N(C_{5}^{3})$.
    \end{proof}

The quotient of a crown-based power by a coordinate minimal normal
subgroup is naturally the crown-based power with one fewer coordinate,
a standard feature of the construction~\cite{DallaVoltaLucchini}.
We extend this observation to \emph{every} normal subgroup contained
in the socle, and also describe the automorphism orbits of such kernels
and their enumerating polynomial. The answer retains the
Boolean--Grassmann dichotomy of the lower normal lattice.

\begin{theorem}\label{thm:crown-quotients}
Let $K\trnl L_k$ with $K\leq V$, put
\[
r=\lambda_{L_k}(K),\qquad m=k-r,
\]
and, in the abelian case, put $q=|D|$.
\begin{enumerate}
\item If $m=0$, then $L_k/K\cong H$.  If $m\geq1$, then
      \[
      L_k/K\cong L_m.
      \]
\item For each $r\in\{0,\ldots,k\}$, the normal subgroups $K\leq V$ of
      length $r$ form a single orbit under $\Aut(L_k)$.
\item The number of kernels $K\leq V$ for which
      $\lambda_{L_k/K}(V/K)=m$ is
      \[
      \begin{cases}
      \binom{k}{m},&A\text{ non-abelian},\\[2mm]
      \qbinom{k}{m}_q,&A\text{ abelian}.
      \end{cases}
      \]
      Equivalently, the socle-quotient enumerator
      \[
      \Psi_{L,k}(t)
      :=\sum_{\substack{K\trnl L_k\\K\leq V}}
      t^{\lambda_{L_k/K}(V/K)}
      \]
      is
      \[
      \Psi_{L,k}(t)=
      \begin{cases}
      (1+t)^k,&A\text{ non-abelian},\\[1mm]
      \displaystyle\sum_{m=0}^k\qbinom{k}{m}_q t^m,
        &A\text{ abelian}.
      \end{cases}
      \]
\item Every homomorphic image of $L_k$ is isomorphic either to $L_m$ for
      some $1\leq m\leq k$, or to a quotient of $H$.
\end{enumerate}
In particular, the normal subgroups contained in $V$ form exactly $k+1$
$\Aut(L_k)$-orbits, and their quotients have the $k+1$ distinct isomorphism
types
\[
H,L_1,\ldots,L_k.
\]
\end{theorem}

\begin{proof} 
(1) If $m=0$, then $K=V$, as $K\leq V$ and $\lambda_{L_{k}}(K)=k=\lambda_{L_{k}}(V)$. So $L_{k}/K=L_{k}/V\cong H$. Now let $m\geq1$. First assume $A$ is non-abelian. By Theorem~\ref{thm:crown-normal-lattice}, $K=A_{I}$ for a unique subset $I\subseteq\{1,\cdots,k\}$ and $r=|I|$. Let $J=\{1,\cdots,k\}\setminus I$, so $|J|=m$, and let $\pi_{J}:L_{k}\to L^{J}$ be the projection onto the $J$-coordinates. As $g_{1}A=\cdots=g_{k}A$ for every $(g_{1},\cdots,g_{k})\in L_{k}$, and as any $(g_{j})_{j\in J}$ with all $g_{j}A$ equal extends to an element of $L_{k}$ by putting $g_{i}=g_{j_{0}}$ for $i\in I$, for a fixed $j_{0}\in J$, relabelling the coordinates in $J$ as $1$ to $m$ we get
        \[
           \im\pi_{J}=\{(g_{j})_{j\in J}:g_{j}A\text{ is independent of }j\}\cong L_{m}.
        \]
        If $(g_{1},\cdots,g_{k})\in\ker\pi_{J}$, then $g_{j}=1$ for $j\in J$, and so $g_{i}A=g_{j}A=A$, i.e., $g_{i}\in A$ for all $i$. Thus $\ker\pi_{J}\leq A_{I}$, and clearly $A_{I}\leq\ker\pi_{J}$. So $\ker\pi_{J}=K$ and $L_{k}/K\cong L_{m}$.\par
        Now assume that $A$ is abelian. Recall that in the proof of Theorem~\ref{thm:crown-normal-lattice}(1), we had shown that $L=AM$ with $A\cap M=1$ (as $A\leq C$ and $C\cap M=1$ there), where $M$ is a core-free maximal subgroup of $L$. So $M\cong L/A=H$ and $L=A\rtimes H$. Identifying $H$ with its diagonal copy in $L^{k}$, we get $L_{k}=A^{k}\rtimes H$, with $H$ acting diagonally on $A^{k}$. As $K\trnl L_{k}$, $K$ is $H$-invariant and $$L_{k}/K\cong(A^{k}/K)\rtimes H.$$ As $A^{k}$ is semisimple, $A^{k}=K\oplus W$ for some submodule $W\cong A^{k}/K$, and $K$, $W$ are semisimple. As all composition factors of $A^{k}$ are isomorphic to $A$, and $K$, $W$ have composition lengths $r$ and $k-r$, $K\cong A^{r}$ and $A^{k}/K\cong A^{k-r}=A^{m}$ as $\F_{p}H$-modules. Now if $T\colon A^{k}/K\to A^{m}$ is an $\F_{p}H$-isomorphism, then $\tilde{T}\colon (x,h)\longmapsto(T(x),h)$ is a bijection from $(A^{k}/K)\rtimes H$ onto $A^{m}\rtimes H$, and as $T(h\cdot x)=h\cdot T(x)$,
        \[
           \tilde{T}((x_{1},h_{1})(x_{2},h_{2}))=(T(x_{1})+h_{1}\cdot T(x_{2}),h_{1}h_{2})=\tilde{T}((x_{1},h_{1}))\tilde{T}((x_{2},h_{2})).
        \]
        Thus we get
        \[
           L_{k}/K\cong(A^{k}/K)\rtimes H\cong A^{m}\rtimes H\cong L_{m}.
        \]

 (2) First assume $A$ is non-abelian. Every $\sigma\in S_{k}$ induces an automorphism $$\varphi_{\sigma}\colon (g_{1},\cdots,g_{k})\longmapsto(g_{\sigma^{-1}(1)},\cdots,g_{\sigma^{-1}(k)})$$ of $L_{k}$, and $\varphi_{\sigma}(A_{I})=A_{\sigma(I)}$. So if $K_{1}=A_{I}$ and $K_{2}=A_{J}$ with $|I|=|J|=r$, then choosing $\sigma\in S_{k}$ with $\sigma(I)=J$, we get $\varphi_{\sigma}(K_{1})=K_{2}$.\par
        Now suppose $A$ is abelian. Any $H$-endomorphism of $A^{k}$ can be represented by a $k\times k$-matrix over $D$, where the $(i,j)$-entry is the $H$-map from the $j$th copy of $A$ to the $i$th copy of $A$. So $$\Aut_{H}(A^{k})\cong GL_{k}(D).$$ As $\Hom_{\F_{p}H}(A,T(K))=T\circ\Hom_{\F_{p}H}(A,K)$, under the identification of Theorem~\ref{thm:crown-normal-lattice}$(4)$, each $T\in\Aut_{H}(A^{k})$ acts on the $D$-subspaces of $D^{k}$ through its matrix. Now let $K,K'\leq V$ be normal subgroups of length $r$, and let $U,U'$ be the corresponding $r$-dimensional $D$-subspaces of $D^{k}$. Extending bases $\{u_{1},\cdots,u_{r}\}$ of $U$ and $\{u_{1}',\cdots,u_{r}'\}$ of $U'$ to bases of $D^{k}$, we get a $D$-linear automorphism of $D^{k}$ with $u_{i}\mapsto u_{i}'$ for $1\leq i\leq r$, which maps $U$ onto $U'$. So there exists $T\in\Aut_{H}(A^{k})$ with $T(K)=K'$. By the same computation as in $(1)$, $\tilde{T}\colon (l,h)\longmapsto(T(l),h)$ is an automorphism of $L_{k}=A^{k}\rtimes H$, and $\tilde{T}(K)=K'$. So in both cases, all normal subgroups $K\leq V$ of length $r$ form a single orbit under $\Aut(L_{k})$.

 (3) By the correspondence theorem, $[1,V/K]_{\N(L_{k}/K)}\cong[K,V]_{\N(L_{k})}$, and as $\N(L_{k})$ is modular, so
        \[
           \lambda_{L_{k}/K}(V/K)=\lambda_{L_{k}}(V)-\lambda_{L_{k}}(K)=k-r.
        \]
        So the kernels $K\leq V$ with $\lambda_{L_{k}/K}(V/K)=m$ are precisely the normal subgroups $K\leq V$ with $\lambda_{L_{k}}(K)=k-m$. If $A$ is non-abelian, then by Theorem~\ref{thm:crown-normal-lattice} these are the $A_{I}$ with $|I|=k-m$, and so their number is $\binom{k}{k-m}=\binom{k}{m}$. If $A$ is abelian, then they correspond to the $(k-m)$-dimensional $D$-subspaces of $D^{k}$, and so their number is $\qbinom{k}{k-m}_{q}=\qbinom{k}{m}_{q}$, by the symmetry of the Gaussian binomial coefficients. Hence,
        \begin{align*}
            \Psi_{L,k}(t)&=\sum_{m=0}^{k}\#\{K\trnl L_{k}:K\leq V,\, \lambda_{L_{k}/K}(V/K)=m\}\,t^{m} \\
            &=\begin{cases}
                \displaystyle{\sum_{m=0}^{k}\binom{k}{m}t^{m}}=(1+t)^{k}, & A\text{ is non-abelian}, \\
                \displaystyle{\sum_{m=0}^{k}\qbinom{k}{m}_{q}t^{m}}, & A\text{ is abelian}.
            \end{cases}
        \end{align*}

  (4) Let $\pi:L_{k}\twoheadrightarrow Q$ be a surjective homomorphism and $N=\ker\pi$, so $Q\cong L_{k}/N$. By Theorem~\ref{thm:crown-normal-lattice}, either $N\leq V$ or $V\leq N$. If $N\leq V$, then by $(1)$, $Q\cong H$ or $Q\cong L_{m}$ with $m=k-\lambda_{L_{k}}(N)\geq1$. If $V\leq N$, then by the correspondence theorem $N/V\trnl L_{k}/V\cong H$, and hence
        \[
           Q\cong L_{k}/N\cong(L_{k}/V)/(N/V)
        \]
        is a quotient of $H$.

    Finally, every automorphism of $L_{k}$ induces an automorphism of $\N(L_{k})$, and so it preserves $\lambda_{L_{k}}$. Also, $A_{1}\cdots A_{r}\leq V$ is a normal subgroup of length $r$ for each $r\in\{0,\cdots,k\}$. So by $(2)$, the normal subgroups contained in $V$ form exactly $k+1$ $\Aut(L_{k})$-orbits. As $A\ne1$, for $m\ne n$ we have $$|L_{m}|=|A|^{m}|H|\ne|A|^{n}|H|=|L_{n}|,$$ and $|L_{m}|>|H|$ for $m\geq1$. So by $(1)$, the quotients coming from the $k+1$ orbits are pairwise non-isomorphic, and they are exactly $H,L_{1},\cdots,L_{k}$.
\end{proof}

For a finite group $X$ and $K\trnl X$, we call the quotient map
$\pi_K:X\to X/K$ \emph{$\M$-closure compatible} if
\[
\pi_K\bigl(\clM^{\,X}(N)\bigr)
 =\clM^{\,X/K}(N/K)
\qquad(K\leq N\trnl X).
\]

\begin{proposition}\label{prop:crown-closure-quotients}
Let $K\trnl L_k$.  The quotient map $\pi_K$ is $\M$-closure compatible if
and only if one of the following holds:
\begin{enumerate}
\item $K=1$;
\item $k=2$ and $K\in\MinNorm(L_2)$;
\item $V\leq K$ and
      \[
      H/(K/V)\cong L_k/K
      \]
      is $\M$-full.
\end{enumerate}
Thus, for $k\geq3$, no quotient by a non-trivial proper subgroup of the socle
is $\M$-closure compatible.
\end{proposition}

\begin{proof}
    First assume $\pi_{K}$ is $\M$-closure compatible. Taking $N=K$, we get $$\pi_{K}(\clM^{L_{k}}(K))=\clM^{L_{k}/K}(1)=1.$$ As $\ker\pi_{K}=K$, $\clM^{L_{k}}(K)\leq K$. But by definition $K\leq\clM^{L_{k}}(K)$. Thus, $K=\clM^{L_{k}}(K)$ and hence $K$ is an $\M$-subgroup. By Theorem~\ref{thm:crown-normal-lattice}, either $K<V$ or $V\leq K$.\par
    Now let $K<V$. If $k=1$, then $V=A$ is the only minimal normal subgroup of $L_{1}=L$, and so $K=1$. Suppose $k\geq2$. Then by Corollary~\ref{cor:crown-M-lattice}, $K=1$ or $K\in\MinNorm(L_{k})$. So let $K\in\MinNorm(L_{k})$ and suppose $k\geq3$. Choose $M\in\MinNorm(L_{k})$ with $M\ne K$ and put $N=MK$. If $A$ is non-abelian, then $K=A_{i}$ and $M=A_{j}$ for some $i\ne j$, and so $N=A_{\{i,j\}}$. If $A$ is abelian, then $K$ and $M$ correspond to distinct $1$-dimensional $D$-subspaces of $D^{k}$, and so $N$ corresponds to a $2$-dimensional subspace. So in both cases $\lambda_{L_{k}}(N)=2\leq k-1$. Hence $N<V$, and by Corollary~\ref{cor:crown-M-lattice}, $\clM^{L_{k}}(N)=V$. Also, as $\lambda_{L_{k}}(K)=1$ and $\N(L_{k})$ is modular by Lemma~\ref{lem:modular}, there is no normal subgroup of $L_{k}$ strictly between $K$ and $N$. Thus, $N/K\in\MinNorm(L_{k}/K)$ and hence $\clM^{L_{k}/K}(N/K)=N/K$. So,
    \[
       V/K=\pi_{K}(\clM^{L_{k}}(N))=\clM^{L_{k}/K}(N/K)=N/K.
    \]
    By the correspondence theorem we get $V=N$, which is a contradiction as $N<V$. Thus, $k=2$ and $K\in\MinNorm(L_{2})$.\par
    Now let $V\leq K$, and let $N/K\trnl L_{k}/K$ with $K\leq N\trnl L_{k}$. Then $V\leq N$, so $N$ contains the socle and hence is an $\M$-subgroup. Therefore,
    \[
       \clM^{L_{k}/K}(N/K)=\pi_{K}(\clM^{L_{k}}(N))=\pi_{K}(N)=N/K.
    \]
    Thus $L_{k}/K$ is $\M$-full, and $$L_{k}/K\cong(L_{k}/V)/(K/V)\cong H/(K/V).$$

    Conversely, if $K=1$, then $\pi_{1}$ is just the identification of $L_{k}$ with $L_{k}/1$, and so it is $\M$-closure compatible. Now assume $(2)$ or $(3)$ holds, and let $K\leq N\trnl L_{k}$. In case $(2)$, $L_{2}$ is $\M$-full by Corollary~\ref{cor:crown-M-full}, and $L_{2}/K\cong L_{1}$ by Theorem~\ref{thm:crown-quotients}, which is again $\M$-full by Corollary~\ref{cor:crown-M-full}. In case $(3)$, $V\leq N$, so $N$ is an $\M$-subgroup, and $L_{k}/K$ is $\M$-full by our assumption. So in both cases $\clM^{L_{k}}(N)=N$ and $\clM^{L_{k}/K}(N/K)=N/K$. Hence,
    \[
       \pi_{K}(\clM^{L_{k}}(N))=\pi_{K}(N)=N/K=\clM^{L_{k}/K}(N/K)
    \]
    and the compatibility holds. Finally, if $k\geq3$ and $1<K<V$, then none of $(1)$--$(3)$ holds, and so $\pi_{K}$ is not $\M$-closure compatible.
\end{proof}

\begin{definition}\label{def:quotient-M-full}
A finite group $X$ is said to be \emph{quotient-$\M$-full} if every homomorphic image of
$X$ is $\M$-full.
\end{definition}

The quotient classification makes this property completely recursive for
crown-based powers.

\begin{corollary}\label{cor:crown-quotient-M-full}
In the abelian case, put $q=|D|$.  Let $K\trnl L_k$ with $K<V$, and put
$m=\lambda_{L_k/K}(V/K)$.  Then
\[
L_k/K\text{ is $\M$-full}\quad\Longleftrightarrow\quad m\leq2.
\]
The number of such kernels $K<V$ is
\[
\begin{cases}
k+\binom{k}{2},&A\text{ non-abelian},\\[2mm]
\qbinom{k}{1}_q+\qbinom{k}{2}_q,&A\text{ abelian}.
\end{cases}
\]
The kernel $K=V$ gives one additional $\M$-full quotient exactly when $H$ is
$\M$-full.  Finally,
\[
L_k\text{ is quotient-$\M$-full}
\quad\Longleftrightarrow\quad
k\leq2\ \text{and}\ H\text{ is quotient-$\M$-full}.
\]
\end{corollary}

\begin{proof}
    As $K<V$, $m=\lambda_{L_{k}/K}(V/K)\geq1$. Also, by the proof of Theorem~\ref{thm:crown-quotients}$(3)$, $m=k-\lambda_{L_{k}}(K)$. Now,
    \[
       L_{k}/K\text{ is }\M\text{-full} \iff L_{m}\text{ is }\M\text{-full} \iff m\leq2,
    \]
    where the first equivalence follows from Theorem~\ref{thm:crown-quotients}$(1)$, as $L_{k}/K\cong L_{m}$, and the second follows from Corollary~\ref{cor:crown-M-full}. So the kernels $K<V$ with $L_{k}/K$ $\M$-full are precisely those with $m\in\{1,2\}$, and by Theorem~\ref{thm:crown-quotients}$(3)$ their number is
    \[
       \begin{cases}
            \binom{k}{1}+\binom{k}{2}=k+\binom{k}{2}, & A\text{ is non-abelian}\\
            \qbinom{k}{1}_{q}+\qbinom{k}{2}_{q}, & A\text{ is abelian}.
        \end{cases}
    \]
    When $K=V$, $L_{k}/K=L_{k}/V\cong H$, and so $L_{k}/V$ is $\M$-full if and only if $H$ is $\M$-full. This gives exactly one additional $\M$-full quotient in our count.\par
    Now suppose $L_{k}$ is quotient-$\M$-full. Then in particular $L_{k}$ is $\M$-full, and so $k\leq2$ by Corollary~\ref{cor:crown-M-full}. As $H$ is a homomorphic image of $L_{k}$, every quotient of $H$ is also a quotient of $L_{k}$, and hence $H$ is also quotient-$\M$-full.\par
    Conversely, assume $k\leq2$ and $H$ is quotient-$\M$-full. By Theorem~\ref{thm:crown-quotients}$(4)$, every homomorphic image of $L_{k}$ is isomorphic either to $L_{m}$ for some $1\leq m\leq k\leq2$, which is $\M$-full by Corollary~\ref{cor:crown-M-full}, or to a quotient of $H$, which is $\M$-full as $H$ is quotient-$\M$-full. Thus, $L_{k}$ is quotient-$\M$-full.
\end{proof}

Everything proved in this section applies to the crowns of an arbitrary
finite group. We use the notation of Detomi and
Lucchini~\cite{DetomiLucchiniCrowns}.

\begin{corollary}\label{cor:arbitrary-crown}
Let $G$ be a finite group, let $B$ be a non-Frattini chief factor of $G$, and let $L_B$ be the associated monolithic primitive group.  Write
\[
A_B=\Soc(L_B),\qquad H_B=L_B/A_B,
\]
and let $\delta_G(B)$ and $R_G(B)$ be the crown multiplicity and crown kernel,
respectively (see \cite[Definitions~6 and~8]{DetomiLucchiniCrowns}).  Define $I_G(B)$ by
\[
I_G(B)/R_G(B)=\Soc(G/R_G(B)).
\]
Then $\delta_G(B)\geq1$ and $I_G(B)/R_G(B)$ is the $B$-crown by \cite[Proposition~7]{DetomiLucchiniCrowns}, and, by \cite[Proposition~9]{DetomiLucchiniCrowns} (see also \cite[Section~5]{LucchiniCrowns}),
\[
G/R_G(B)\cong(L_B)_{\delta_G(B)}.
\]
Consequently:
\begin{enumerate}
\item the complete normal-subgroup lattice of $G/R_G(B)$ is the vertical
      gluing of $\N(H_B)$ with a Boolean lattice of rank $\delta_G(B)$ when
      $B$ is non-abelian, and with the subspace lattice of
      $D_B^{\,\delta_G(B)}$ when $B$ is abelian, with $\N(H_B)$ on top and the two glued at $I_G(B)/R_G(B)$, where $A_B$ is then an
      $\F_pH_B$-module for a prime $p$ and
      \[
      D_B=\End_{\F_pH_B}(A_B);
      \]
\item
      \[
      G/R_G(B)\text{ is $\M$-full}
      \quad\Longleftrightarrow\quad
      \delta_G(B)\leq2;
      \]
\item
      \[
      G/R_G(B)\text{ is quotient-$\M$-full}
      \quad\Longleftrightarrow\quad
      \delta_G(B)\leq2
      \ \text{and}\
      H_B\text{ is quotient-$\M$-full};
      \]
\item if $R_G(B)\leq K\trnl G$, then the quotient map
      \[
      G/R_G(B)\longrightarrow G/K,
      \]
      identifying $(G/R_G(B))/(K/R_G(B))$ with $G/K$, is $\M$-closure compatible if and only if one of the following holds:
      \begin{enumerate}[label=\upshape(\alph*)]
      \item $K=R_G(B)$;
      \item $\delta_G(B)=2$ and $K/R_G(B)$ is minimal normal in
            $G/R_G(B)$;
      \item $I_G(B)\leq K$ and $G/K$ is $\M$-full;
      \end{enumerate}
\item the number of normal subgroups $K$ with
      \[
      R_G(B)\leq K<I_G(B)
      \]
      for which $G/K$ is $\M$-full equals
      \[
      \begin{cases}
      \delta_G(B)+\binom{\delta_G(B)}{2},&B\text{ non-abelian},\\[2mm]
      \qbinom{\delta_G(B)}{1}_{q_B}
       +\qbinom{\delta_G(B)}{2}_{q_B},&B\text{ abelian},
      \end{cases}
      \]
      where $q_B=|D_B|$ in the abelian case.
\item if $\delta_G(B)\geq2$, then the abstract lattice
      $\N(G/R_G(B))$ determines $\delta_G(B)$, detects whether $B$ is
      abelian, and, in the abelian case, determines $q_B$.
\end{enumerate}
\end{corollary}

\begin{proof}
    For $(1)$ through $(6)$, apply Theorem~\ref{thm:crown-normal-lattice}, Corollaries~\ref{cor:crown-M-lattice}, \ref{cor:crown-M-full}, and \ref{cor:crown-quotient-M-full}, Proposition~\ref{prop:crown-closure-quotients}, and Theorem~\ref{thm:crown-lattice-rigidity}, with $L=L_B$ and $k=\delta_G(B)$.
\end{proof}

\section{Finite abelian groups}
\label{sec:abelian}

For finite abelian groups Theorem~A takes a particularly transparent form:
apart from $0$ and the subgroups above the socle, the $\M$-subgroups are the
non-trivial cyclic $p$-subgroups. This allows us to count them. For cyclic
groups $C_n$ the classification may be compared with the description of ring
$\M$-ideals of $\Z_n$ in~\cite{BlecherGoswami}.

\begin{theorem}\label{thm:abelian-classification}
Let $G$ be a finite abelian group and let $H\leq G$.  Then $H$ is an
$\M$-subgroup if and only if one of the following holds:
\begin{enumerate}
\item $H=0$;
\item $H$ is a non-trivial cyclic $p$-group for some prime $p$;
\item $\Soc(G)\leq H$.
\end{enumerate}
\end{theorem}

\begin{proof}
    Let $N\trnl G$ be a minimal normal subgroup and $0\ne x\in N$. As $G$ is abelian, $\langle x\rangle\trnl G$, and by minimality of $N$, $N=\langle x\rangle$ is cyclic. If $|N|=n$ is not a prime, then there exists $d\in\mathds{N}$ with $d\,|\,n$ and $1<d<n$, and the subgroup of $N$ of order $d$ is normal in $G$, which contradicts the minimality of $N$. Hence $|N|=p$ for some prime $p$. Conversely, a subgroup of $G$ of prime order is normal and has no non-trivial proper subgroups, and so it is minimal normal. Thus the minimal normal subgroups of $G$ are precisely its subgroups of prime order.\par
    Now let $H\ne0$ with $|\supp(H)|=1$. If $x,y\in H$ with $\ord(x)=p\ne q=\ord(y)$ for distinct primes $p,q$, such $x,y$ exist by Cauchy's theorem if $pq\,|\,|H|$, then $H$ would contain two minimal normal subgroups, contradicting $|\supp(H)|=1$. So $H$ is a $p$-group. If $H$ is not cyclic, then as $H$ is a finite abelian $p$-group we have
    \[
       H\cong C_{p^{\alpha_{1}}}\oplus\cdots\oplus C_{p^{\alpha_{n}}}
    \]
    with $n\geq2$. Let $x_{i}$ be a generator of the $i$th summand, and put
    \[
       x=(p^{\alpha_{1}-1}x_{1},0,\cdots,0) \quad\text{and}\quad y=(0,p^{\alpha_{2}-1}x_{2},0,\cdots,0).
    \]
    Then $\ord(x)=p=\ord(y)$ and $\langle x\rangle\cap\langle y\rangle=0$, so $\langle x\rangle$ and $\langle y\rangle$ are distinct minimal normal subgroups of $G$ contained in $H$. This again contradicts $|\supp(H)|=1$. Thus, $H$ is cyclic. Conversely, if $H$ is a non-trivial cyclic $p$-group, then $H$ has a unique subgroup of order $p$, and hence $|\supp(H)|=1$. Thus, for $H\ne0$,
    \[
       |\supp(H)|=1 \iff H \text{ is a cyclic }p\text{-group}.
    \]
    Now by Corollary~\ref{cor:finite-classification}, $H\in\N^\M(G)$ if and only if $H=0$, or $|\supp(H)|=1$, or $\Soc(G)\leq H$. So $H$ is an $\M$-subgroup of $G$ if and only if either $H=0$ or $H$ is a non-trivial cyclic $p$-group or $\Soc(G)\leq H$.
\end{proof}

Consequently $|\N^{\M}(G)|$ is determined by the number of cyclic
$p$-subgroups, which can be read off from the invariants of $G$.

\begin{proposition}\label{prop:abelian-count}
Let $G$ be a non-trivial finite abelian group, with notation as in Definition~\ref{def:abelian-notation}.  For every prime $p$ dividing $|G|$,
\[
c_p(G)=\sum_{k=1}^{\lambda_{p,1}}
\frac{p^{u_p(k)}-p^{u_p(k-1)}}{p^{k-1}(p-1)}.
\]
If $G\cong C_{p^a}$, then $|\N^{\M}(G)|=a+1$.  For every other non-trivial
finite abelian group,
\[
|\N^{\M}(G)|
=1+|\Sub(G/\Soc(G))|+\sum_{p\mid |G|}c_p(G).
\]
\end{proposition}
\begin{proof}
    Following the notation in Definition~\ref{def:abelian-notation}, let $G_{p}\cong C_{p^{\lambda_{p,1}}}\oplus\cdots\oplus C_{p^{\lambda_{p,r_{p}}}}$ and put $G_{p}[p^{k}]=\{x\in G_{p}:p^{k}x=0\}$. The condition $p^{k}x=0$ is checked coordinate-wise, and the elements $x_{i}\in C_{p^{\lambda_{p,i}}}$ with $p^{k}x_{i}=0$ form the whole of $C_{p^{\lambda_{p,i}}}$ if $k\geq\lambda_{p,i}$, and the unique subgroup of order $p^{k}$ if $k<\lambda_{p,i}$. Thus,
    \[
       |G_{p}[p^{k}]|=\prod_{i=1}^{r_{p}}p^{\min\{k,\lambda_{p,i}\}}=p^{\sum_{i=1}^{r_{p}}\min\{k,\lambda_{p,i}\}}=p^{u_{p}(k)}.
    \]
    As $G_{p}$ is a $p$-group, $p^{k}x=0$ if and only if $\ord(x)\mid p^{k}$. So for $k\geq1$,
    \[
       x\in G_{p}[p^{k}]\setminus G_{p}[p^{k-1}] \iff \ord(x)\mid p^{k} \,\text{ and }\, \ord(x)\nmid p^{k-1} \iff \ord(x)=p^{k}.
    \]
    Thus, $|\{x\in G_{p}:\ord(x)=p^{k}\}|=p^{u_{p}(k)}-p^{u_{p}(k-1)}$.\par
    Every cyclic $p$-subgroup of $G$ lies in $G_{p}$, and as the largest order of an element of $G_{p}$ is $p^{\lambda_{p,1}}$, the non-trivial ones have orders $p,p^{2},\cdots,p^{\lambda_{p,1}}$. Let $c_{p,k}(G)$ be the number of cyclic subgroups of $G_{p}$ of order $p^{k}$. Each of them has $\varphi(p^{k})=p^{k-1}(p-1)$ generators, and each element of order $p^{k}$ generates exactly one of them. So $c_{p,k}(G)\varphi(p^{k})=p^{u_{p}(k)}-p^{u_{p}(k-1)}$, and hence
    \[
       c_{p}(G)=\sum_{k=1}^{\lambda_{p,1}}c_{p,k}(G)=\sum_{k=1}^{\lambda_{p,1}}\frac{p^{u_{p}(k)}-p^{u_{p}(k-1)}}{p^{k-1}(p-1)}.
    \]
    Now if $G\cong C_{p^{a}}$, then $u_{p}(k)=\min\{k,a\}$, and so
    \[
       c_{p}(G)=\sum_{k=1}^{a}\frac{p^{k}-p^{k-1}}{p^{k-1}(p-1)}=\sum_{k=1}^{a}1=a.
    \]
    As every subgroup of $C_{p^{a}}$ is either $0$ or a non-trivial cyclic $p$-group, by Theorem~\ref{thm:abelian-classification}, $$|\N^\M(G)|=1+c_{p}(G)=a+1.$$ For any other non-trivial finite abelian group $G$, by the correspondence theorem, the subgroups of $G$ containing $\Soc(G)$ are in bijection with $\Sub(G/\Soc(G))$. Moreover, no non-trivial cyclic $p$-subgroup $C$ of $G$ contains $\Soc(G)$. Indeed, if $\Soc(G)\leq C$, then $\Soc(G)$ is a $p$-group, and so by Cauchy's theorem $G$ is a $p$-group. Then $\Soc(G)=G_{p}[p]\leq C$ is cyclic of order $p^{u_{p}(1)}=p^{r_{p}}$. So $r_{p}=1$ and $G\cong C_{p^{\lambda_{p,1}}}$, which is a contradiction. Also, $\Soc(G)\ne0$ as $G\ne0$. So the three types of $\M$-subgroups in Theorem~\ref{thm:abelian-classification} are pairwise disjoint, and hence
    \[
       |\N^\M(G)|=1+\sum_{p\mid|G|}c_{p}(G)+|\Sub(G/\Soc(G))|.
\qedhere    \]
\end{proof}
The $\M$-full abelian groups are those of rank at most $2$, in the following
sense.

\begin{corollary}\label{cor:abelian-Mfull}
Let $G$ be a non-trivial finite abelian group and put $r(G)=\sum_p r_p$.  Then $G$ is $\M$-full if
and only if $r(G)\leq2$.
\end{corollary}

\begin{proof}
    Following Definition~\ref{def:abelian-notation} and collecting the cyclic summands of all the $G_{p}$ together, we get
    \[
       G\cong\bigoplus_{p\mid|G|}\bigoplus_{i=1}^{r_{p}}C_{p^{\lambda_{p,i}}}.
    \]
    So $G$ is a direct sum of exactly $r(G)$ nontrivial cyclic $p$-groups. By the fundamental theorem of finite cyclic groups, each $C_{p^{\lambda_{p,i}}}$ is monolithic.\par
    If $r(G)=1$, then $G\cong C_{p^{a}}$ for some prime $p$ and some $a\geq1$. So every non-zero subgroup of $G$ is a cyclic $p$-group, and by Theorem~\ref{thm:abelian-classification} we get $\N^\M(G)=\N(G)$. So $G$ is $\M$-full. If $r(G)\geq2$, then as each $C_{p^{\lambda_{p,i}}}$ is monolithic, by Corollary~\ref{cor:product-monolithic},
    \[
       G\text{ is }\M\text{-full} \iff r(G)=2.
    \]
    Combining the two cases, $G$ is $\M$-full if and only if $r(G)\leq2$.
\end{proof}

\section{\texorpdfstring{$\M$}{M}-depth}
\label{sec:depth}

A chief series of $G$ is a maximal chain in $\N(G)$, and its length
$\elln(G)$ is an invariant of $G$. In this section we show that the same holds
for $\N^{\M}(G)$: all maximal $\M$-chains of a finite group have the same
length $\md(G)$, which we compute (Theorem~D). The key observation is that a
normal subgroup over a single minimal normal subgroup meets the socle only in
that minimal normal subgroup, so a maximal $\M$-chain must cross the socle in
one jump.

\begin{lemma}\label{lem:socle-complement}
Let $G$ be a non-trivial finite group, let $A\in\MinNorm(G)$, and put $S=\Soc(G)$.  There are minimal normal subgroups
\[
A=A_1,A_2,\ldots,A_{\sigma(G)}
\]
such that
\[
S=A_1\times A_2\times\cdots\times A_{\sigma(G)}.
\]
Consequently, if $T\trnl G$ contains exactly the atom $A$, then
\[
T\cap S=A,
\]
and every $G$-normal subgroup in the interval $[1,T]$ is an $\M$-subgroup.
\end{lemma}

\begin{proof}
    Let $N_{1}=A=A_{1}$. Then, $\lambda_{G}(N_{1})=1$. Suppose we have chosen $A=A_{1},A_{2},\cdots,A_{i}\in\MinNorm(G)$ such that $N_{i}=A_{1}\times\cdots\times A_{i}<\Soc(G)$ and $\lambda_{G}(N_{i})=i$. As $$\Soc(G)=\displaystyle{\bigvee\{N:N\in\MinNorm(G)\}},$$ so there exists $A_{i+1}\in\MinNorm(G)$ such that $A_{i+1}\nleq N_{i}$. Now $N_{i}\cap A_{i+1}\trnl G$, and as $A_{i+1}$ is minimal normal, this forces $N_{i}\cap A_{i+1}=1$. Hence, $[N_{i},A_{i+1}]\leq N_{i}\cap A_{i+1}=1$, and therefore $N_{i+1}=N_{i}A_{i+1}=N_{i}\times A_{i+1}$. By Lemma~\ref{lem:relative-length}, we get $$\lambda_{G}(N_{i+1})=\lambda_{G}(N_{i})+1=i+1,$$ since $\lambda_{G}(A_{i+1})=1$ and $\lambda_{G}(N_{i}\cap A_{i+1})=0$. As $G$ is finite, there exists a natural number $k$ such that $N_{k}=\Soc(G)$. So, $\sigma(G)=\lambda_{G}(N_{k})=k$. Thus, we have minimal normal subgroups $A=A_{1},A_{2},\cdots,A_{\sigma(G)}$ such that 
    \[
       S=A_{1}\times\cdots\times A_{\sigma(G)}.
    \]
    Now suppose $T\trnl G$ contains exactly the atom $A$. Then $A\leq T\cap S$. Put $C=A_{2}\times\cdots\times A_{\sigma(G)}$, so that $S=AC$ and $A\cap C=1$. Then by Lemma~\ref{lem:modular},
    \[
       T\cap S=AC\cap(T\cap S)=A(C\cap T\cap S)=A(C\cap T).
    \]
    If $C\cap T\ne1$, then as $C\cap T\trnl G$, there exists $B\in\MinNorm(G)$ with $B\leq C\cap T$. Then $B\leq T$, and $B\ne A$ as $A\cap C=1$. This contradicts the fact that $T$ contains exactly the atom $A$. Thus, $C\cap T=1$ and hence $T\cap S=A$.\par
    Furthermore, let $1\ne M\in[1,T]_{\N(G)}$. As $G$ is finite, there exists $K\in\MinNorm(G)$ such that $K\leq M\leq T$. As $T$ contains only the atom $A$, it follows that $K=A$, and hence $|\supp(M)|=1$. Thus, $M\in\N^\M(G)$ by Corollary~\ref{cor:finite-classification}. Also, $1\in\N^\M(G)$ trivially. Thus, every $G$-normal subgroup in the interval $[1,T]$ is an $\M$-subgroup.
\end{proof}
We can now prove Theorem~D.

\begin{theorem}\label{thm:depth-formula}
Let $G$ be a non-trivial finite group and put $S=\Soc(G)$.
\begin{enumerate}
\item If $G$ has one minimal normal subgroup, then
      \[
      \md(G)=\elln(G).
      \]
\item If $G$ has at least two minimal normal subgroups, then
      \[
      \md(G)=\elln(G)-\sigma(G)+2.
      \]
\end{enumerate}
\end{theorem}

\begin{proof}
    For $(1)$, observe that since $G$ has one minimal normal subgroup, every non-trivial normal subgroup of $G$ contains it, and so by Corollary~\ref{cor:finite-classification} it follows that $G$ is $\M$-full. Now if
    \[
       1=E_{0}<E_{1}<\cdots<E_{r}=G
    \]
    is a chief series of $G$, then $E_{i}\in\N(G)$ and as $G$ is $\M$-full so $E_{i}\in\N^\M(G)$ for all $i$. Conversely, every maximal $\M$-chain is a maximal chain of $\N(G)$ and hence is a chief series of $G$. Thus, $\elln(G)=\md(G)$.\par
    For $(2)$, consider a maximal $\M$-chain,
    \[
       1=E_{0}<E_{1}<\cdots<E_{r}=G.
    \]
    Let $E_{j}$ be the first element of this chain containing $S$, and put $T=E_{j-1}$. Then $S\nleq T$. Also $T\ne1$. Indeed, if $T=1$, then for any $A\in\MinNorm(G)$ we have $1<A<S\leq E_{1}$, as $G$ has at least two minimal normal subgroups, and $A\in\N^\M(G)$. This contradicts the maximality of the chain. As $T\in\N^\M(G)$, $T\ne1$ and $S\nleq T$, by Corollary~\ref{cor:finite-classification}, $|\supp(T)|=1$. Let $\supp(T)=\{A\}$. Then, $A\in\MinNorm(G)$ and by Lemma~\ref{lem:socle-complement}, we get $T\cap S=A$. Now $TS\in\N^\M(G)$ as $S\leq TS$. Thus we have
    \[
       T<TS\leq E_{j}.
    \]
    So maximality of the chain forces $TS=E_{j}$. Thus, our chain becomes
    \[
       1=E_{0}<E_{1}<\cdots<T<TS<\cdots<G
    \]
    By Lemma~\ref{lem:socle-complement}, $[1,T]_{\N(G)}=[1,T]_{\N^\M(G)}$, and every normal subgroup containing the socle is an $\M$-subgroup. So the parts of the chain from $1$ to $T$ and from $TS$ to $G$ cannot be refined in $\N(G)$, and as $\N(G)$ is modular by Lemma~\ref{lem:modular}, they have lengths $\lambda_{G}(T)$ and $\lambda_{G}(G)-\lambda_{G}(TS)$. Thus, we get
    \[
       r=\lambda_{G}(T)+1+\lambda_{G}(G)-\lambda_{G}(TS).
    \]
    Now by Lemma~\ref{lem:relative-length}, as $\lambda_{G}(T\cap S)=\lambda_{G}(A)=1$ and $\lambda_{G}(S)=\sigma(G)$, we get
    \[
       \lambda_{G}(TS)+1=\lambda_{G}(T)+\sigma(G).
    \]
    Substituting this value back, we get $r=\elln(G)-\sigma(G)+2$, since $\lambda_{G}(G)=\elln(G)$. So every maximal $\M$-chain has length $\elln(G)-\sigma(G)+2$, and as an $\M$-chain of maximum length is maximal, $\md(G)=\elln(G)-\sigma(G)+2$.
\end{proof}

The proof of Theorem~\ref{thm:depth-formula} describes all chains of maximal
length.

\begin{corollary}\label{cor:maximum-chain-structure}
Let $G$ be a finite group with at least two minimal normal subgroups and put
$S=\Soc(G)$.  Let $T\trnl G$ contain exactly one minimal normal subgroup.
Concatenate a maximal chain in $[1,T]_{\N(G)}$, the step
\[
T<TS,
\]
and a maximal chain in $[TS,G]_{\N(G)}$.  The resulting chain consists of
$\M$-subgroups and has length $\md(G)$.

Conversely, every $\M$-chain of length $\md(G)$ arises in this way, with $T$
its last nonessential non-trivial term.
\end{corollary}

\begin{proof}
    The argument follows precisely as in Theorem~\ref{thm:depth-formula}$(2)$. Indeed, let $\supp(T)=\{A\}$. As $G$ has at least two minimal normal subgroups, $S\nleq T$, and so $T<TS$. By Lemma~\ref{lem:socle-complement}, $T\cap S=A$ and the terms of the first chain are $\M$-subgroups, and the terms of the last chain contain $S$, so they are $\M$-subgroups too. As $\N(G)$ is modular, the resulting chain has length $\lambda_{G}(T)+1+\elln(G)-\lambda_{G}(TS)$, which is $$\elln(G)-\sigma(G)+2=\md(G)$$ by Lemma~\ref{lem:relative-length} and Theorem~\ref{thm:depth-formula}. Conversely, an $\M$-chain of length $\md(G)$ is maximal, and so by the proof of Theorem~\ref{thm:depth-formula}$(2)$ it has the above form with $T=E_{j-1}$. Finally, as every non-trivial normal subgroup of $G$ contains a minimal normal subgroup, a normal subgroup of $G$ is essential if and only if it contains $S$. So $T=E_{j-1}$ is the last non-essential non-trivial term of the chain.
\end{proof}

In particular, the $\M$-depth detects $\M$-fullness.

\begin{corollary}\label{cor:depth-Mfull}
For a non-trivial finite group $G$,
\[
\md(G)=\elln(G)
\quad\Longleftrightarrow\quad
G\text{ is $\M$-full}.
\]
\end{corollary}

\begin{proof}
    First suppose $\md(G)=\elln(G)$. Now we have two cases. Firstly, if $G$ contains a unique minimal normal subgroup. Then, every non-trivial $N\in\N(G)$ contains exactly one minimal normal subgroup, and hence $\N(G)=\N^\M(G)$ by Corollary~\ref{cor:finite-classification}. Secondly, if $G$ contains at least two minimal normal subgroups. Then by Theorem~\ref{thm:depth-formula}, we have $$\md(G)=\elln(G)-\sigma(G)+2.$$ As $\md(G)=\elln(G)$, it follows that $\sigma(G)=2$. By Theorem~\ref{thm:Mfull-classification}, $G$ is $\M$-full.\par
    Conversely, suppose $G$ is $\M$-full. If $G$ has a unique minimal normal subgroup then by Theorem~\ref{thm:depth-formula}, $\md(G)=\elln(G)$. Now suppose $G$ has at least two minimal normal subgroups. Then, $\md(G)=\elln(G)-\sigma(G)+2$. As $G$ is $\M$-full, by Theorem~\ref{thm:Mfull-classification}, $\sigma(G)\leq2$. But $\sigma(G)\ne1$ as $G$ has at least two minimal normal subgroups. Hence, $\sigma(G)=2$. Thus, $\md(G)=\elln(G)$.
\end{proof}

Since chief length and socle length are additive over direct products, we
obtain the product formula of Theorem~D.

\begin{corollary}\label{cor:product-depth}
Let $G_1,\ldots,G_t$ be non-trivial finite groups and let $t\geq2$.  Then
\[
\md(G_1\times\cdots\times G_t)
=\sum_{i=1}^t\bigl(\elln(G_i)-\sigma(G_i)\bigr)+2.
\]
In particular, if all $G_i$ are monolithic, then
\[
\md(G_1\times\cdots\times G_t)
=\sum_{i=1}^t\elln(G_i)-t+2.
\]
\end{corollary}

\begin{proof}
    Let $G=G_{1}\times\cdots\times G_{t}$, and for $X\leq G_{i}$ let $\widehat{X}$ be its coordinate copy in $G$. As $G_{j}^{*}$ centralizes $G_{i}^{*}$ for $j\ne i$, $\widehat{X}$ is normal in $G$ if and only if $X\trnl G_{i}$. So $X\mapsto\widehat{X}$ is a lattice isomorphism from $\N(G_{i})$ onto $[1,G_{i}^{*}]_{\N(G)}$, and hence
    \[
       \lambda_{G}(\widehat{X})=\lambda_{G_{i}}(X)\quad\text{for all }X\trnl G_{i}.
    \]
    In particular, $\lambda_{G}(G_{1}^{*})=\elln(G_{1})$. As $G/G_{1}^{*}\cong G_{2}\times\cdots\times G_{t}$, by Lemma~\ref{lem:relative-length},
    \[
       \elln(G)=\lambda_{G}(G_{1}^{*})+\elln(G/G_{1}^{*})=\elln(G_{1})+\elln(G_{2}\times\cdots\times G_{t}).
    \]
    Hence by induction on $t$, $\elln(G)=\displaystyle{\sum_{i=1}^{t}\elln(G_{i})}$. Also, by Theorem~\ref{thm:product-atoms}, $\Soc(G)=\widehat{\Soc(G_{1})}\cdots\widehat{\Soc(G_{t})}$, and each $\widehat{\Soc(G_{i})}$ intersects $\widehat{\Soc(G_{1})}\cdots\widehat{\Soc(G_{i-1})}$ trivially. So applying Lemma~\ref{lem:relative-length}(\ref{eq:modular-rank}) repeatedly, we get
    \[
       \sigma(G)=\sum_{i=1}^{t}\lambda_{G}(\widehat{\Soc(G_{i})})=\sum_{i=1}^{t}\sigma(G_{i})\geq t\geq2
    \]
    as $\sigma(G_{i})\geq1$ for each $i$. Thus $G$ has at least two minimal normal subgroups. Hence by Theorem~\ref{thm:depth-formula}, we have
    \[
       \md(G)=\elln(G)-\sigma(G)+2=\sum_{i=1}^{t}(\elln(G_{i})-\sigma(G_{i}))+2.
    \]
    Now if each $G_{i}$ is monolithic, then $\sigma(G_{i})=1$ for all $i$. Hence $\md(G)=\displaystyle{\sum_{i=1}^{t}\elln(G_{i})}-t+2$.
\end{proof}

For crown-based powers the $\M$-depth no longer depends on the crown
multiplicity, once this is at least $2$.

\begin{corollary}\label{cor:crown-depth}
Let $L$ be finite monolithic primitive with socle $A$, and let $L_k$ be its
crown-based power.  Then
\[
\md(L_k)=
\begin{cases}
\elln(L/A)+1,&k=1,\\
\elln(L/A)+2,&k\geq2.
\end{cases}
\]
\end{corollary}

\begin{proof}
    Let $V=\Soc(L_{k})$. Now if $k=1$, then $L_{1}=L$ and $V=A$. By Corollary~\ref{cor:crown-M-full}, $L_{1}$ is $\M$-full. Hence by Corollary~\ref{cor:depth-Mfull}, we get $\md(L_{1})=\elln(L_{1})$. By Lemma~\ref{lem:relative-length}, we have $\elln(L)=\lambda_{L}(A)+\elln(L/A)$. As $A$ is minimal normal in this case, we have $\lambda_{L}(A)=1$. Thus, $\md(L_{1})=\elln(L/A)+1$.\par
    For $k\geq2$, by Proposition~\ref{prop:crown-M-chains}, every maximal $\M$-chain of $L_{k}$ has length $\elln(L/A)+2$. Thus, $\md(L_{k})=\elln(L/A)+2$.
\end{proof}

We now evaluate $\md$ for several classical families. For nilpotent groups
all chief factors have prime order, so $\elln(G)=\Omega(|G|)$.

\begin{corollary}\label{cor:nilpotent-depth}
Let $G$ be a non-trivial finite nilpotent group.  Then
\[
\md(G)=
\begin{cases}
\Omega(|G|),&\rho(G)=1,\\
\Omega(|G|)-\rho(G)+2,&\rho(G)\geq2.
\end{cases}
\]
\end{corollary}

\begin{proof}
    Let
    \[
       1=G_{0}<G_{1}<\cdots<G_{r}=G
    \]
    be a chief series of $G$. So $G_{i}/G_{i-1}$ is a minimal normal subgroup of $G/G_{i-1}$. As $G$ is finite nilpotent, $G/G_{i-1}$ is also finite nilpotent. By the proof of Corollary~\ref{cor:nilpotent-Mfull} we get $[G_{i}:G_{i-1}]=p_{i}$ for some prime $p_{i}$, for all $i$. Hence, $|G|=\displaystyle{\prod_{i=1}^{r}[G_{i}:G_{i-1}]}=p_{1}p_{2}\cdots p_{r}$, and so $\elln(G)=r=\Omega(|G|)$.\par
    Now suppose $\rho(G)=1$. Then by Corollary~\ref{cor:nilpotent-Mfull} we have that $G$ is $\M$-full. So by Corollary~\ref{cor:depth-Mfull}, $\md(G)=\elln(G)=\Omega(|G|)$. Again, by the proof of Corollary~\ref{cor:nilpotent-Mfull} we have $\sigma(G)=\rho(G)$. Thus for $\rho(G)\geq2$ we have $\sigma(G)\geq2$. Hence $G$ has at least two minimal normal subgroups. So by Theorem~\ref{thm:depth-formula}, \[\md(G)=\elln(G)-\sigma(G)+2=\Omega(|G|)-\rho(G)+2.\qedhere\]
\end{proof}

For abelian groups $\rho(G)$ is the rank $r(G)$, and we obtain the following.

\begin{theorem}\label{thm:abelian-depth}
Let $G$ be a non-trivial finite abelian group.
\begin{enumerate}
\item If $G$ is a cyclic $p$-group, then
      \(\md(G)=\Omega(|G|)\).
\item Otherwise,
      \[
      \md(G)=\Omega(|G|)-r(G)+2.
      \]
\end{enumerate}
\end{theorem}

\begin{proof}
    As $G$ is abelian, it is nilpotent and its Sylow $p$-subgroup $G_{p}$ satisfies $Z(G_{p})=G_{p}$. Also, as each cyclic summand of $G_{p}$ has a unique subgroup of order $p$, $G_{p}[p]\cong C_{p}^{r_{p}}$. Hence,
    \[
       \rho(G)=\sum_{p\mid|G|}\dim_{\F_{p}}G_{p}[p]=\sum_{p\mid|G|}r_{p}=r(G).
    \]
    Now $r(G)=1$ if and only if $G$ is a cyclic $p$-group. So the result follows from Corollary~\ref{cor:nilpotent-depth}.
\end{proof}

For cyclic groups the formula involves only the prime factorisation of the
order.

\begin{corollary}\label{cor:cyclic-depth}
Let $n\geq2$ and write
$n=p_1^{a_1}\cdots p_s^{a_s}$ with distinct primes $p_1,\ldots,p_s$.  Then
\[
\md(C_n)=
\begin{cases}
 a_1,&s=1,\\
 \Omega(n)-\omega(n)+2,&s\geq2.
\end{cases}
\]
\end{corollary}

\begin{proof}
    First let $s=1$. Then $C_{n}=C_{p_{1}^{a_{1}}}$ is a cyclic $p_{1}$-group. So by Theorem~\ref{thm:abelian-depth}, $\md(C_{n})=\Omega(p_{1}^{a_{1}})=a_{1}$.
    Now let $s\geq2$. By the Chinese Remainder Theorem, $$C_{n}\cong C_{p_{1}^{a_{1}}}\oplus\cdots\oplus C_{p_{s}^{a_{s}}}.$$ So the Sylow $p_{i}$-subgroup of $C_{n}$ is cyclic, i.e., $r_{p_{i}}=1$ for all $i$, and hence $r(C_{n})=s=\omega(n)\geq2$. In particular, $C_{n}$ is not a cyclic $p$-group, and so by Theorem~\ref{thm:abelian-depth}, $\md(C_{n})=\Omega(n)-\omega(n)+2$.
\end{proof}

Non-nilpotent examples can be handled in the same way when the normal
structure is known. For dihedral groups of order $2n$ with $n$ odd, all
non-trivial proper normal subgroups lie in the rotation subgroup.

\begin{proposition}\label{prop:dihedral-depth}
Let $D_{2n}$ be the dihedral group of Definition~\ref{def:arithmetic-notation}, with odd $n\geq3$.  Then
\[
\md(D_{2n})=\md(C_n)+1.
\]
\end{proposition}

\begin{proof}
    Recall from Definition~\ref{def:arithmetic-notation} that
    \[
       D_{2n}=\langle r,s : r^{n}=s^{2}=1,\ srs=r^{-1}\rangle.
    \]
    Suppose $R=\langle r\rangle$. Then $R\cong C_{n}$, and every subgroup of $R$ is $\langle r^{d}\rangle$ for some $d\mid n$. As
    \[
       r\langle r^{d}\rangle r^{-1}=\langle r^{d}\rangle \quad\text{and}\quad s\langle r^{d}\rangle s^{-1}=\langle r^{-d}\rangle=\langle r^{d}\rangle,
    \]
    so every subgroup of $R$ is normal in $D_{2n}$. Conversely, let $N\trnl D_{2n}$. If $N$ contains no reflection, then $N\leq R$. Otherwise, let $sr^{i}\in N$ for some $i$. As $rs=sr^{-1}$, we get $r(sr^{i})r^{-1}=sr^{i-2}\in N$, and so $(sr^{i-2})(sr^{i})=r^{2}\in N$. As $n\geq3$ is odd, $\ord(r^{2})=n$, and hence $R=\langle r^{2}\rangle\leq N$. Then $s=(sr^{i})r^{-i}\in N$, and thus $N=D_{2n}$. Hence,
    \[
       \N(D_{2n})=\Sub(C_{n})\sqcup \{D_{2n}\}.
    \]
    This also implies that the minimal normal subgroups of $D_{2n}$ coincide with those of $C_{n}$, and so $\Soc(D_{2n})=\Soc(C_{n})$. So for $E\in\Sub(C_{n})$, we have $\supp_{D_{2n}}(E)=\supp_{C_{n}}(E)$, and thus by Corollary~\ref{cor:finite-classification}, $E\in\N^\M(D_{2n})\iff E\in\N^\M(C_{n})$. Also $\Soc(D_{2n})\leq D_{2n}$, and so $D_{2n}\in\N^\M(D_{2n})$. Thus we get,
    \[
       \N^\M(D_{2n})=\N^\M(C_{n})\sqcup \{D_{2n}\}.
    \]
    So every maximal $\M$-chain of $D_{2n}$ looks like
    \[
       1=G_{0}<G_{1}<\cdots<G_{r}=C_{n}<D_{2n},
    \]
    where $1=G_{0}<G_{1}<\cdots<G_{r}=C_{n}$ is a maximal $\M$-chain of $C_{n}$. Thus, we get $\md(D_{2n})=\md(C_{n})+1$.
\end{proof}

Finally, in powers of simple groups the $\M$-depth is as small as it can be.

\begin{corollary}\label{cor:simple-power-depth}
If $T$ is a non-abelian finite simple group and $k\geq1$, then
\[
\md(T^k)=
\begin{cases}1,&k=1,\\2,&k\geq2.\end{cases}
\]
The same formula holds for $C_p^k$.
\end{corollary}

\begin{proof}
    We have by Proposition~\ref{prop:simple-powers} $(1)$,
    \[
       \MinNorm(T^{k})=\{T_{1},\cdots,T_{k}\},
    \]
    where $T_{i}$ is the $i$th coordinate factor. Moreover,
    \[
       \N^\M(T^{k})=\begin{cases}
           \{1,T\}, & k=1, \\
           \{1,T_{1},\cdots,T_{k},T^{k}\}, & k\geq2.
       \end{cases}
    \]
    Thus, when $k=1$, $\md(T^{k})=1$. When $k\geq2$, the $\M$-subgroups other than $1$ and $T^{k}$ are $T_{1},\cdots,T_{k}$, which are pairwise incomparable. So an $\M$-chain contains at most one $T_{i}$, and as $1<T_{1}<T^{k}$ is an $\M$-chain, we have $\md(T^{k})=2$.\par
    Again by Proposition~\ref{prop:simple-powers} $(2)$, we have
    \[
       \N^\M(C_{p}^{k})=\{0\}\cup \{U\leq C_{p}^{k}:\dim_{\F_{p}}(U)=1\}\cup \{C_{p}^{k}\}.
    \]
    So if $k=1$, then $\md(C_{p})=\elln(C_{p})=1$. If $k\geq2$, then in the same way, the one-dimensional subspaces are pairwise incomparable, and $0<U<C_{p}^{k}$ for every one-dimensional subspace $U$. Hence, $\md(C_{p}^{k})=2$. Thus, the same formula, for $T^{k}$, holds for $C_{p}^{k}$ too.
\end{proof}

\textbf{Concluding remarks.}
\label{sec:remarks}
The results above show that, for finite groups, the $\M$-condition is a
condition on the socle, and that the defect $\N(G)\setminus\N^{\M}(G)$ is
measured by the socle length $\sigma(G)$ and by the geometry of
$\MinNorm(G)$. We mention two directions suggested by this picture.

First, Theorems~\ref{thm:dichotomy} and~\ref{thm:Mfull-udim} hold for arbitrary groups, and Theorem~\ref{thm:atomic-classification} for all normally atomic groups, whereas
Example~\ref{ex:Z3} shows that the closure-theoretic part of Theorem~A fails
without atomicity. It would be interesting to describe $\N^{\M}(G)$ for
classes of infinite groups that are not normally atomic, such as
finitely generated nilpotent or residually finite groups, where essential
normal subgroups are abundant but minimal normal subgroups may not exist.

Second, Corollary~\ref{cor:arbitrary-crown} only uses the quotients
$G/R_G(B)$. Since crowns are defined for profinite groups as well, one may ask
for the analogue of Theorems~B and~D for the lattice of closed normal
subgroups of a profinite group, where $\sigma$ and $\elln$ may be infinite.

\bigskip
\bigskip

\noindent\textbf{Declarations.}
\smallskip

\noindent\textbf{Competing interests.}
The authors declare that they have no competing financial or non-financial interests that are directly or indirectly related to this work.

\smallskip
\noindent\textbf{Data availability.}
No datasets were generated or analyzed during the present study.

\end{document}